\documentclass[pdflatex,sn-mathphys-num]{sn-jnl}
\usepackage{mathrsfs}

\usepackage{subcaption}
\usepackage{graphicx}%
\usepackage{multirow}%
\usepackage{amsmath,amssymb,amsfonts}%
\usepackage{amsthm}%
\usepackage{mathrsfs}%
\usepackage[title]{appendix}%
\usepackage{xcolor}%
\usepackage{textcomp}%
\usepackage{manyfoot}%
\usepackage{booktabs}%
\usepackage{array}%
\usepackage{placeins}%
\usepackage{algorithm}%
\usepackage{algorithmicx}%
\usepackage{algpseudocode}%
\usepackage{listings}%

\usepackage{amsthm}
\newtheorem{lemma}{Lemma}
\newtheorem{corollary}{Corollary}
\usepackage[T1]{fontenc}
\usepackage{graphicx}
\usepackage{xcolor}
  \usepackage{amsmath, amssymb, amsthm}
\usepackage{tikz} 
\usepackage{pgfplots} 
\pgfplotsset{compat=1.18}
\newcommand{\Adag}[1]{#1^{\dagger}}
\newcommand{\Perp}[1]{#1^{\perp}}

\newcommand{\expM}[1]{e^{#1}}
\newcommand{\trace}[1]{\textnormal{Tr}\left({#1}\right)}
\newcommand{\grad}[1]{\nabla{#1}}
\newcommand{\Rn}{\mathbb{R}^n}
\newcommand{\divg}[1]{{\rm div_{#1}}}

\usepackage{ifthen}

\theoremstyle{thmstyleone}%
\newtheorem{theorem}{Theorem}
\newtheorem{proposition}[theorem]{Proposition}%

\theoremstyle{thmstyletwo}%
\newtheorem{remark}{Remark}%

\theoremstyle{thmstylethree}%

\begin{document}

\title[Wasserstein KL Divergence for Gaussian Measures on \(\mathbb{R}^n\)]{Information Geometry for Wasserstein KL Divergence of Gaussian Measures on \(\mathbb{R}^n\)}


\author*[1]{\fnm{Adwait} \sur{Datar}}\email{adwait.datar@tuhh.de}

\author[1,2,3]{\fnm{Nihat} \sur{Ay}}\email{nihat.ay@tuhh.de}

\affil*[1]{\orgdiv{Institute for Data Science Foundations}, \orgname{Hamburg University of Technology}, \orgaddress{\city{Hamburg}, \postcode{21073}, \country{Germany}}}

\affil[2]{\orgname{Santa Fe Institute}, \orgaddress{ \city{Santa Fe}, \postcode{87501}, \state{New Mexico}, \country{USA}}}

\affil[3]{\orgname{Leipzig University}, \orgaddress{ \city{Leipzig}, \postcode{04109}, \country{Germany}}}



\abstract{
We study the Wasserstein Kullback--Leibler divergence (WKL divergence) on the manifold of nondegenerate Gaussian measures over $\mathbb R^n$. 
In the canonical-divergence construction, the Fisher--Rao metric recovers forward KL along intrinsic mixture geodesics and reverse KL along geodesics of the conjugate exponential connection. 
Replacing the Fisher--Rao metric by the Otto metric and following the latter route produces the $e_1$-connection underlying WKL divergence. 
We establish its geodesic completeness, classify its forward limits, prove that every ordered pair is joined by a unique $e_1$-connector generated by a quadratic potential, and derive an explicit WKL divergence formula with separate mean and covariance contributions. 
WKL divergence is nonnegative and separating. 
At equal covariances WKL divergence equals one half of the squared Euclidean mean distance. 
Finally, WKL divergence extends finitely and continuously to singular targets from a nondegenerate source, but diverges when the target remains nondegenerate and the source covariance becomes singular. 
Path-dependent joint limits at Dirac pairs preclude a continuous extension to the full positive-semidefinite covariance product, although a lower-semicontinuous extended-real extension exists.
}




\maketitle

\section{Introduction}\label{sec:introduction}
Classical information geometry equips statistical models with intrinsic objects such as the Fisher--Rao metric and the Kullback--Leibler divergence (KL divergence), without requiring a metric on the sample space \cite{AN06,ay2017information}. 
In this article the sample space is $\mathbb R^n$, endowed with its Euclidean geometry, and we ask for an information divergence that reflects this ground geometry. The classical KL divergence does not do so: if $x\neq y$, then the KL divergence between $\delta_x$ and $\delta_y$ is infinite, independently of the Euclidean distance $\lVert x-y\rVert$.

Our starting point is Ay's canonical-divergence construction \cite{ay2015novel}; see also \cite{felice2021towards}. For a Riemannian metric $g$, an affine connection $\nabla$, and its $\nabla$-geodesic $\gamma$ from $\mu$ to $\nu$, the canonical divergence is the weighted energy
\begin{align*}
D_g^{(\nabla)}(\mu\Vert\nu)
=\int_0^1 t\,\lVert\dot\gamma(t)\rVert_{\gamma(t),g}^2\,dt.
\end{align*}
With the Fisher--Rao metric, the intrinsic Gaussian mixture\footnote{Here ``mixture'' is intrinsic to the Gaussian family, not a literal density mixture.} geodesic yields the forward divergence $D_{\rm KL}(\mu\Vert\nu)$. Its Fisher--Rao conjugate is the usual exponential connection, denoted $e_0$ here to distinguish it from $e_1$; the $e_0$-geodesic yields the reverse divergence $D_{\rm KL}(\nu\Vert\mu)$. Both calculations are included for completeness in Appendices~\ref{sec:app_mixture_canonical} and~\ref{sec:app_reverse_KL}.

We follow the same conjugate-connection route after replacing the Fisher--Rao metric by the Otto metric. 
Computing the canonical divergence for the intrinsic mixture geodesic with the norm induced by the Otto metric gives an exact integral representation, but not an explicit closed form in general; see Appendix~\ref{sec:app_mixture_canonical}. 
Instead, conjugating the mixture connection with respect to the Otto metric gives the $e_1$-connection \cite{ay2024information,ay2025correction}. 
Its Gaussian geodesics can be computed explicitly, and the corresponding canonical divergence yields a closed endpoint formula; see Appendix~\ref{sec:WKL_divergence_formula}. 
This is the Otto analogue of the reverse-KL route above. 
Since the Otto metric induces the $2$-Wasserstein distance as its Riemannian distance \cite{Otto,KZ22}, we call the resulting functional the \emph{Wasserstein KL divergence}, abbreviated WKL divergence, following our preliminary work \cite{datar2026wkl}.

The Gaussian WKL divergence construction and its central closed formula first appeared in our Geometric Science of Information conference paper \cite{datar2026wkl}. 
That paper derived the closed formula for WKL divergence, including commuting and univariate specializations; compared the equal-covariance case with reverse KL; and presented preliminary equal-covariance Dirac-limit observations and numerical illustrations.
The conference paper applied the analytic formulation in \cite{ay2024information,ay2025correction}, which assumes a compact sample-space manifold, to Gaussian manifold over the non-compact sample space $\Rn$ without supplying the corresponding noncompact justification.
Relative to that conference paper, the present journal article supplies an information geometric foundation for the WKL divergence on the Gaussian manifold over the noncompact sample space $\Rn$.  
The main geometric advances here place the previously derived $e_1$ endpoint connector within a rigorous framework, prove geodesic completeness of the Gaussian $e_1$-connection, give the forward limits of its geodesics. 
For completeness, the appendices describe the intrinsic Gaussian mixture connection, verify its Fisher--Rao/$e_0$ and Otto/$e_1$ metric-conjugacy relations, and recover the classical forward-KL identity; they also derive the Otto/mixture integral representation and its closed univariate specialization.
The present article also strengthens the analytic study of WKL divergence. A general boundary theorem replaces the preliminary Dirac observation: it characterizes both one-sided singular-covariance limits, analyzes joint Dirac limits, constructs a lower-semicontinuous extension, and proves that no continuous extension exists on the full product of positive-semidefinite covariance cones.

Subsequent work, not part of the present article, places WKL divergence in a broader family of state-space-aware KL divergences, extends its Gaussian formula to elliptic families, introduces Kalman--Wasserstein and other metric-dependent variants, and applies these divergences to low-noise stochastic control \cite{stein2026wklcontrol}. That work takes the Gaussian WKL divergence formula as its starting point; the present article is concerned with the intrinsic Gaussian geometry, global connector, local expansion, and singular boundary behavior of WKL divergence itself.

The remainder of the paper is organized as follows. Section~\ref{sec:gaussian_geometry} develops the Gaussian manifold, the Fisher--Rao and Otto metrics, the $e_0$- and $e_1$-connections, and their global properties. Section~\ref{sec:gaussian_wkl} derives the explicit WKL divergence formula and studies its transport interpretation, local behavior, comparisons, univariate specialization, and singular-covariance limits. Appendix~\ref{sec:app_mixture_duality} develops the intrinsic Gaussian mixture connection, its metric conjugates, and the corresponding canonical divergences; Appendix~\ref{sec:app_exponential_canonical} derives the reverse-KL divergence and WKL divergence identities along the $e_0$- and $e_1$-geodesics.

\section{Gaussian Information Geometry: Metrics, Connections, and Canonical Divergences}\label{sec:gaussian_geometry}
In this section, we introduce the information geometric constructions leading to the Wasserstein KL divergence. 
We follow the development from \cite{ay2024information,ay2025correction} and verify the required constructions directly in the Gaussian setting; see those references for further details.
On one hand, the results in the present paper can be seen as a special case for instantiation of the results from \cite{ay2024information,ay2025correction} to the case of finite
dimensional Gaussian measures. 
On the other hand, these earlier results assume a compact sample manifold structure which excludes Gaussian measures defined on $\Rn$.
The following calculations can thus be also seen as an extension of earlier results in the Gaussian setting, without assuming compactness of the underlying manifold.

\noindent\textbf{Notation and conventions.}
Throughout, $I$ denotes the identity matrix of the size determined by context, $A^\top$ the transpose of $A$, and $\trace{A}$ its trace. We write $\lVert\cdot\rVert$ for the Euclidean norm and $\lVert\cdot\rVert_F$ for the Frobenius norm. For a symmetric matrix $A$, the relations $A\succ0$ and $A\succeq0$ mean positive definiteness and positive semidefiniteness, respectively. The matrix exponential is denoted by $e^A$, and $\grad{f}$ denotes the Euclidean gradient of $f$.

\subsection{The Gaussian Manifold and Its Tangent and Cotangent Spaces}\label{sec:tgt_cotgt_spaces}

Let $\lambda$ denote the Lebesgue measure on $\mathbb{R}^n$, and let $C^\infty(\mathbb{R}^n)
$ denote the space of smooth real-valued functions on $\mathbb{R}^n$. 
Let the space of finite signed measures on $\mathbb{R}^n$ be denoted by $\mathcal{M}(\mathbb{R}^n)$ which is a Banach space endowed with the total variation norm.
Let
\begin{align*}
    \mathcal{S}^\infty
    :=
    \left\{
    f\lambda \in \mathcal{M}(\mathbb{R}^n)
    \;\middle|\;
    f \in C^\infty(\mathbb{R}^n),
    \quad
    \int_{\mathbb{R}^n} |f(x)|\, d\lambda(x) < \infty
    \right\}.
\end{align*}
be the space of finite signed measures on $\mathbb{R}^n$ admitting a smooth density with respect to Lebesgue measure.
We further define the linear subspace
\begin{align*}
    \mathcal{S}_0^\infty
    :=
    \left\{
    \mu \in \mathcal{S}^\infty
    \;\middle|\;
    \mu(\mathbb{R}^n)=0
    \right\}.
\end{align*}
Finally define
\begin{align*}
    \mathcal{P}_+^\infty
    :=
    \left\{
    \mu \in \mathcal{S}^\infty
    \;\middle|\;
    \mu(\mathbb{R}^n)=1,
    \quad
    \mu > 0
    \right\},
\end{align*}
to be the space of smooth positive probability measures on $\mathbb{R}^n$.

Let
\[
S_n:=\{A\in\mathbb{R}^{n\times n}\mid A=A^\top\},
\qquad
S_n^+:=\{A\in S_n\mid A\succ0\}.
\]
Thus, $S_n^+$ is an open cone in the vector space $S_n$. Define the parametrization
\begin{align*}
    \varphi :
    \mathbb{R}^n \times S_n^+
    \longrightarrow
    \mathcal{P}_+^\infty, \quad 
    (m,\Sigma)
    \longmapsto
    p(\;\cdot\;;m,\Sigma)\lambda,
\end{align*}
where
\begin{align*}
    p(x;m,\Sigma)
    :=
    \frac{1}{(2\pi)^{n/2}\det(\Sigma)^{1/2}}
    \exp\left(
    -\frac{1}{2}(x-m)^\top \Sigma^{-1}(x-m)
    \right)
\end{align*}
is the multivariate Gaussian density with mean $m \in \mathbb{R}^n$ and covariance matrix $\Sigma \in S_n^+$.
The parametrization $\varphi$ is injective. We equip its image
\[
\mathcal{G}:=\varphi(\mathbb{R}^n\times S_n^+)\subset\mathcal{P}_+^\infty
\]
with the smooth manifold structure transported from $\mathbb{R}^n\times S_n^+$, so that $\varphi$ is a global chart. We refer to $\mathcal{G}$ as the manifold of nondegenerate Gaussian measures and write $\mathcal{N}(m,\Sigma)$ for the Gaussian measure with parameters $(m,\Sigma)$.

Let us now introduce the tangent space at a point $\mu \in \mathcal{G}$. To this end, consider a smooth curve $(-\varepsilon,\varepsilon)\ni t \mapsto (m_t,\Sigma_t) \in \mathbb{R}^n \times S_n^+$ and the $\varphi$ image of this curve in $\mathcal{G}$ given by
\begin{align*}
   (-\varepsilon,\varepsilon)\ni t \mapsto \mu_t := \varphi(m_t,\Sigma_t) = p(\;\cdot\; ;m_t,\Sigma_t) \lambda \in \mathcal{G}.
\end{align*}
First note that
\begin{align*}
   &\frac{d}{dt} p(x;m_t,\Sigma_t) \\
   &= p(x;m_t,\Sigma_t) \frac{d}{dt}\left(\ln p(x;m_t,\Sigma_t) \right)\\
   &=p(x;m_t,\Sigma_t)\Bigl(\frac{1}{2}(x-m_t)^\top \Sigma_t^{-1}\dot{\Sigma}_t\Sigma_t^{-1}(x-m_t)\\
   &\hspace{3.5cm}+\dot{m}_t^\top\Sigma_t^{-1}(x-m_t)-\frac{1}{2}\trace{\Sigma_t^{-1}\dot{\Sigma}_t}\Bigr)\\
   &=p(x;m_t,\Sigma_t)\cdot g(x;m_t,\Sigma_t,\dot{m}_t,\dot{\Sigma}_t),
\end{align*}
where the quadratic function $g$ is defined for $(m,\Sigma)\in \Rn \times S_n^{+}$ and $(m_v,\Sigma_v)\in \Rn \times S_n$
\begin{align*}
g(x;m,\Sigma,m_v,\Sigma_v)
&=\frac{1}{2}(x-m)^\top\Sigma^{-1}\Sigma_v\Sigma^{-1}(x-m)\\
&\quad+m_v^\top\Sigma^{-1}(x-m)-\frac{1}{2}\trace{\Sigma^{-1}\Sigma_v}.
\end{align*}
\begin{samepage}
Let us compute the velocity of the curve $\mu_t$ at $t=0$ as
\begin{align*}
\frac{d}{dt}\Big|_{t=0} \mu_t
= \lim_{t\rightarrow 0}\frac{1}{t} \left(\mu_t -\mu_0\right) &= \lim_{t\rightarrow 0}\frac{1}{t}\left(\big(p(\;\cdot\; ;m_t,\Sigma_t) - p(\;\cdot\; ; m_0,\Sigma_0) \big)\lambda\right)\\
&=
\left(\frac{d}{dt}\Big|_{t=0} p(\;\cdot\; ;m_t,\Sigma_t)\right)\lambda \\
&=g(\;\cdot\; ; m_0,\Sigma_0,\dot{m}_0,\dot{\Sigma}_0) \mu_0,
\end{align*}
\end{samepage}
where the derivative is taken in the total-variation norm of the ambient space $\mathcal M(\mathbb R^n)$, equivalently in $L^1(\lambda)$ for the densities. More generally, on each compact parameter neighborhood, all first parameter derivatives of the Gaussian density are uniformly dominated by an integrable function of the form $C(1+\lVert x\rVert^k)e^{-c\lVert x\rVert^2}$ and vary pointwise continuously. Dominated convergence therefore shows that $\varphi$ is continuously differentiable as an $L^1(\lambda)$-valued, equivalently $\mathcal M(\mathbb R^n)$-valued, map.
Moreover,
\begin{align*}
\left(\left.\frac{d}{dt}\right|_{t=0}\mu_t\right)(\Rn)
&=\int_{\Rn}g(x;m_0,\Sigma_0,\dot m_0,\dot\Sigma_0)\,d\mu_0(x)\\
&=\mathbb{E}_{\mu_0}\!\left[g(\,\cdot\,;m_0,\Sigma_0,\dot m_0,\dot\Sigma_0)\right]=0,
\end{align*}
so $\left.\frac{d}{dt}\right|_{t=0}\mu_t\in\mathcal{S}_0^\infty$. As usual, tangent vectors are equivalence classes of curves with the same velocity, and we identify them with these velocities. The Fr\'echet differential (and hence also the G\^ateaux differential) is
\[
d\varphi_{(m,\Sigma)}(m_v,\Sigma_v)
=g(\,\cdot\,;m,\Sigma,m_v,\Sigma_v)\mu.
\]
Since $\varphi$ is a global chart, its differential identifies $\Rn\times S_n$ with the tangent space; explicitly,
\[
T_\mu\mathcal{G}
=\operatorname{im}d\varphi_{(m,\Sigma)}
=\left\{g(\,\cdot\,;m,\Sigma,m_v,\Sigma_v)\mu\ \middle|\ (m_v,\Sigma_v)\in\Rn\times S_n\right\}.
\]
This map is injective, since a quadratic polynomial $g(\,\cdot\,;m,\Sigma,m_v,\Sigma_v)$ that vanishes $\mu$-almost surely vanishes identically, which forces $m_v=0$ and $\Sigma_v=0$. Thus $\dim T_\mu\mathcal{G}=n+n(n+1)/2$.
Motivated by this, let
\[
\begin{aligned}
\mathcal Q&:=\left\{x\mapsto\frac12x^\top Ax+b^\top x+c:
A\in S_n,\ b\in\mathbb R^n,\ c\in\mathbb R\right\},\\
\mathcal Q_\mu^0&:=\{f\in\mathcal Q:\mathbb E_\mu f=0\}.
\end{aligned}
\]
Write $\mathcal Q/\mathbb R$ for the quotient by constant functions and
$f+\mathbb R$ for the class of $f$.
The displayed form of $g$ and the injectivity above give
\begin{equation}\label{eq:gaussian_tangent_scores}
T_\mu\mathcal G=\{f\mu:f\in\mathcal Q_\mu^0\}.
\end{equation}
Every tangent vector is a signed measure $a=h\mu$ with $h\in\mathcal Q_\mu^0$. Since Gaussian measures have finite moments of every order, quadratic functions are integrable against every such tangent measure. Moreover, tangent measures have total mass zero, so for every $c\in\mathbb R$,
\[
\int_{\mathbb R^n}(f+c)\,da
=\int_{\mathbb R^n}f\,da+c\,a(\mathbb R^n)
=\int_{\mathbb R^n}f\,da.
\]
Thus integration depends only on the equivalence class modulo constants and defines the well-defined pairing
\begin{equation}\label{eq:quadratic_cotangent_pairing}
(f+\mathbb R,a)\longmapsto\int_{\mathbb R^n}f\,da,
\qquad
(f+\mathbb R,a)\in\mathcal Q/\mathbb R\times T_\mu\mathcal G.
\end{equation}
Every class $f+\mathbb R$ has the centered representative $f-\mathbb E_\mu f\in\mathcal Q_\mu^0$. It is unique: two centered representatives in the same class differ by a constant whose $\mu$-expectation is zero. Consequently, the canonical map
\[
\mathcal Q/\mathbb R\longrightarrow\mathcal Q_\mu^0,
\qquad
f+\mathbb R\longmapsto f-\mathbb E_\mu f,
\]
is an isomorphism.

\begin{proposition}
Let $\mu=\mathcal{N}(m,\Sigma)\in\mathcal{G}$.
Every linear functional $\ell\in T^*_\mu\mathcal{G}$ can be represented uniquely in the form
\[
\ell(a)=\int_{\mathbb{R}^n} f\,da,
\qquad a\in T_\mu\mathcal{G},
\]
for a unique equivalence class $f+\mathbb{R}\in \mathcal{Q}/\mathbb{R}$.
In particular,
\[
T_\mu^*\mathcal{G}\cong \mathcal{Q}/\mathbb{R}.
\]
\end{proposition}

\begin{proof}
Under \eqref{eq:gaussian_tangent_scores}, the pairing \eqref{eq:quadratic_cotangent_pairing} is the $L^2(\mu)$ product
\[
(f,h)\longmapsto\mathbb E_\mu[fh],
\qquad f,h\in\mathcal Q_\mu^0.
\]
It is positive definite because $\mu$ has full support. Hence its finite-dimensional Riesz map is an isomorphism from $\mathcal Q_\mu^0$ onto $(T_\mu\mathcal G)^*$. Composing with $\mathcal Q/\mathbb R\cong\mathcal Q_\mu^0$ proves existence and uniqueness of the representing class. Continuity is automatic in finite dimension.
\end{proof}
\subsection{Affine Vector Fields and Their Induced Flows}
Let $\mathcal{T}(\mathbb{R}^n)$ denote the space of smooth vector fields on $\mathbb{R}^n$, namely
\begin{align*}
    \mathcal{T}(\mathbb{R}^n)
=
\left\{
X : \mathbb{R}^n \to \mathbb{R}^n
\;\middle|\;
X \in C^\infty(\mathbb{R}^n,\mathbb{R}^n)
\right\}.
\end{align*}
In this work, we are particularly interested in the class of affine vector fields. A vector field $X \in \mathcal{T}(\mathbb{R}^n)$ is called \emph{affine} if it is of the form $X(x) = Ax + b$, where $A \in \mathbb{R}^{n \times n}$ and $ b \in \mathbb{R}^n$.
The set of affine vector fields is denoted by $\mathcal{T}_a(\mathbb{R}^n)$, which is naturally isomorphic to $\mathbb{R}^{n\times n} \times \mathbb{R}^n$.
For any $A\in\mathbb R^{n\times n}$ and $t\in\mathbb R$, write
\begin{align}\label{eq:Kt_definition}
K_t(A):=\int_0^t e^{sA}\,ds.
\end{align}

\noindent Given a vector field $X \in \mathcal{T}(\mathbb{R}^n)$, let $\varphi_t : \mathbb{R}^n \to \mathbb{R}^n
$ denote its associated flow map, defined as the solution to the ordinary differential equation
\begin{equation}
\frac{d}{dt}\varphi_t(x)
=
X(\varphi_t(x)),
\qquad
\varphi_0(x)=x.
\label{eq:flow_equation}
\end{equation}
For affine vector fields, the flow $\varphi_t$ forms a one-parameter family of affine transformations on $\mathbb{R}^n$. 
The flow induces a natural action on measures through push-forward. Given a measure $\mu$ on $\mathbb{R}^n$, the push-forward measure $\mu_t := (\varphi_t)_*(\mu)$ is defined by
\begin{equation}
(\varphi_t)_*(\mu)(A)
=
\mu(\varphi_t^{-1}(A))
\label{eq:pushforward_measure}
\end{equation}
for every measurable set $A \subseteq \mathbb{R}^n$.
Equivalently, for every $\mu$-integrable test function $f : \mathbb{R}^n \to \mathbb{R}$,
\begin{align}
\label{eq:pushforward_prop_measures}
\int_{\mathbb{R}^n} f(x)\, d\mu_t(x)
=
\int_{\mathbb{R}^n} f(\varphi_t(x))\, d\mu(x).    
\end{align}

\vspace{3mm}
\noindent Fix $\mu=\mathcal{N}(m,\Sigma)\in\mathcal{G}$. For an affine vector field
\[
X(x)=Ax+b,
\qquad
A\in\mathbb{R}^{n\times n},\ \ b\in\mathbb{R}^n,
\]
let $\varphi_t$ denote its flow and define the evolved measure
\[
\mu_t := (\varphi_t)_*\mu.
\]
The affine flow $\varphi_t$ has invertible linear part $e^{tA}$ and is therefore a diffeomorphism. Hence $\mu_t$ is again a nondegenerate Gaussian measure. Both $\mu_t$ and $\mu$ have strictly positive Lebesgue densities, so $\mu_t\sim\mu$ and the Radon--Nikodym derivative below exists. Define the $\mu$-divergence of $X$ by
\begin{align}
\label{eq:div_mu_defn}
    \operatorname{div}_\mu(X)
:=-
\left.\frac{d}{dt}\right|_{t=0}\frac{d\mu_t}{d\mu}.
\end{align}
For affine $X$, the following proposition shows that this divergence is a centered quadratic function and hence belongs to $\mathcal Q_\mu^0$.

\begin{proposition}\label{prop:div_mu}
Let $\mu=\mathcal{N}(m,\Sigma)\in\mathcal{G}$ and let $X(x)=Ax+b$ be an affine vector field on $\mathbb{R}^n$. Then
\begin{align*}
    \operatorname{div}_\mu(X)
&=-g(\;\cdot\; ;m,\Sigma,m_v,\Sigma_v),
\end{align*}
where $m_v=Am+b$ and $\Sigma_v=A\Sigma + \Sigma A^\top$.
In particular, $\operatorname{div}_\mu(X)\in \mathcal{Q}_\mu^0.$
\end{proposition}
\begin{proof}
The flow of the affine vector field $X(x)=Ax+b$ is given by
\begin{align}\label{eq:affine_flow_solution}
\varphi_t(x)=e^{tA}x+K_t(A)b.
\end{align}
Hence, if $\mu=\mathcal N(m,\Sigma)$, then $\mu_t=(\varphi_t)_*\mu$ is again Gaussian, with
\begin{align}\label{eq:affine_gaussian_evolution}
\mu_t=\mathcal N(m_t,\Sigma_t),\qquad
m_t=e^{tA}m+K_t(A)b,\qquad
\Sigma_t=e^{tA}\Sigma e^{tA^\top}.
\end{align}
Differentiating at $t=0$ gives
\[
\dot m_0=Am+b,
\qquad
\dot\Sigma_0=A\Sigma+\Sigma A^\top.
\]
By the formula computed in the previous subsection for the derivative of a Gaussian density with respect to its parameters,
\[
\left.\frac{d}{dt}\right|_{t=0}\frac{d\mu_t}{d\mu}
=
\left.\frac{d}{dt}\right|_{t=0}\frac{p(\cdot;m_t,\Sigma_t)}{p(\cdot;m,\Sigma)}
=
g(\cdot\,;m,\Sigma,\dot m_0,\dot\Sigma_0).
\]
Therefore
\[
\operatorname{div}_\mu(X)
=
-g(\cdot\,;m,\Sigma,Am+b,A\Sigma+\Sigma A^\top).
\]
Since this is of the form $-g(\cdot\,;m,\Sigma,m_v,\Sigma_v)$, it lies in $\mathcal Q_\mu^0$.
\end{proof}
\noindent Accordingly, we obtain a linear map
\[
\operatorname{div}_\mu : \mathcal{T}_a(\mathbb{R}^n)\longrightarrow \mathcal{Q}_\mu^0,
\qquad
X\longmapsto \operatorname{div}_\mu(X),
\]
which associates to each affine vector field its centered quadratic weighted divergence.
Note that
\[
A\Sigma+\Sigma A^\top \in S_n.
\]
Since $\Sigma\in S_n^+$ is symmetric positive definite, the Lyapunov operator
\begin{equation}\label{eq:lyapunov_operator}
    \mathcal{L}_\Sigma : S_n \to S_n,
\qquad
B \longmapsto B\Sigma+\Sigma B,
\end{equation}
is an isomorphism. Indeed, diagonalize $\Sigma=U\operatorname{diag}(\lambda_1,\ldots,\lambda_n)U^\top$. For $C\in S_n$, the equation $B\Sigma+\Sigma B=C$ becomes
\[
(U^\top BU)_{ij}
=\frac{(U^\top CU)_{ij}}{\lambda_i+\lambda_j},
\]
which has a unique symmetric solution because every $\lambda_i+\lambda_j$ is positive. This is also the standard spectral criterion for the Sylvester equation \cite{higham1993perturbation}. Hence there exists a unique matrix $A_{\mathrm{symm}}\in S_n$ such that
\[
A_{\mathrm{symm}}\Sigma+\Sigma A_{\mathrm{symm}}^\top
=
A\Sigma+\Sigma A^\top.
\]
Define
\[
A_{\mathrm{df}} := A-A_{\mathrm{symm}},
\qquad
b_{\mathrm{df}} := -A_{\mathrm{df}}m,
\qquad
b_{\mathrm{symm}} := b+A_{\mathrm{df}}m.
\]
Then
\[
X_{\mathrm{grad}}(x):=A_{\mathrm{symm}}x+b_{\mathrm{symm}},
\qquad
X_{\mathrm{divfree}}(x):=A_{\mathrm{df}}x+b_{\mathrm{df}}=A_{\mathrm{df}}(x-m),
\]
so that
\[
X=X_{\mathrm{grad}}+X_{\mathrm{divfree}}.
\]
Moreover, one checks that
\[
\operatorname{div}_\mu(X_{\mathrm{divfree}})=0
\qquad \text{and} \qquad
\operatorname{div}_\mu(X_{\mathrm{grad}})=\operatorname{div}_\mu(X).
\]

This decomposition is in fact an orthogonal decomposition with respect to the $L^2(\mu)$ inner product on vector fields, as summarized in the next proposition.

\begin{proposition}
Fix $\mu=\mathcal{N}(m,\Sigma)\in\mathcal{G}$ and let $X\in\mathcal{T}_a(\mathbb{R}^n)$. Then there exist unique vector fields
\[
X_{\mathrm{grad}}\in\{\nabla f\mid f\in\mathcal{Q}\}
\qquad\text{and}\qquad
X_{\mathrm{divfree}}\in\mathcal{T}_a(\mathbb{R}^n)
\]
such that
\[
X=X_{\mathrm{grad}}+X_{\mathrm{divfree}},
\qquad
\operatorname{div}_\mu(X_{\mathrm{divfree}})=0,
\qquad
\operatorname{div}_\mu(X_{\mathrm{grad}})=\operatorname{div}_\mu(X),
\]
and
\[
\langle\!\langle X_{\mathrm{divfree}},X_{\mathrm{grad}}\rangle\!\rangle_\mu^O
:=\mathbb{E}_\mu\!\left[X_{\mathrm{divfree}}(x)^\top X_{\mathrm{grad}}(x)\right]=0.
\]
\end{proposition}

\begin{proof}
The preceding Lyapunov construction gives existence and uniqueness with $X_{\mathrm{grad}}$ restricted to affine gradient fields. It remains only to verify orthogonality. To this end, observe that
\begin{align*}
    \langle\!\langle X_{\mathrm{divfree}},X_{\mathrm{grad}}\rangle\!\rangle_{\mu}^{O}&=\mathbb{E}_{\mu}\!\left[
X_{\mathrm{divfree}}(x)^\top X_{\mathrm{grad}}(x)
\right] \\
&=\mathbb{E}_{\mu}\!\left[\left(A_{\mathrm{df}}(x-m)\right)^\top \left(A_{\mathrm{symm}}x+b_{\mathrm{symm}}\right)\right]\\
&=\mathbb{E}_{\mu}\!\left[(x-m)^\top A_{\mathrm{df}}^\top A_{\mathrm{symm}}x\right]\\
&=\mathbb{E}_{\mu}\!\left[(x-m)^\top A_{\mathrm{df}}^\top A_{\mathrm{symm}}\left(x-m\right)\right]\\
&=\trace{A_{\mathrm{df}}^\top A_{\mathrm{symm}}\Sigma}=-\trace{A_{\mathrm{df}}^\top A_{\mathrm{symm}}\Sigma}\\
&=0,
\end{align*}
where we have used the cyclic property of the trace operator and $A_{\mathrm{df}}\Sigma+\Sigma A_{\mathrm{df}}^\top=0$.
\end{proof}

Consequently, it is enough to consider gradient vector fields, which can be written in the form
\[
X_{\mathrm{grad}}(x)=\nabla f(x),
\qquad
f\in \mathcal{Q}.
\]
Since two functions differing by a constant have the same gradient, we may restrict attention to $\mathcal{Q}/\mathbb{R}$ and arrive at the following proposition.
\begin{proposition}
Let $\mu=\mathcal{N}(m,\Sigma)\in \mathcal{G}$.
The map
\begin{align}
    \label{eq:Laplace_op_defn}
    \Delta_\mu:\mathcal Q/\mathbb R\to \mathcal Q_\mu^0,
\qquad
f+\mathbb R\longmapsto \operatorname{div}_\mu(\nabla f),
\end{align}
is an isomorphism.
More explicitly, if $f\in \mathcal{Q}$ with
\[
f(x)=\frac12 x^\top A x+b^\top x,
\qquad A\in S_n,
\]
then
\[
\Delta_\mu(f+\mathbb R)
=
-g(\cdot\,;m,\Sigma,Am+b,A\Sigma+\Sigma A).
\]
Conversely, for every \(g(\cdot\,;m,\Sigma,m_v,\Sigma_v)\in\mathcal Q_\mu^0\), there exists a unique \(f+\mathbb R\in\mathcal Q/\mathbb R\) such that
\[
\Delta_\mu(f+\mathbb R)=-g(\cdot\,;m,\Sigma,m_v,\Sigma_v).
\]
\end{proposition}

\begin{proof}
Let \(f(x)=\frac12 x^\top A x+b^\top x\) with \(A\in S_n\). Then
\[
\nabla f(x)=Ax+b,
\]
and Proposition~\ref{prop:div_mu} gives
\[
\Delta_\mu(f+\mathbb R)
=
\operatorname{div}_\mu(\nabla f)
=
-g(\cdot\,;m,\Sigma,Am+b,A\Sigma+\Sigma A).
\]

For surjectivity, let \(g(\cdot\,;m,\Sigma,m_v,\Sigma_v)\in\mathcal Q_\mu^0\) be arbitrary. By the Lyapunov isomorphism established above, let \(A=\mathcal L_\Sigma^{-1}(\Sigma_v)\) and define
\[
b:=m_v-Am.
\]
Then with
\[
f(x)=\frac12 x^\top A x+b^\top x
\]
we obtain
\[
\Delta_\mu(f+\mathbb R)
=
-g(\cdot\,;m,\Sigma,m_v,\Sigma_v).
\]
Injectivity follows because if \(\Delta_\mu(f+\mathbb R)=0\), then \(A=0\) and \(b=0\), so \(f\) is constant modulo \(\mathbb R\).
\end{proof}

\subsection{Riemannian Metrics on the Gaussian Manifold}
\subsubsection{Fisher--Rao Metric}
We now introduce the Fisher--Rao metric by translating the $L_2(\mu)$ product defined for functions from the cotangent space to the tangent space.
Consider the \(L^2(\mu)\) product for functions  \(f,g \in \mathcal{Q}_\mu^0\) defined by
\[
\langle\!\langle f,g\rangle\!\rangle_{\mu}^{FR}
= \int_{\mathbb{R}^n} f(x)g(x)\,\mu(dx),
\]
which can be translated to functions in $(f+\mathbb{R}),(g+\mathbb{R}) \in \mathcal{Q}/\mathbb{R}$ as
\begin{align*}
    \langle\!\langle f+\mathbb{R},g+\mathbb{R}\rangle\!\rangle_{\mu}^{FR}
= \int_{\mathbb{R}^n} \left(f(x)-\mathbb{E}_{\mu}[f(x)]\right)\left(g(x)-\mathbb{E}_{\mu}[g(x)]\right)\,\mu(dx).
\end{align*}
Note that since
\begin{align*}
   \langle\!\langle f+\mathbb{R},g+\mathbb{R}\rangle\!\rangle_{\mu}^{FR}= \int_{\mathbb{R}^n} \left(f-\mathbb{E}_{\mu}[f]\right)\left(g-\mathbb{E}_{\mu}[g]\right)d\mu=\int_{\mathbb{R}^n} f da,
\end{align*}
where $a=\left(g-\mathbb{E}_{\mu}[g]\right)\mu$, 
we have that $\langle\!\langle \;\cdot\;,g+\mathbb{R}\rangle\!\rangle_{\mu}^{FR}$ acts on $f+\mathbb{R}$ in the same way as the measure $a$ acts on $f+\mathbb{R}$.
We thus have the isomorphism
\[
\phi_\mu : \mathcal{Q}/\mathbb{R} \to T_\mu \mathcal{G},
\qquad
f+\mathbb{R} \mapsto \left(f-\mathbb{E}_{\mu}[f]\right)\mu,
\]
with the inverse
\[
\phi_\mu^{-1} : T_\mu \mathcal{G} \to \mathcal{Q}/\mathbb{R},
\qquad
a \mapsto \frac{da}{d\mu}+\mathbb{R}.
\]
Write $a,b \in T_{\mu} \mathcal{G}$ in Gaussian coordinates as $(m_a,\Sigma_a)$ and $(m_b,\Sigma_b)$, respectively, and let
$f_a+\mathbb R$ and $f_b+\mathbb R$ be their cotangent representatives:
\[
f_a+\mathbb R:=\phi_\mu^{-1}(a)
=\frac{da}{d\mu}+\mathbb R,
\qquad
f_b+\mathbb R:=\phi_\mu^{-1}(b)
=\frac{db}{d\mu}+\mathbb R.
\]
Choose the centered representatives
\[
f_a=\frac{da}{d\mu}
=g(\,\cdot\,;m,\Sigma,m_a,\Sigma_a),
\qquad
f_b=\frac{db}{d\mu}
=g(\,\cdot\,;m,\Sigma,m_b,\Sigma_b).
\]
The Fisher--Rao metric on the tangent space is the pullback of the cotangent inner product through $\phi_\mu^{-1}$. Therefore,
\begin{align*}
\langle\!\langle a,b\rangle\!\rangle_{\mu}^{FR}
&:=
\left\langle\!\left\langle f_a+\mathbb R,f_b+\mathbb R\right\rangle\!\right\rangle_{\mu}^{FR}\\
&=\int_{\mathbb{R}^n}f_a(x)f_b(x)\,d\mu(x)\\
&=m_a^\top\Sigma^{-1}m_b
+\frac12\operatorname{Tr}\!\bigl(\Sigma^{-1}\Sigma_a\Sigma^{-1}\Sigma_b\bigr).
\end{align*}
Indeed, with $z=x-m$, the first centered representative is
\[
f_a(x)
=m_a^\top\Sigma^{-1}z
+\frac12\left(z^\top\Sigma^{-1}\Sigma_a\Sigma^{-1}z
-\operatorname{Tr}(\Sigma^{-1}\Sigma_a)\right)
\]
and $f_b$ is given by the analogous expression with $a$ replaced by $b$. The mixed linear--quadratic terms have zero expectation, while for centered Gaussian $z$,
\[
\operatorname{Cov}(z^\top Cz,z^\top Dz)
=2\operatorname{Tr}(C\Sigma D\Sigma),
\qquad C,D\in S_n,
\]
which gives the displayed coordinate formula.

\subsubsection{Otto Metric}\label{sec:Otto_metric}
We now follow a similar recipe for the Otto metric by translating the $L_2(\mu)$ product defined for gradient vector fields from the cotangent space to the tangent space.
For $f,g \in \mathcal{Q}/\mathbb{R}$, let
\[
\langle\!\langle f+\mathbb{R},g+\mathbb{R}\rangle\!\rangle_{\mu}^{O}
= \int_{\mathbb{R}^n} \nabla f(x)^\top \nabla g(x)\,\mu(dx).
\]
Since 
\begin{alignat*}{2}
\langle\!\langle f+\mathbb{R},g+\mathbb{R}\rangle\!\rangle_{\mu}^{O}
&= \langle\!\langle \nabla f,\nabla g\rangle\!\rangle_{\mu}^{O}
 = \int_{\mathbb{R}^n} \nabla f^{\top}\nabla g \, d\mu &\\
&= \int_{\mathbb{R}^n} \left.\frac{d}{dt}\right|_{t=0} (f\circ \varphi_t)\, d\mu
&\text{(where } \varphi_t \text{ is the flow of } \nabla g)\\
&= \left.\frac{d}{dt}\right|_{t=0}\int_{\mathbb{R}^n} (f\circ \varphi_t)\, d\mu 
&\text{(because $f\circ \varphi_t$ is quadratic)} \\
&= \left.\frac{d}{dt}\right|_{t=0}\int_{\mathbb{R}^n} f\, d\mu_t 
&\text{(because of \eqref{eq:pushforward_prop_measures})}\\
&= \left.\frac{d}{dt}\right|_{t=0}\int_{\mathbb{R}^n} f\,\left(\frac{d\mu_t}{d\mu}\right) d\mu & \\
&= \int_{\mathbb{R}^n} f\, \left.\frac{d}{dt}\right|_{t=0} \left(\frac{d\mu_t}{d\mu}\right)\, d\mu &\\
&= \int_{\mathbb{R}^n} f\bigl(-\operatorname{div}_{\mu}(\nabla g)\bigr)\, d\mu
& \text{(using definition \eqref{eq:div_mu_defn})}\\
&= \int_{\mathbb{R}^n} f\bigl(-\Delta_{\mu}(g+\mathbb{R})\bigr)\, d\mu & \text{(using definition \eqref{eq:Laplace_op_defn})}\\
&= \int_{\mathbb{R}^n} f\,da \qquad &\text{(where $a= -\Delta_{\mu}(g+\mathbb{R})\mu$)},\\
&= \int_{\mathbb{R}^n} (f+\mathbb{R})\,da \qquad &\text{since $a\in \mathcal{S}_0^{\infty}$}.
\end{alignat*}
The differentiated integrands are dominated by a fixed polynomial in $x$. Gaussian integrability then justifies differentiation under the integral sign and makes every integral above finite.

Thus, $\langle\!\langle \;\cdot\;,g+\mathbb{R}\rangle\!\rangle_{\mu}^{O}$ acts on $f+\mathbb{R}$ in the same way as the measure $a$ acts on $f+\mathbb{R}$.
With this motivation, define the isomorphism $\phi_{\mu}$ and its inverse as
\begin{align*}
\phi_{\mu} : \mathcal{Q}/\mathbb{R} &\to T_{\mu}\mathcal{G},
& g+\mathbb{R} &\mapsto -\Delta_{\mu}(g+\mathbb{R})\,\mu, \\
\phi_{\mu}^{-1} : T_{\mu}\mathcal{G} &\to \mathcal{Q}/\mathbb{R},
& a &\mapsto -\Delta_{\mu}^{-1}\!\left(\frac{da}{d\mu}\right),
\end{align*}
respectively.
Let $a,b\in T_\mu\mathcal{G}$ have Gaussian coordinate velocities $(m_a,\Sigma_a)$ and $(m_b,\Sigma_b)$, respectively. Write
\[
f_a+\mathbb{R}:=\phi_\mu^{-1}(a)
=-\Delta_\mu^{-1}\!\left(\frac{da}{d\mu}\right),
\qquad
f_b+\mathbb{R}:=\phi_\mu^{-1}(b)
=-\Delta_\mu^{-1}\!\left(\frac{db}{d\mu}\right),
\]
and choose the representatives
\[
f_a(x)=\frac12x^\top A_ax+b_a^\top x,
\qquad
f_b(x)=\frac12x^\top A_bx+b_b^\top x.
\]
Here $A_a,A_b\in S_n$ are the unique symmetric solutions of
\begin{align*}
A_a\Sigma+\Sigma A_a&=\Sigma_a,
& A_am+b_a&=m_a,\\
A_b\Sigma+\Sigma A_b&=\Sigma_b,
& A_bm+b_b&=m_b.
\end{align*}
The Otto metric on the tangent space is the pullback of the cotangent inner product through $\phi_\mu^{-1}$ and is therefore
\begin{align}\label{eq:otto_gaussian_metric}
\langle\!\langle a,b\rangle\!\rangle_\mu^O
&=\left\langle\!\left\langle\phi_\mu^{-1}(a),\phi_\mu^{-1}(b)\right\rangle\!\right\rangle_\mu^O\notag\\
&=\int_{\mathbb{R}^n}(A_ax+b_a)^\top(A_bx+b_b)\,d\mu(x)\\
&=\operatorname{Tr}(A_a^\top A_b\Sigma)+m_a^\top m_b.
\end{align}
The covariance part of \eqref{eq:otto_gaussian_metric} is the Lyapunov
representation of the Gaussian Wasserstein Riemannian metric; see
\cite[Propositions~6--7]{malago2018wasserstein}, which also develops the
corresponding Riemannian exponential, Levi--Civita connection, and parallel
transport.  The $e_1$-connection introduced below uses the same Otto metric,
but is defined as the metric conjugate of the intrinsic mixture connection;
it is not the Levi--Civita connection studied there.

\subsection{Mixture and Exponential Connections and Their Geodesics}\label{sec:connections_and_geodesics}
The intrinsic Gaussian mixture connection is affine in moment coordinates. Its construction, geodesics, and metric-conjugacy properties are developed in Appendix~\ref{sec:app_mixture_duality}. Here we derive the \(e_0\)- and \(e_1\)-connections associated with the Fisher--Rao and Otto metrics, respectively, following \cite{ay2024information,ay2025correction}. The construction is the same in both cases: identify a tangent vector $a$ with its cotangent representative $f_a+\mathbb{R}=\phi_{\mu}^{-1}(a)$, transport that representative trivially as
\begin{align*}
    T_{\mu}^*\mathcal{G}\ni (\mu,f_a+\mathbb{R}) \mapsto (\nu,f_a+\mathbb{R}) \in T_{\nu}^*\mathcal{G},
\end{align*}
and map $f_a+\mathbb{R}$ at $\nu$ back to the tangent vector $\phi_{\nu}(f_a+\mathbb{R})$ at $\nu$. Here $\phi_\mu$ denotes the Fisher--Rao or Otto identification according to the subsection; their definitions are different.

The phrase ``mixture connection'' requires some care after restriction to $\mathcal G$: an ambient density mixture generally leaves the Gaussian family. The intrinsic connection used here remains within $\mathcal G$ because it is defined by affine interpolation of Gaussian moment coordinates. Appendix~\ref{sec:app_mixture_duality} proves that it is conjugate to the $e_0$- and $e_1$-connections with respect to the Fisher--Rao and Otto metrics, respectively. The $e_1$ parallel transport depends only on its endpoints and therefore defines a curvature-flat connection, but this connection is not torsion-free \cite{ay2026torsion}. Thus, ``flat'' in this context refers only to vanishing curvature and not to dual flatness.
\subsubsection{The \texorpdfstring{$e_0$}{e0}-Connection for the Fisher--Rao Metric}

\noindent Following the above described construction for the Fisher--Rao metric gives the parallel transport map
\begin{align*}
    \Pi_{\mu,\nu}^{(e_0)}: T_{\mu}\mathcal{G}\ni (\mu,a) \mapsto (\nu,(\phi_{\nu}\circ \phi_{\mu}^{-1})(a)) \in T_{\nu}\mathcal{G},
\end{align*}
where
\begin{align*}
    (\phi_{\nu}\circ \phi_{\mu}^{-1})(a)=\left(\frac{da}{d\mu}-\mathbb{E}_{\nu}\left[\frac{da}{d\mu}\right]\right)\nu.
\end{align*}
For $\mu=\mathcal{N}(m_\mu,\Sigma_\mu)$ and $\nu=\mathcal{N}(m_\nu,\Sigma_\nu)$, set
\[
g_\mu:=g(\,\cdot\,;m_\mu,\Sigma_\mu,m_v,\Sigma_v).
\]
Then the transport map has the coordinate representation
\begin{align*}
\Pi_{\mu,\nu}^{(e_0)}:(\mu,g_\mu\mu)
\longmapsto\bigl(\nu,(g_\mu-\mathbb{E}_\nu[g_\mu])\nu\bigr).
\end{align*}
The corresponding $e_0$-geodesic $\gamma_0$ is described by
\begin{align}\label{eq:e0-geodesic_ode}
    \dot{\gamma}_0(t)=\Pi_{\mu,\gamma_0(t)}^{(e_0)}a=\left(\frac{da}{d\mu}-\mathbb{E}_{\gamma_0(t)}\left[\frac{da}{d\mu}\right]\right)\gamma_0(t).
\end{align}
In order to write the geodesic equations in coordinates, consider the curve $t\mapsto \mu_t$ as defined in Section~\ref{sec:tgt_cotgt_spaces} and let $(m_a,\Sigma_a)\in \Rn \times S_n$ represent the initial velocity.
The geodesic equation leads to
\begin{align*}
    g(x;m_t,\Sigma_t,\dot{m}_t,\dot{\Sigma}_t)\cdot \mu_t=\left(g(x;m_0,\Sigma_0,m_a,\Sigma_a)-\mathbb{E}_{\mu_t}[g(x;m_0,\Sigma_0,m_a,\Sigma_a)]\right)\mu_t.
\end{align*}
Equating the quadratic and linear terms in $x$ gives
\begin{align*}
\Sigma_t^{-1}\dot{\Sigma}_t\Sigma_t^{-1}&=\Sigma_0^{-1}\Sigma_a\Sigma_0^{-1},\\
    -\Sigma_t^{-1}\dot{\Sigma}_t\Sigma_t^{-1}m_t+\Sigma_t^{-1}\dot{m}_t &=  -\Sigma_0^{-1}\Sigma_a\Sigma_0^{-1}m_0+\Sigma_0^{-1}m_a,
\end{align*}
which can be rearranged to obtain
\begin{align*}
\dot{\Sigma}_t&=\Sigma_t\Sigma_0^{-1}\Sigma_a\Sigma_0^{-1}\Sigma_t,\\
    \dot{m}_t &=  \Sigma_t(\Sigma_0^{-1}\Sigma_a\Sigma_0^{-1}m_t-\Sigma_0^{-1}\Sigma_a\Sigma_0^{-1}m_0)+\Sigma_t \Sigma_0^{-1}m_a.
\end{align*}
These equations can be solved by transforming to the "natural" exponential family coordinates $P_t=\Sigma_t^{-1}$ and $\eta_t=\Sigma_t^{-1}m_t$ to get
\begin{alignat*}{2}
\dot{P}_t&=-P_0\Sigma_a P_0 \quad &\implies& \quad P_t=P_0-t\cdot P_0\Sigma_aP_0,\\
\dot{\eta}_t&=P_0m_a-P_0\Sigma_a\eta_0 \quad &\implies& \quad \eta_t=\eta_0+t\cdot (P_0m_a-P_0\Sigma_a\eta_0),
\end{alignat*}
which are straight-line trajectories in the natural coordinates, on the maximal interval containing zero for which $P_t\succ0$.

If one starts instead with two points $\mu_0=\mathcal{N}(m_0,\Sigma_0)$ and $\mu_1=\mathcal{N}(m_1,\Sigma_1)$ in $\mathcal G$, put
\[
P_j:=\Sigma_j^{-1},
\qquad
\eta_j:=P_jm_j,
\qquad j\in\{0,1\}.
\]
The unique unit-time $e_0$ connector is the straight line in natural coordinates,
\begin{equation}\label{eq:e0_endpoint_connector}
\begin{aligned}
P_t&=(1-t)P_0+tP_1,
&
\eta_t&=(1-t)\eta_0+t\eta_1,
&
\gamma_0(t)&=\mathcal N(P_t^{-1}\eta_t,P_t^{-1}).
\end{aligned}
\end{equation}
In coordinate free notation, the $e_0$-geodesic can be written as
\begin{align}\label{eq:e0-geodesic}
    \gamma_0(t)=\frac{\left(\frac{d\mu_1}{d\mu_0}\right)^t}{Z(t)}\mu_0,
    \qquad
    Z(t)=\int_{\Rn}\left(\frac{d\mu_1}{d\mu_0}\right)^t d\mu_0.
\end{align}

\subsubsection{The \texorpdfstring{$e_1$}{e1}-Connection for the Otto Metric}
We now investigate what happens when we use the Otto metric instead of the Fisher--Rao metric.

The parallel transport map defined on the tangent space is
\begin{align*}
    \Pi_{\mu,\nu}^{(e_1)}: T_{\mu}\mathcal{G}\ni (\mu,a) \mapsto (\nu,(\phi_{\nu}\circ \phi_{\mu}^{-1})(a)) \in T_{\nu}\mathcal{G},
\end{align*}
where
\begin{align*}
    (\phi_{\nu}\circ \phi_{\mu}^{-1})(a)=\left(\left(\Delta_{\nu} \circ \Delta_{\mu}^{-1}\right)\left(\frac{da}{d\mu}\right)\right)\nu.
\end{align*}
This can be computed explicitly for $\mu=\mathcal{N}(m_\mu,\Sigma_\mu)$ and $\nu=\mathcal{N}(m_\nu,\Sigma_\nu)$. As in the preceding $e_0$ calculation, for a tangent vector with coordinate velocity $(m_v,\Sigma_v)$ set
\[
g_\mu:=g(\,\cdot\,;m_\mu,\Sigma_\mu,m_v,\Sigma_v)\in\mathcal{Q}_\mu^0.
\]
The cotangent representative of $g_\mu\mu$ is
\[
-\Delta_\mu^{-1}(g_\mu)=f+\mathbb{R}.
\]
Choose the representative
\begin{align}\label{eq:potential}
f(x)=\frac12x^\top Ax+b^\top x,
\qquad A\in S_n,\quad b\in\Rn,
\end{align}
where $A$ and $b$ are uniquely determined by
\[
A\Sigma_\mu+\Sigma_\mu A=\Sigma_v,
\qquad
Am_\mu+b=m_v.
\]
By Proposition~\ref{prop:div_mu},
\begin{align*}
-\Delta_\nu(f+\mathbb{R})
=g(\,\cdot\,;m_\nu,\Sigma_\nu,
Am_\nu+b,
A\Sigma_\nu+\Sigma_\nu A).
\end{align*}
Therefore, the transport map in coordinates can be written as
\begin{align*}
   \Pi_{\mu,\nu}^{(e_1)}: (\mu,g(\;\cdot\; ;m_{\mu},\Sigma_{\mu},m_v,\Sigma_v)\mu) \mapsto \Big(\nu,g(\;\cdot\; ;m_{\nu},\Sigma_{\nu},Am_{\nu}+b,A\Sigma_{\nu}+\Sigma_{\nu}A)\nu\Big),
\end{align*}
where $A,b$ satisfy
\begin{align}
    A\Sigma_{\mu}+\Sigma_{\mu}A&=\Sigma_v, \label{eq:A_mu_b_mu}\\
    Am_{\mu}+b&=m_v. \label{eq:A_mu_b_mu_mean}
\end{align}

\noindent The corresponding $e_1$-geodesic $\gamma_1$ is described by
\begin{align}\label{eq:e1-geodesic}
    \dot{\gamma}_1(t)=\Pi_{\mu,\gamma_1(t)}^{(e_1)}a=\left(\left(\Delta_{\gamma_1(t)} \circ \Delta_{\mu}^{-1}\right)\left(\frac{da}{d\mu}\right)\right)\gamma_1(t).
\end{align}
In order to write the geodesic equations in coordinates, consider the curve $t\mapsto \mu_t$ as defined in Section~\ref{sec:tgt_cotgt_spaces} and let $(m_a,\Sigma_a)\in \Rn \times S_n$ represent the initial velocity.
The geodesic equation leads to
\begin{align*}
    g(x;m_t,\Sigma_t,\dot{m}_t,\dot{\Sigma}_t) \mu_t=\left(g(\;\cdot\; ;m_t,\Sigma_t,Am_t+b,A\Sigma_t+\Sigma_tA)\right)\mu_t.
\end{align*}
Equating the quadratic and linear terms in $x$ gives
\begin{align*}
\dot{\Sigma}_t&=A\Sigma_t+\Sigma_tA,\\
    \dot{m}_t &= Am_t+b,
\end{align*}
where $A$ and $b$ are given in \eqref{eq:A_mu_b_mu}--\eqref{eq:A_mu_b_mu_mean}.

This coordinate calculation has an equivalent pushforward interpretation. Because the $e_1$ transport keeps the cotangent class $f+\mathbb R$ fixed, the vector field $\nabla f(x)=Ax+b$ is autonomous. Let $\varphi_t$ be its forward gradient flow\footnote{The terminology ``gradient flow'' here uses the positive-gradient convention.},
\begin{align}\label{eq:dynamics}
\frac{d}{dt}\varphi_t(x)
=\nabla f(\varphi_t(x))
=A\varphi_t(x)+b,
\qquad
\varphi_0(x)=x.
\end{align}
The flow map is
\begin{align}\label{eq:e1_flow_map}
\varphi_t(x)=e^{tA}x+K_t(A)b.
\end{align}
Since $\varphi_{t+s}=\varphi_s\circ\varphi_t$, the definition of $\Delta_\rho$, applied at $\rho=(\varphi_t)_*\mu$, gives
\[
\frac{d}{dt}(\varphi_t)_*\mu
=-\Delta_{(\varphi_t)_*\mu}(f+\mathbb R)(\varphi_t)_*\mu,
\]
which is precisely the $e_1$-geodesic equation because $f+\mathbb R$ is the fixed transported cotangent representative. Hence, writing $\mu_0=\mu=\mathcal N(m_0,\Sigma_0)$,
\begin{align}\label{eq:e1_global_geodesic}
\gamma_1(t)=(\varphi_t)_*\mu_0
=\mathcal N(m_t,\Sigma_t),
\qquad
m_t=e^{tA}m_0+K_t(A)b,
\qquad
\Sigma_t=e^{tA}\Sigma_0e^{tA}.
\end{align}
For a symmetric matrix $A$, write $A^\dagger$ for its Moore--Penrose pseudoinverse and
\[
\Perp{A}:=\operatorname{proj}_{\ker A}=I-AA^\dagger=I-A^\dagger A
\]
for its orthogonal null-space projector. Diagonalizing $A$ gives the useful identity
\begin{align}\label{eq:Kt_symmetric_formula}
K_t(A)=(e^{tA}-I)A^\dagger+t\Perp{A},
\end{align}
which also covers singular $A$. Conversely, every Gaussian $e_1$-geodesic has this form: the Lyapunov equation in \eqref{eq:A_mu_b_mu} assigns a unique symmetric $A$, and \eqref{eq:A_mu_b_mu_mean} then assigns a unique $b$, to each initial tangent vector.

\subsubsection{Geodesic Completeness and Endpoint Connectivity}
The explicit equations also settle the global-in-time question. We call the $e_1$-connection geodesically complete when every geodesic with prescribed initial position and velocity is defined for every $t\in\mathbb R$. This is distinct from metric completeness and does not assign a point to $t=\infty$.

\begin{proposition}[Completeness of the Gaussian $e_1$-connection]\label{prop:e1_completeness}
Let $\mu_0=\mathcal N(m_0,\Sigma_0)\in\mathcal G$ and let $(m_a,\Sigma_a)\in\mathbb R^n\times S_n$ be an initial velocity. Let $A\in S_n$ and $b\in\mathbb R^n$ be the unique solutions of
\begin{align}\label{eq:e1_initial_generator}
A\Sigma_0+\Sigma_0A=\Sigma_a,
\qquad
Am_0+b=m_a.
\end{align}
Then the curve in \eqref{eq:e1_global_geodesic} is the unique $e_1$-geodesic with these initial data, and it belongs to $\mathcal G$ for every $t\in\mathbb R$. Hence the $e_1$-connection on $\mathcal G$ is geodesically complete.

A limit as $t\to+\infty$ exists in the closure of the Gaussian family if and only if
\begin{align}\label{eq:e1_forward_limit_condition}
A\preceq0
\qquad\textnormal{and}\qquad
\Perp{A}b=0.
\end{align}
When these conditions hold,
\begin{align}\label{eq:e1_forward_limit}
\gamma_1(t)\Longrightarrow
\mathcal N\!\left(\Perp{A}m_0-A^\dagger b,\,\Perp{A}\Sigma_0\Perp{A}\right),
\end{align}
where the Gaussian on the right may be degenerate. In particular, if $A\prec0$, the limit is the Dirac measure $\delta_{-A^{-1}b}$. The limit lies in $\mathcal G$ only for the constant geodesic $A=0$, $b=0$.
\end{proposition}

\begin{proof}
The Lyapunov isomorphism associated with $\Sigma_0\succ0$ gives the unique symmetric $A$ in \eqref{eq:e1_initial_generator}, after which $b$ is unique. Formula \eqref{eq:e1_global_geodesic} gives the asserted geodesic. Since $e^{tA}$ is invertible for every $t\in\mathbb R$,
\[
z^\top\Sigma_tz=(e^{tA}z)^\top\Sigma_0(e^{tA}z)>0
\qquad(z\ne0),
\]
which proves that the solution remains in $\mathcal G$ for every $t\in\mathbb R$. Uniqueness follows from \eqref{eq:e1_initial_generator}.

For the final assertion, diagonalize $A$. A positive eigenvalue makes the variance in its eigendirection grow exponentially so that the limit does not exist. If $A\preceq0$, then $e^{tA}\to\Perp{A}$; when $\Perp{A}b\ne0$, the covariance remains bounded but the projected mean has the unbounded term $t\Perp{A}b$, again precluding the existence of a limit. When $\Perp{A}b=0$, $K_t(A)b\to-A^\dagger b$ and the mean and covariance converge to those in \eqref{eq:e1_forward_limit}. These observations also prove necessity. A nondegenerate limiting covariance would require $\Perp{A}=I$, hence $A=0$, and then \eqref{eq:e1_forward_limit_condition} requires $b=0$.
\end{proof}

The corresponding statement for the Fisher--Rao $e_0$-connection is different. For the same initial Gaussian and coordinate velocity, put
\[
C:=\Sigma_0^{-1/2}\Sigma_a\Sigma_0^{-1/2}.
\]
The natural precision coordinate along the $e_0$-geodesic is
\begin{align}\label{eq:e0_precision_initial_geodesic}
\Sigma_t^{-1}
=\Sigma_0^{-1/2}(I-tC)\Sigma_0^{-1/2}.
\end{align}
If $c_1,\ldots,c_n$ are the eigenvalues of $C$, its maximal interval is $(\tau_-,\tau_+)$, where
\begin{align}\label{eq:e0_maximal_interval}
\tau_-&=\max_{c_i<0}\frac1{c_i},
&
\tau_+&=\min_{c_i>0}\frac1{c_i},
\end{align}
with the value $-\infty$, respectively $+\infty$, when the corresponding index set is empty. Thus an $e_0$-geodesic is defined on all of $\mathbb R$ exactly when $\Sigma_a=0$, and the $e_0$-connection is not geodesically complete. It is defined for every $t\geq0$ exactly when $\Sigma_a\preceq0$. For example, $\Sigma_0=I$ and $\Sigma_a=I$ give $\Sigma_t^{-1}=(1-t)I$: the precision approaches the boundary at $t=1$, while the covariance diverges. By contrast, Proposition~\ref{prop:e1_completeness} applies to every $e_1$ initial velocity.


The preceding discussion treats prescribed initial data. We now determine the unique quadratic $e_1$ generator joining two prescribed nondegenerate Gaussian endpoints in unit time.

For $B\succ0$, $\log B$ is the principal symmetric logarithm; for $B\succeq0$, $B^{1/2}$ is the unique symmetric positive-semidefinite square root. The endpoint calculation below first appeared in \cite{datar2026wkl}.

\begin{lemma}[Unique quadratic $e_1$ endpoint generator]\label{lemm:formulae_for_A_b}
Let $\mu=\mathcal{N}(m_0,\Sigma_0)$ and $\nu=\mathcal{N}(m_1,\Sigma_1)$ be two nondegenerate Gaussian measures with $\Sigma_0\succ 0$ and $\Sigma_1\succ 0$.
Define
\begin{align}\label{eq:e1_endpoint_parameters}
R:=\Sigma_0^{-\frac{1}{2}}
\left(\Sigma_0^{\frac{1}{2}}\Sigma_1\Sigma_0^{\frac{1}{2}}\right)^{\frac{1}{2}}
\Sigma_0^{-\frac{1}{2}},
\end{align}
and let the quadratic function $f$ in \eqref{eq:potential} have coefficients
\begin{align}
A&=\log R,\label{eq:soln_A}\\
b&=K_1(A)^{-1}(m_1-Rm_0).\label{eq:soln_b}
\end{align}
Then the forward gradient flow \eqref{eq:dynamics} satisfies $(\varphi_1)_*\mu=\nu$.
The resulting connector is
\begin{align}\label{eq:e1_endpoint_connector}
\gamma_1(t)
=\mathcal N\!\left(R^tm_0+K_t(A)b,\,R^t\Sigma_0R^t\right),
\qquad 0\leq t\leq1.
\end{align}
Furthermore, $f$ is unique among functions of the form \eqref{eq:potential}; if additive constants are admitted, it is unique up to such a constant.
\end{lemma}

\begin{proof}
By \eqref{eq:e1_global_geodesic}, the endpoint conditions reduce to
\[
\Sigma_1=e^A\Sigma_0e^A,
\qquad
m_1=e^Am_0+K_1(A)b.
\]
Set
\[
C:=\left(\Sigma_0^{\frac12}\Sigma_1\Sigma_0^{\frac12}\right)^{\frac12}.
\]
Both $C$ and the matrix $R$ in \eqref{eq:e1_endpoint_parameters} are positive definite. Consequently, the principal logarithm $A=\log R$ is real and symmetric and satisfies $e^A=R$. Moreover,
\begin{align}
\label{eq:riccatti}
R\Sigma_0R
=\Sigma_0^{-\frac12}C^2\Sigma_0^{-\frac12}
=\Sigma_1.
\end{align}
Thus the covariance at time one is $\Sigma_1$.

To prove uniqueness, suppose that $X\succ0$ also satisfies $X\Sigma_0X=\Sigma_1$. Then
\[
Y:=\Sigma_0^{\frac12}X\Sigma_0^{\frac12}\succ0,
\qquad
Y^2=\Sigma_0^{\frac12}\Sigma_1\Sigma_0^{\frac12}.
\]
Uniqueness of the positive-definite square root gives $Y=C$ and hence $X=R$. Since the matrix exponential is a bijection from symmetric matrices to positive-definite matrices, with inverse given by the principal logarithm, $A=\log R$ is the unique symmetric generator.

It remains to consider the mean. If $A=U\operatorname{diag}(\lambda_1,\ldots,\lambda_n)U^\top$, then \eqref{eq:Kt_definition} gives
\[
K_1(A)
=U\operatorname{diag}\bigl(k(\lambda_1),\ldots,k(\lambda_n)\bigr)U^\top,
\qquad
k(\lambda):=
\begin{cases}
\dfrac{e^\lambda-1}{\lambda},&\lambda\neq0,\\[5pt]
1,&\lambda=0.
\end{cases}
\]
Since $k(\lambda)>0$ for every $\lambda\in\mathbb{R}$, $K_1(A)\succ0$ and is invertible. Hence the mean endpoint condition above has the unique solution given in \eqref{eq:soln_b}. This completes the proof.
\end{proof}

\begin{proposition}[Global endpoint connectivity]\label{prop:e_geodesic_connectivity}
Every ordered pair of points
\[
\mu_0=\mathcal N(m_0,\Sigma_0),
\qquad
\mu_1=\mathcal N(m_1,\Sigma_1)
\qquad(\Sigma_0,\Sigma_1\succ0)
\]
is joined, for $0\leq t\leq1$, by unique $e_1$- and $e_0$-geodesics satisfying $\gamma_i(0)=\mu_0$ and $\gamma_i(1)=\mu_1$. The $e_1$ endpoint generator and connector are given in Lemma~\ref{lemm:formulae_for_A_b}, while the $e_0$ connector is given in \eqref{eq:e0_endpoint_connector}. Thus both connections are globally geodesically connected on $\mathcal G$, although only the $e_1$-connection is geodesically complete.
\end{proposition}

\begin{proof}
For $e_1$, existence and uniqueness follow from Lemma~\ref{lemm:formulae_for_A_b}. The symmetry requirement is intrinsic: Gaussian tangent vectors are represented by gradient fields, so nonsymmetric affine maps that happen to push one Gaussian to another do not supply additional $e_1$-geodesics. For $e_0$, the natural coordinates $(P,\eta)=(\Sigma^{-1},\Sigma^{-1}m)$ are affine. The cone $S_n^+$ is convex, so the connector in \eqref{eq:e0_endpoint_connector} remains in $\mathcal G$ throughout $[0,1]$ and is the only affine-coordinate line with the prescribed endpoints. The completeness distinction follows from Proposition~\ref{prop:e1_completeness} and \eqref{eq:e0_maximal_interval}.
\end{proof}

\subsubsection{Geodesics in the Univariate Case}
Let
\[
 \mu=\mathcal N(m_0,v_0),
 \qquad
 \nu=\mathcal N(m_1,v_1),
 \qquad
 v_j>0,\quad \sigma_j:=\sqrt{v_j}>0,
 \qquad
 d=m_1-m_0.
\]
\begin{samepage}
In the common mean--variance chart $(m,v)=(m,\sigma^2)$, the geodesics of
the intrinsic Gaussian mixture, $e_0$, and $e_1$ connections joining these
endpoints have the respective coordinate representations
\begin{align}
 m_t^{(m)}&=m_0+td,
 &v_t^{(m)}&=(1-t)v_0+tv_1+t(1-t)d^2,
 \label{eq:univariate_m_connector}\\
 m_t^{(e_0)}&=v_t^{(e_0)}
 \left(\frac{(1-t)m_0}{v_0}+\frac{tm_1}{v_1}\right),
 &v_t^{(e_0)}&=\left(\frac{1-t}{v_0}+\frac{t}{v_1}\right)^{-1},
 \label{eq:univariate_e0_connector}\\
 m_t^{(e_1)}&=(1-q_t)m_0+q_tm_1,
 &v_t^{(e_1)}&=v_0r^{2t},
 \label{eq:univariate_e1_connector}
\end{align}
\end{samepage}
where
\[
 r:=\frac{\sigma_1}{\sigma_0},
 \qquad
 q_t:=
 \begin{cases}
  \dfrac{r^t-1}{r-1},&r\ne1,\\[5pt]
  t,&r=1.
 \end{cases}
\]
The first curve is the $m$-geodesic given in Appendix~\ref{sec:app_mixture_duality}, while the other two are the $e_0$- and $e_1$-geodesics given in \eqref{eq:e0_endpoint_connector} and \eqref{eq:e1_endpoint_connector}. Figure~\ref{fig:univariate_connectors_mixture_coordinates} and Figure~\ref{fig:univariate_connectors_e0_coordinates} uses the endpoints
\[
 \mu=\mathcal N(0,1),
 \qquad
 \nu=\mathcal N(3,16)
\]
and plots the three connectors in the expectation coordinates $(m,M)$ in the left panel and the $e_0$-natural coordinates $(\eta,P)$ in the right panel, respectively, where
\[
 M:=m^2+v,
 \qquad
 P:=v^{-1},
 \qquad
 \eta:=Pm=\frac{m}{v}.
\]
The expectation coordinates $(m,M)$ are affine for the intrinsic mixture connection, and $(P,\eta)$ are affine for the $e_0$-connection which can be seen from the respective straight-line segments in each panel. The $e_1$-connector is not affine in either coordinate system. 
\begin{figure}[H]
\centering
\begin{subfigure}[t]{0.485\textwidth}
\centering
\begin{tikzpicture}
\begin{axis}[
    width=\linewidth,
    height=5.0cm,
    xmin=-0.15, xmax=3.2,
    ymin=0, ymax=26,
    xlabel={mean $m$},
    ylabel={raw moment $M=m^2+\sigma^2$},
    axis lines=left,
    samples=160,
    domain=0:1,
    grid=major,
    grid style={gray!18},
    tick label style={font=\scriptsize},
    label style={font=\small},
    clip=false
]
    \addplot[black, solid, line width=1.2pt]
        ({3*x}, {1+24*x});

    \addplot[blue!70!black, dashed, line width=1.2pt]
        ({3*x/(16-15*x)},
         {(3*x/(16-15*x))^2+16/(16-15*x)});

    \addplot[red!75!black, dash dot, line width=1.2pt]
        ({exp(x*ln(4))-1},
         {(exp(x*ln(4))-1)^2+exp(x*ln(16))});

    \addplot[only marks, mark=o, mark size=2pt, black, fill=white]
        coordinates {(1.5,13) (0.176470588,1.913494810) (1,5)};

    \addplot[only marks, mark=*, mark size=2.2pt, black]
        coordinates {(0,1) (3,25)};
    \node[anchor=south east, font=\scriptsize] at (axis cs:0.1,1.45)
        {$\mu$};
    \node[anchor=south east, font=\scriptsize] at (axis cs:2.96,25.2)
        {$\nu$};
\end{axis}
\end{tikzpicture}
\caption{Mixture expectation coordinates $(m,M)$.}
\label{fig:univariate_connectors_mixture_coordinates}
\end{subfigure}
\hfill
\begin{subfigure}[t]{0.485\textwidth}
\centering
\begin{tikzpicture}
\begin{axis}[
    width=\linewidth,
    height=5.0cm,
    xmin=-0.01, xmax=0.27,
    ymin=0.03, ymax=1.05,
    xlabel={$\eta=m/\sigma^2$},
    ylabel={precision $P=1/\sigma^2$},
    axis lines=left,
    samples=160,
    domain=0:1,
    grid=major,
    grid style={gray!18},
    tick label style={font=\scriptsize},
    label style={font=\small},
    clip=false
]
    \addplot[black, solid, line width=1.2pt]
        ({3*x/(1+24*x-9*x^2)},
         {1/(1+24*x-9*x^2)});

    \addplot[blue!70!black, dashed, line width=1.2pt]
        ({3*x/16}, {1-15*x/16});

    \addplot[red!75!black, dash dot, line width=1.2pt]
        ({exp(-x*ln(4))-exp(-x*ln(16))},
         {exp(-x*ln(16))});

    \addplot[only marks, mark=o, mark size=2pt, black, fill=white]
        coordinates {(0.139534884,0.093023256)
                     (0.09375,0.53125) (0.25,0.25)};

    \addplot[only marks, mark=*, mark size=2.2pt, black]
        coordinates {(0,1) (0.1875,0.0625)};
    \node[anchor=south west, font=\scriptsize] at (axis cs:0.006,1.005)
        {$\mu$};
    \node[anchor=south west, font=\scriptsize] at (axis cs:0.192,0.070)
        {$\nu$};
\end{axis}
\end{tikzpicture}
\caption{$e_0$ affine coordinates $(\eta,P)$.}
\label{fig:univariate_connectors_e0_coordinates}
\end{subfigure}

\medskip
\begin{tikzpicture}[baseline=-0.5ex]
  \draw[black, solid, line width=1.2pt] (0,0)--(0.65,0);
  \node[anchor=west, font=\scriptsize] at (0.75,0)
      {intrinsic mixture $\gamma^{(m)}$};
  \draw[blue!70!black, dashed, line width=1.2pt] (4.2,0)--(4.85,0);
  \node[anchor=west, font=\scriptsize] at (4.95,0)
      {$e_0$-connection $\gamma_0$};
  \draw[red!75!black, dash dot, line width=1.2pt] (7.65,0)--(8.30,0);
  \node[anchor=west, font=\scriptsize] at (8.40,0)
      {$e_1$-connection $\gamma_1$};
\end{tikzpicture}

\caption{The three connection geodesics from
$\mu=\mathcal N(0,1)$ to $\nu=\mathcal N(3,16)$ in the
mixture-affine (left) and $e_0$-affine (right) coordinate charts. In each
panel, the geodesic native to that chart is a straight segment traversed
linearly. Open circles mark the shared numerical value $t=\tfrac12$
of the three endpoint-normalized parameters.}
\label{fig:univariate_gaussian_connectors_affine_charts}
\end{figure}

\subsection{Canonical Divergences}
Following \cite{ay2015novel,ay2024information,ay2025correction}, we use the canonical-divergence functional
\begin{align}\label{eq:KL_energy_formula}
    D_g^{(\nabla)}(\mu \| \nu) = \int_0^1 t \lVert \dot{\gamma}(t) \rVert_{\gamma(t),g}^2 dt,
\end{align}
where $\gamma$ is the selected $\nabla$-geodesic from $\mu$ to $\nu$ and $g$ is the chosen Riemannian metric; we suppress the subscript $g$ when the metric is clear. Thus the construction is not restricted to the exponential connections used below. In the Fisher--Rao/$e_0$ case, this is the classical canonical divergence and equals the reversed KL divergence. In the Otto/$e_1$ case, the resulting canonical divergence is $D_{\rm WKL}$. Appendix~\ref{sec:app_mixture_canonical} instead evaluates the same functional along the intrinsic Gaussian mixture geodesic. With the Fisher--Rao metric it gives $D_{\rm KL}(\mu\Vert\nu)$ in the forward order, whereas with the Otto metric it yields a different canonical divergence that generally differs from $D_{\rm WKL}(\nu\Vert\mu)$.

For the $e_{0}$-geodesic given by \eqref{eq:e0-geodesic}, the canonical divergence \eqref{eq:KL_energy_formula} reduces to the KL-divergence (see \cite{ay2024information,ay2025correction} and Appendix~\ref{sec:app_reverse_KL} for a self-contained derivation in the Gaussian setting), i.e.,
\begin{align}\label{eq:KL_energy_formula_specialize_to_KL}
    D^{(\rm e_0)}(\mu \| \nu) 
    &:=\int_0^1 t \langle\dot{\gamma}_0(t),\dot{\gamma}_0(t) \rangle_{\gamma_0(t)}^{\rm FR} dt= D_{\rm KL}(\nu \| \mu).
\end{align}
Applying \eqref{eq:KL_energy_formula} to the $e_{1}$-geodesic \eqref{eq:e1-geodesic} yields the following identity.
\begin{align*}
     D^{(\rm e_1)}(\mu \| \nu):=\int_0^1 t \langle\dot{\gamma}_1(t),\dot{\gamma}_1(t) \rangle_{\gamma_1(t)}^{\rm O} dt=\int_{\Rn} \int_0^1 \left(f \circ \varphi_1 - f \circ \varphi_t\right) \, dt \,  d\mu.
\end{align*}
This reduction is proved in Appendix~\ref{sec:WKL_divergence_formula}.
This defines the canonical divergence that we call the {\em Wasserstein KL divergence\/}, abbreviated by ${\rm WKL}$.
\begin{equation}
\label{eq:contrast_func}
    D_{\rm WKL}(\mu \| \nu) \; := \; \int_{\Rn} \int_0^1 \left(f \circ \varphi_1 - f \circ \varphi_t\right) \, dt \,  d\mu.
\end{equation}

We verify that $D_{\rm WKL}$ is non-negative and vanishes exactly on the diagonal. Along the gradient flow,
\begin{align*}
\frac{d}{dt}f(\varphi_t(x))
&=\left\langle\nabla f(\varphi_t(x)),\frac{d}{dt}\varphi_t(x)\right\rangle
=\lVert\nabla f(\varphi_t(x))\rVert^2\geq0.
\end{align*}
Thus $t\mapsto f(\varphi_t(x))$ is non-decreasing and
\[
h(t,x):=f(\varphi_1(x))-f(\varphi_t(x))\geq0.
\]
It follows immediately that $D_{\rm WKL}(\mu\|\nu)\geq0$. If $D_{\rm WKL}(\mu\|\nu)=0$, then $h=0$ for $(dt\otimes\mu)$-almost every $(t,x)$. By Fubini's theorem, for $\mu$-almost every $x$, the continuous function $t\mapsto h(t,x)$ vanishes for almost every $t\in[0,1]$, and hence for every $t\in[0,1]$. For such $x$, the function $t\mapsto f(\varphi_t(x))$ is constant, so
\[
0=\frac{d}{dt}f(\varphi_t(x))
=\left\lVert\frac{d}{dt}\varphi_t(x)\right\rVert^2
\qquad\text{for every }t\in[0,1].
\]
Therefore $\varphi_t(x)=x$ for every $t\in[0,1]$ and for $\mu$-almost every $x$. In particular, $\nu=(\varphi_1)_*\mu=\mu$. Conversely, if $\mu=\nu$, the constant potential generates the constant $e_1$-geodesic and gives $D_{\rm WKL}(\mu\|\mu)=0$. Hence $D_{\rm WKL}$ has the usual separation properties of a divergence.

It can further be shown that for product distributions $\mu_1\otimes\mu_2$ and $\nu_1\otimes\nu_2$, we have the additivity property $D_{\rm WKL}(\mu_1 \otimes \mu_2 \| \nu_1 \otimes \nu_2)=D_{\rm WKL}(\mu_1\|\nu_1)+D_{\rm WKL}(\mu_2\|\nu_2)$ (see \cite[Theorem~3.3]{stein2026wklcontrol}).

\section{The Wasserstein KL Divergence on the Gaussian Manifold}\label{sec:gaussian_wkl}

With the Gaussian endpoint generator explicit by Lemma~\ref{lemm:formulae_for_A_b}, we now evaluate WKL divergence in endpoint coordinates and compare it with KL and other transport--information contrasts.
\subsection{Time-Integral Evaluation and Closed Formula}

The endpoint-generator and quadratic-flow calculations recalled in Section~\ref{sec:connections_and_geodesics}, as well as the time-integral evaluation and central closed formula below (including its commuting specialization), first appeared in \cite{datar2026wkl}. They are reproduced here with harmonized notation and the details needed for the boundary and comparison results.

For the quadratic potential $f$, coefficients $A,b$, and forward gradient flow $\varphi_t$ introduced in \eqref{eq:potential}--\eqref{eq:dynamics}, the following lemma evaluates the inner integral on the right-hand side of \eqref{eq:contrast_func}.
\begin{lemma}\label{lemm:time_integral}
Consider the gradient flow dynamics \eqref{eq:dynamics} with a quadratic function $f$ given in \eqref{eq:potential} and let $M=2Ae^{2A}-e^{2A}+I$. The following identity holds for all $x_0\in \Rn$:
\begin{align}
    \int_0^1 \left(f \circ \varphi_1 - f \circ \varphi_t\right)(x_0) dt
    &=\frac{1}{4}(x_0+\Adag{A}b)^\top M(x_0+\Adag{A}b)+  \frac{1}{2} b^\top\Perp{A}b.
\end{align}
\end{lemma}
\begin{proof}
By \eqref{eq:affine_flow_solution} and \eqref{eq:Kt_symmetric_formula}, setting
$y=x_0+\Adag{A}b$ rewrites the trajectory as
\begin{align*}
    x(t)=\expM{At}y+t\Perp{A}b-\Adag{A}b.
\end{align*}
Here we use the fact that $\Perp{A}$, $\Adag{A}$, and $e^A$ are simultaneously diagonalizable by an orthogonal matrix and therefore commute.

We now compute $f\circ \varphi_t$
\begin{align*}
    f(x(t))&=\frac{1}{2}x(t)^\top Ax(t) + b^\top x(t)=\frac{1}{2} \left(y^\top Ae^{2At}y-b^\top\Adag{A}b\right)+tb^\top\Perp{A}b+y^\top\Perp{A}b
\end{align*}
Therefore,
\begin{align*}
    (f\circ\varphi_1) (x_0) - (f \circ \varphi_t) (x_0)
    =\frac{1}{2} y^\top A\left(e^{2A}-e^{2At}\right)y + (1-t)b^\top\Perp{A}b.
\end{align*}
Integrating with respect to time, we get,
\begin{align*}
    \int_0^1 \left(f \circ \varphi_1 - f \circ \varphi_t\right)(x_0) dt&=\int_0^1\left(\frac{1}{2} y^\top A\left(e^{2A}-e^{2At}\right)y + (1-t) b^\top\Perp{A}b \right) dt    \\
    &=\frac{1}{4} y^\top \left(2Ae^{2A}-e^{2A}+I\right)y+  \frac{1}{2} b^\top\Perp{A}b.
\end{align*}
Plugging in the definitions of $y$ and $M$, we get the desired identity.
\end{proof}
We now present a main result, which provides a formula for the WKL-divergence between two nondegenerate Gaussian measures.
\begin{theorem}\label{theom:Wasserstein_KL_main_theorem}
Let $\mu=\mathcal{N}(m_0,\Sigma_0)$ and $\nu=\mathcal{N}(m_1,\Sigma_1)$ be two nondegenerate Gaussian measures with $\Sigma_0\succ 0$ and $\Sigma_1\succ 0$.
Then the following identity holds
\begin{align}\label{eq:D_main_formula}  
    D_{\rm WKL}(\mu \| \nu) &= \frac{1}{4}\trace{\Sigma_0-\Sigma_1+\Sigma_0R^2\log(R^2)} +\frac{1}{4}\left\lVert \sqrt{Q+2\Perp{\log(R)}} (m_1-m_0) \right\rVert^2
\end{align}
where
\begin{align*}
    R&=\Sigma_0^{-\frac{1}{2}}\left(\Sigma_0^{\frac{1}{2}}\Sigma_1\Sigma_0^{\frac{1}{2}}\right)^{\frac{1}{2}}\Sigma_0^{-\frac{1}{2}} \hspace{0.2cm} \textnormal{and}\hspace{0.2cm}
    Q=\Adag{(R-I)}\left(\log(R^2)R^2-R^2+I\right)\Adag{(R-I)}.
\end{align*}
Furthermore, if $\Sigma_0 \Sigma_1 = \Sigma_1 \Sigma_0$, then we get the following simplification:
\begin{align*}
  D_{\rm WKL}(\mu \| \nu) &= \frac{1}{4}\left(\left\lVert \sqrt{Q}\left(\sqrt{\Sigma_1}-\sqrt{\Sigma_0}\right) \right\rVert_F^2 +\left\lVert \sqrt{Q+2\Perp{\log(R)}} (m_1-m_0) \right\rVert^2\right).
\end{align*}
\end{theorem}
\begin{proof}
Note that the outer integral in the right-hand-side of \eqref{eq:contrast_func} corresponds to taking an expectation with respect to $\mu$.
Using Lemma \ref{lemm:formulae_for_A_b}, we define $A$ and $b$ according to \eqref{eq:soln_A} and \eqref{eq:soln_b}, so that $(\varphi_1)_*\mu=\nu$.
Using Lemma \ref{lemm:time_integral} along with the well-known properties of the expectation and trace operators\footnote{We mainly use the linearity and the cyclic property of the trace operators}, we get that
\begin{align}
    D_{\rm WKL}(\mu \| \nu)&=\int_{\Rn}\int_0^1 \left(f \circ \varphi_1 - f \circ \varphi_t\right)(x_0) dt d\mu(x_0)\nonumber \\
    &=\mathbb{E}\left[\frac{1}{4}(x_0+\Adag{A}b)^\top M(x_0+\Adag{A}b) +  \frac{1}{2} b^\top\Perp{A}b\right]\nonumber  \\
    &=\frac{1}{4}\left(\trace{(M\Sigma_0)}+(m_0+\Adag{A}b)^\top M(m_0+\Adag{A}b)\right) +  \frac{1}{2} b^\top\Perp{A}b \label{eq:D_of_mu_nu}
\end{align}
where $A$ and $b$ are given by \eqref{eq:soln_A} and \eqref{eq:soln_b}.
We now substitute the expression \eqref{eq:soln_b} for $b$ into \eqref{eq:D_of_mu_nu} term by term.
Let us first compute and simplify $m_0 + \Adag{A}b$ as
\begin{align}
    m_0 + \Adag{A} b & = m_0 + \Adag{A}K_1(A)^{-1}(m_1-e^Am_0) \nonumber \\
    & = m_0 + \Adag{(\expM{A}-I)}(m_1-e^Am_0) \nonumber 
     = \Perp{(\expM{A}-I)}m_0+ \Adag{(\expM{A}-I)}(m_1-m_0). \nonumber 
\end{align}
Using the fact that $M\Perp{(\expM{A}-I)}=0$, we get that
\begin{align}
    (m_0+\Adag{A}b)^\top M(m_0+\Adag{A}b)
    &=(m_1-m_0)^\top \Adag{(\expM{A}-I)}M \Adag{(\expM{A}-I)}(m_1-m_0).\label{eq:term2}
\end{align}
Since $\Perp{A}K_1(A)^{-1}=\Perp{A}$ and $\Perp{A}\expM{A}=\Perp{A}$, the last term can be computed as
\begin{align}
    b^\top\Perp{A}b=(m_1-e^Am_0)^\top\Perp{A}(m_1-e^Am_0)
    &=(m_1-m_0)^\top\Perp{A}(m_1-m_0).\label{eq:term3}
\end{align}
Finally, the cyclic property of the trace operator gives us
\begin{align}
    \trace{M\Sigma_0} &=\trace{\left(\log(R^2)R^2-R^2+I\right)\Sigma_0} \nonumber\\
    &=\trace{\log(R^2)R^2\Sigma_0-\Sigma_1+\Sigma_0} \label{eq:term1}.
\end{align}
Plugging in \eqref{eq:term2}, \eqref{eq:term3} and \eqref{eq:term1} in \eqref{eq:D_of_mu_nu} and using $M=2Ae^{2A}-e^{2A}+I=\log(R^2)R^2-R^2+I $ as well as $A=\log(R)$ gives the desired result. 
We now show that 
\begin{align}\label{eq:mat_in_mean_norm}
 Q+2\Perp{\log(R)}=\Adag{(\expM{A}-I)}M \Adag{(\expM{A}-I)}+2\Perp{A} \succ 0  
\end{align}
ensuring the square root is well defined.
Consider the orthogonal eigenvalue decomposition $A=U\Lambda U^\top$ where $\Lambda$ is a diagonal matrix consisting of eigenvalues $\lambda_i$, $i \in \{1,2,\cdots,n\}$ of $A$ and $U$ is an orthogonal matrix consisting of eigenvectors as columns.
Since $A$ and $e^{A}$ have the same eigenvectors, observe that
$U^\top MU=2\Lambda e^{2\Lambda}-e^{2\Lambda}+I$, which is a diagonal matrix consisting of entries $2\lambda_i e^{2\lambda_i}-e^{2\lambda_i}+1$, $i \in \{1,2,\cdots,n\}$ on the diagonal.
Observe that the function $h:\lambda\mapsto (e^{-2\lambda}-1)$ is strictly convex, which implies that $e^{-2\lambda}-1=h(\lambda)> h(0)+h'(0)\lambda=-2\lambda$ for all $\lambda \in \mathbb{R}\setminus \{0\}$ which implies that $2\lambda e^{2\lambda}-e^{2\lambda}+1> 0$ for all $\lambda \in \mathbb{R}\setminus \{0\}$.
Hence each diagonal entry of $U^\top M U$ is nonnegative, with equality if and only if $\lambda_i = 0$.
That directly implies $M \succeq 0$ and $M$ has an eigenvalue at $0$ if and only if $A$ has an eigenvalue at $0$.
Therefore, $Q=\Adag{(\expM{A}-I)}M \Adag{(\expM{A}-I)}\succeq 0$.
Furthermore, note that $z^\top Qz=0$ only if $z$ belongs to the null space of $A$ which implies that $z^\top\Perp{A}z=\lVert z\rVert^2$.
This shows that $z^\top\left(Q+2\Perp{\log(R)}\right)z=0$ if and only if $z=0$ thereby showing that $Q+2\Perp{\log(R)}\succ 0$.
This ensures the square root is well defined.

Finally, when $\Sigma_0 \Sigma_1 = \Sigma_1 \Sigma_0$, the first term can be reformulated using $M\Perp{(\expM{A}-I)}=0$ as
\begin{align*}
    \trace{M\Sigma_0} &= \trace{(\expM{A}-I)\Adag{(\expM{A}-I)}M\Adag{(\expM{A}-I)}(\expM{A}-I)\Sigma_0} \\
    &=\left\lVert \sqrt{Q}\left(\sqrt{\Sigma_1}-\sqrt{\Sigma_0}\right) \right\rVert_F^2.
\end{align*}
\end{proof}

Let us now show that the WKL-divergence \eqref{eq:D_main_formula} is the sum of two non-negative terms.
Furthermore, the first term is zero if and only if $\Sigma_1=\Sigma_0$.
This can be seen by observing that the first term $\frac{1}{4}\trace{\Sigma_0-\Sigma_1+\Sigma_0R^2\log(R^2)}=\frac{1}{4}\trace{M\Sigma_0}$ (see \eqref{eq:term1}) where $M=2Ae^{2A}-e^{2A}+I\succeq 0$.
Therefore, $\trace{M \Sigma_0}=\trace{\sqrt{\Sigma_0}M \sqrt{\Sigma_0}}=0$ implies $\sqrt{\Sigma_0}M \sqrt{\Sigma_0}=0$ and thus, $M=0$.
Recall from the proof of Theorem~\ref{theom:Wasserstein_KL_main_theorem} that $M$ has eigenvalues $2\lambda_i e^{2\lambda_i}-e^{2\lambda_i}+1$ where $\lambda_i$, $i\in\{1,2,\cdots,n\}$ are the eigenvalues of $A$.
Thus $M=0$ implies that $2\lambda_i e^{2\lambda_i}-e^{2\lambda_i}+1=0$ which in turn implies that $\lambda_i=0$ for all $i\in\{1,2,\cdots,n\}$ (see proof of Theorem~\ref{theom:Wasserstein_KL_main_theorem}). 
Therefore, 
\begin{align*}
A=\log \left(\Sigma_0^{-\frac{1}{2}}\left(\Sigma_0^{\frac{1}{2}}\Sigma_1\Sigma_0^{\frac{1}{2}}\right)^{\frac{1}{2}}\Sigma_0^{-\frac{1}{2}}\right)=0    
\end{align*}
which gives us $\Sigma_1=\Sigma_0$.
Similarly, owing to the fact that $Q+2\Perp{\log(R)}\succ 0$ (see proof of Theorem~\ref{theom:Wasserstein_KL_main_theorem}), we see that the second term $\frac{1}{4}\left\lVert \sqrt{Q+2\Perp{\log(R)}} (m_1-m_0) \right\rVert^2$ is zero if and only if $m_1=m_0$.
Finally, note that when $\Sigma_1=\Sigma_0$, we get that $R=I$ and $Q=0$.
This gives us that $D_{\rm WKL}(\mu \| \nu) = \frac{1}{2}\left\lVert (m_1-m_0) \right\rVert^2$ showing that the WKL-divergence indeed captures the geometry of the sample space.
The structurally matched classical comparator is the reverse KL divergence
\begin{align}  \label{eq:D_main_formula_KL}
    D_{\rm KL}(\nu \| \mu) &= \frac{1}{2}\left(\log \left(\frac{\det \Sigma_0}{\det \Sigma_1}\right)-n  + \trace{\Sigma_0^{-1}\Sigma_1}\right) +\frac{1}{2}\left\lVert \Sigma_0^{-\frac{1}{2}} (m_1-m_0) \right\rVert^2,
\end{align}
which is distinct from $D_{\rm WKL}(\mu \| \nu)$. The reversed order is dictated by the canonical-divergence construction, not chosen for numerical convenience. Indeed, for the same ordered endpoints $\mu\to\nu$,
\[
D^{(e_0)}(\mu\Vert\nu)=D_{\rm KL}(\nu\Vert\mu),
\qquad
D^{(e_1)}(\mu\Vert\nu)=D_{\rm WKL}(\mu\Vert\nu).
\]
Thus the endpoints and weighted-energy prescription are held fixed while the metric--connection pair changes from Fisher--Rao/$e_0$ to Otto/$e_1$.

There are two distinct ways to obtain the forward order $D_{\rm KL}(\mu\Vert\nu)$: reverse the ordered endpoints within the $e_0$ construction, or retain the endpoints and integrate along the intrinsic mixture geodesic. Appendix~\ref{sec:app_reverse_KL} derives the $e_0$ identity directly, while Appendix~\ref{sec:app_mixture_canonical} treats the mixture-geodesic construction.

\subsection{Gaussian Transport Examples}\label{subsec:gaussian_transport_examples}

\begin{figure}[ht!]
    \centering
    \begin{subfigure}[t]{0.32\textwidth}
        \centering
        \small\textbf{Isotropic scaling}
        \vspace{0.4em}
        \scriptsize
        \[
        \mu=\mathcal{N}\!\left(
        \begin{bmatrix}0\\0\end{bmatrix},
        \begin{bmatrix}1&0\\0&1\end{bmatrix}\right)
        \]
        \[\Downarrow\]
        \[
        \nu=\mathcal{N}\!\left(
        \begin{bmatrix}0\\0\end{bmatrix},
        \begin{bmatrix}3&0\\0&3\end{bmatrix}\right)
        \]
        \includegraphics[width=\linewidth]{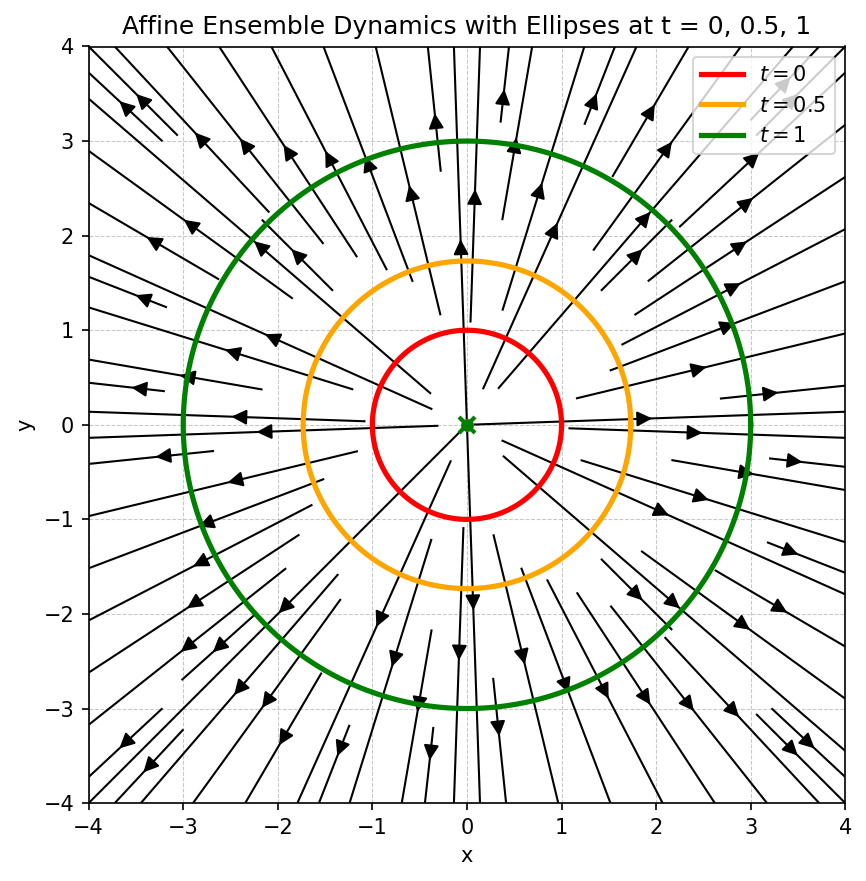}
        \[
        \begin{bmatrix}\dot{x}\\\dot{y}\end{bmatrix}
        =\begin{bmatrix}\log(\sqrt3)x\\\log(\sqrt3)y\end{bmatrix}.
        \]
    \end{subfigure}
    \hfill
    \begin{subfigure}[t]{0.32\textwidth}
        \centering
        \small\textbf{Anisotropic scaling}
        \vspace{0.4em}
        \scriptsize
        \[
        \mu=\mathcal{N}\!\left(
        \begin{bmatrix}0\\0\end{bmatrix},
        \begin{bmatrix}1&0\\0&1\end{bmatrix}\right)
        \]
        \[\Downarrow\]
        \[
        \nu=\mathcal{N}\!\left(
        \begin{bmatrix}0\\0\end{bmatrix},
        \begin{bmatrix}0.5&0\\0&3\end{bmatrix}\right)
        \]
        \includegraphics[width=\linewidth]{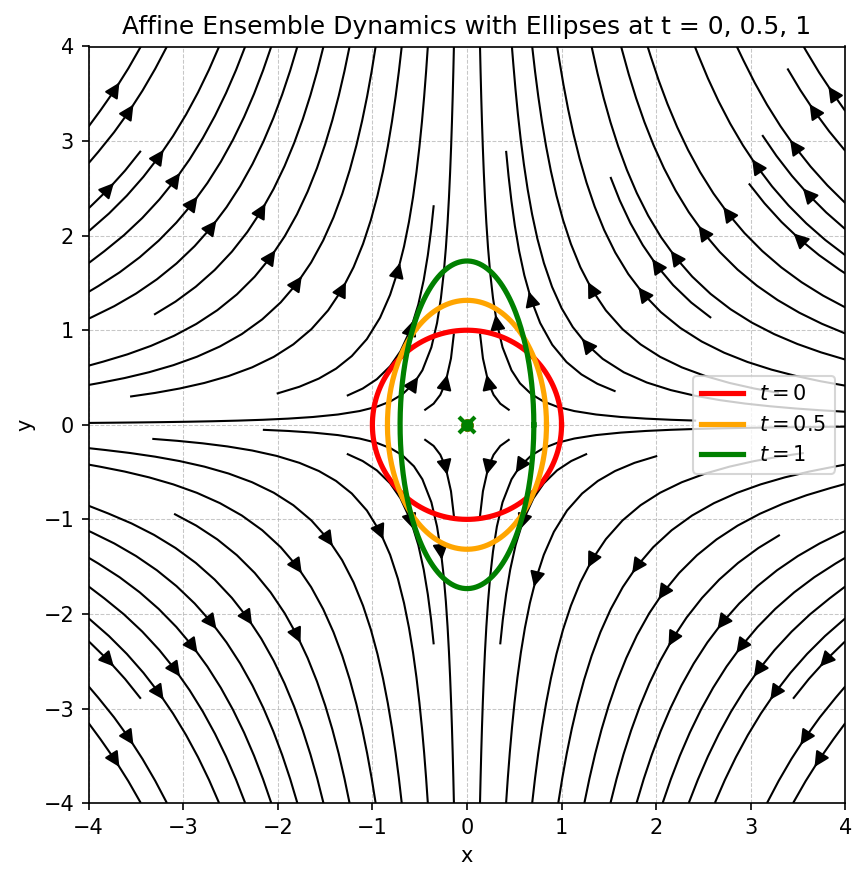}
        \[
        \begin{bmatrix}\dot{x}\\\dot{y}\end{bmatrix}
        =\begin{bmatrix}\log(\sqrt{0.5})x\\\log(\sqrt3)y\end{bmatrix}.
        \]
    \end{subfigure}
    \hfill
    \begin{subfigure}[t]{0.32\textwidth}
        \centering
        \small\textbf{Scaling and translation}
        \vspace{0.4em}
        \scriptsize
        \[
        \mu=\mathcal{N}\!\left(
        \begin{bmatrix}0\\0\end{bmatrix},
        \begin{bmatrix}1&0\\0&1\end{bmatrix}\right)
        \]
        \[\Downarrow\]
        \[
        \nu=\mathcal{N}\!\left(
        \begin{bmatrix}2\\2\end{bmatrix},
        \begin{bmatrix}1&0\\0&3\end{bmatrix}\right)
        \]
        \includegraphics[width=\linewidth]{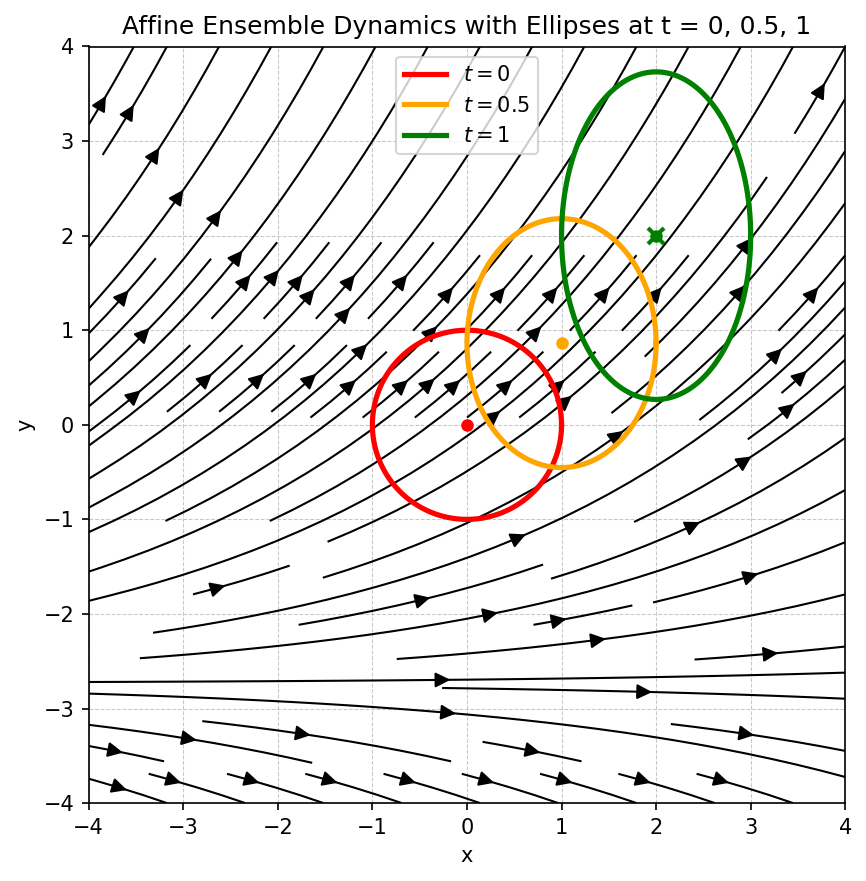}
        \[
        \begin{aligned}
        c&:=\frac{2\log(\sqrt3)}{\sqrt3-1},\\[-1mm]
        \begin{bmatrix}\dot{x}\\\dot{y}\end{bmatrix}
        &=\begin{bmatrix}2\\ \log(\sqrt3)y+c\end{bmatrix}.
        \end{aligned}
        \]
    \end{subfigure}
    \caption{Examples of Gaussian transport induced by affine vector fields. Left: isotropic expansion. Middle: contraction in one direction and expansion in the other. Right: combined translation and anisotropic scaling. The colors are also identified by the embedded labels $t=0,0.5,1$.}
    \label{fig:gaussian_transport_examples}
\end{figure}

To illustrate the affine flows underlying the Gaussian $e_1$-geodesics, we consider three examples in dimension two. In each case the vector field has the form $x\mapsto Ax+b$, so the mean and covariance evolve according to \eqref{eq:e1_global_geodesic}. Figure~\ref{fig:gaussian_transport_examples} displays the initial, intermediate, and terminal Gaussian ellipses together with the autonomous vector field.

\subsection{Relation to Optimal Transport and Related Divergences}
The $2$-Wasserstein distance $W_2$ is a Riemannian distance on the space of probability measures (more precisely, on the density manifold \cite{lafferty1988density}), arising as the geodesic distance associated with the Otto metric.
In this sense, it plays a role analogous to that of the Fisher$-$Rao distance in information geometry.
While our proposed WKL-divergence is not a distance, it is closely related to Wasserstein geometry, and therefore invites comparison with $W_2$, especially in the Gaussian case where both quantities admit closed-form expressions that have a close resemblance.

We review the case of zero-mean nondegenerate Gaussian measures $\mu = \mathcal{N}(0, \Sigma_0)$ and $\nu = \mathcal{N}(0, \Sigma_1)$ and point the reader to \cite{peyre2019computational,malago2018wasserstein} for further details.
The 2-Wasserstein distance between them arises as the minimal cost in an optimal transport problem. 
Specifically, the task is to find a linear transport map represented by a matrix $R \in \mathbb{R}^{n \times n}$ such that if $x \sim \mu$, then $Rx \sim \nu$.
This requirement leads to the algebraic condition $\Sigma_1 = R \Sigma_0 R^\top$, which generally admits multiple solutions.
To ensure uniqueness, we look at the linear maps $R$ that minimize the cost $\mathbb{E}_{x \sim \mu}\left[\lVert x-Rx\rVert ^2\right]$ subject to the constraint $\Sigma_1 = R \Sigma_0 R^\top$.
This leads to the map 
\begin{align}\label{eq:brenier_map}
 R=\Sigma_0^{-\frac{1}{2}}\left(\Sigma_0^{\frac{1}{2}}\Sigma_1\Sigma_0^{\frac{1}{2}}\right)^{\frac{1}{2}}\Sigma_0^{-\frac{1}{2}},   
\end{align}
which is the unique optimal transport map in this setting.
The square of the associated 2-Wasserstein distance is then given by the minimal cost under this map.
Equivalently, one may require $R$ to be symmetric positive semidefinite, so that the linear map is the gradient of a convex quadratic potential. In the present nondegenerate setting this apparently weaker condition already forces $R\succ0$: if $R$ were singular, then $R\Sigma_0R^\top$ would be singular and could not equal $\Sigma_1\succ0$. The unique positive-semidefinite solution is therefore the positive-definite matrix in \eqref{eq:brenier_map}, and its quadratic potential is strictly convex.
For general nondegenerate Gaussian measures (not necessarily zero mean) $\mu = \mathcal{N}(m_0, \Sigma_0)$ and $\nu = \mathcal{N}(m_1, \Sigma_1)$, the squared 2-Wasserstein distance is given by
\begin{align*}
    W_2(\mu,\nu)^2=\left \lVert m_0-m_1\right \rVert^2 + \trace{\Sigma_0 + \Sigma_1-2\left(\Sigma_0^{\frac{1}{2}}\Sigma_1\Sigma_0^{\frac{1}{2}}\right)^{\frac{1}{2}}}
\end{align*}
where the second term reduces to $\left \lVert \sqrt{\Sigma_0} - \sqrt{\Sigma_1}\right \rVert_F^2$ when $\Sigma_0$ and $\Sigma_1$ commute.

The WKL divergence as well as the 2-Wasserstein distance separate a mean contribution from a covariance contribution. When $\Sigma_0$ and $\Sigma_1$ commute, the resemblance is particularly explicit: the WKL divergence formula weights the covariance displacement by $Q$ and the mean displacement by $Q+2\Perp{\log(R)}$.

A more direct connection is obtained from the endpoint transport map. By \eqref{eq:soln_A} and \eqref{eq:brenier_map}, $e^A=R$, and the time-one affine flow map is
\[
\varphi_1(x)=m_1+R(x-m_0),
\]
which is exactly the optimal transport map between the two nondegenerate Gaussian measures. Although the generating potential $f$ need not be convex because $A$ may have negative eigenvalues, symmetry of $A$ implies $e^A=R\succ0$. Hence the time-one map itself is the gradient of a strictly convex quadratic function.

The endpoint maps therefore coincide, but the interpolating curves generally do not. For zero means, the $e_1$-geodesic is given by $\Sigma_t=R^t\Sigma_0R^t$, whereas the $W_2$ displacement geodesic is given by $\Sigma_t=((1-t)I+tR)\Sigma_0 ((1-t)I+tR)$ with means being constant at $m_0=m_1=0$ \cite[Corollary~1]{malago2018wasserstein}.


Several constructions combine optimal transport and information divergences, but they modify different ingredients. The closest terminological neighbor is Li's transport KL, the transport Bregman divergence generated by negative Boltzmann--Shannon entropy \cite{li2021transport}.
For the Gaussian map $R$ from $\mu$ to $\nu$ in \eqref{eq:brenier_map},
\[
D_{\rm TKL}(\nu\Vert\mu)=\operatorname{Tr}(R)-\log\det R-n.
\]
Although the Gaussian formula in \cite{li2021transport} is stated for centered measures, the same expression holds for arbitrary means: the affine Brenier map $x\mapsto m_1+R(x-m_0)$ has constant Jacobian $R$, on which Li's expression depends. It is therefore independent of $m_0$ and $m_1$ and vanishes between translated copies of the same Gaussian, in agreement with \cite[Theorem~2(i)]{li2021transport}. WKL divergence instead contains the mean displacement and separates nondegenerate Gaussian measures. 
A number of other constructions have appeared in the literature \cite{amari2018information,feydy2019sinkhorn,janati2020entropic,mallasto2022entropy,kainth2025bregman} some of which also have closed form expressions for Gaussian measures.
Each of these constructions differs from WKL, and a detailed comparison is left for future work.

\subsection{Univariate Gaussian Measures}
The univariate formula, equal-covariance comparison, and scalar local-curvature
calculation in this subsection appeared in \cite{datar2026wkl}. We retain them
to fix the argument orientation and to place the strengthened boundary results
below in context. The main formula is the following direct application of
Theorem~\ref{theom:Wasserstein_KL_main_theorem}.
\begin{corollary}\label{corr:univariate}
Let $\mu=\mathcal{N}(m_0,\sigma_0^2)$ and $\nu=\mathcal{N}(m_1,\sigma_1^2)$ be two nondegenerate univariate Gaussian measures with $\sigma_0>0$ and $\sigma_1>0$.
Then the following identity holds
\begin{align}\label{eq:D_main_formula_univariate}    
    D_{\rm WKL}\left(\mu \Vert \nu\right) & = \begin{cases}
        \frac{\sigma_0^2-\sigma_1^2+\sigma_1^2\log\left(\frac{\sigma_1^2}{\sigma_0^2}\right)}{4\left(\sigma_1-\sigma_0\right)^2}\left(\left(\sigma_1-\sigma_0\right)^2 + \left(m_1-m_0 \right)^2\right) \hspace{0.0cm}\textnormal{ if } \sigma_0\neq \sigma_1, \\
        \frac{1}{2}(m_1-m_0)^2 \hspace{5.00cm}\textnormal{ if } \sigma_0= \sigma_1.
    \end{cases}
\end{align}
Furthermore, we have that
\begin{align*}
   \lim_{\sigma_1 \rightarrow \sigma_0}
   \frac{\sigma_0^2-\sigma_1^2+\sigma_1^2\log\left(\frac{\sigma_1^2}{\sigma_0^2}\right)}
        {4\left(\sigma_1-\sigma_0\right)^2}
   \Bigl[\left(\sigma_1-\sigma_0\right)^2
        +\left(m_1-m_0\right)^2\Bigr]
   = \frac{1}{2}(m_1-m_0)^2.
\end{align*}
\end{corollary}
\begin{proof}
    Substituting $A\leftarrow a$, $\Sigma_0\leftarrow \sigma_0^2$, $\Sigma_1\leftarrow \sigma_1^2$, we get
\begin{align*}
    R=\frac{\sigma_1}{\sigma_0} \quad \textnormal{  and  } \quad
    Q=\left(\sigma_0^2-\sigma_1^2+\sigma_1^2\log\left(\frac{\sigma_1^2}{\sigma_0^2}\right) \right)\Adag{\left(\left(\sigma_1-\sigma_0\right)^2\right)}.
\end{align*}
This leads to
\begin{align*}    
    &D_{\rm WKL}\left(\mu \Vert \nu\right)
    = \frac{Q}{4} \left((\sigma_1-\sigma_0)^2 + (m_1-m_0)^2\right)  +  \frac{1}{2}\Perp{\log(R)} (m_1-m_0)^2
\end{align*}
which gives the desired result.
The final limit can be shown via successive applications of L'Hopital's rule.
\end{proof}
\begin{figure}[H]
	\centering
	\begin{minipage}{0.49\textwidth}
		\centering
%
%
\begin{tikzpicture}

\begin{axis}[%
width=0.44*4.521in,
height=0.45*3.566in,
at={(0.758in,0.481in)},
scale only axis,
xmin=0,
xmax=20,
xlabel style={font=\color{white!15!black}},
xlabel={$\sigma$},
ymin=0,
ymax=9,
ylabel style={font=\color{white!15!black}},
ylabel={$D_{\rm WKL}$ / $D_{\rm KL}$},
axis background/.style={fill=white},
xmajorgrids,
ymajorgrids,
legend style={legend cell align=left, align=left, draw=white!15!black}
]
\addplot [color=blue!70!black, line width=1.5pt, solid]
  table[row sep=crcr]{%
19.99	0.00500500375250157\\
19.89	0.00505545710052857\\
19.79	0.00510667721025293\\
19.69	0.00515867969783546\\
19.59	0.00521148057902676\\
19.49	0.0052650962814994\\
19.39	0.00531954365762777\\
19.29	0.00537483999773181\\
19.19	0.00543100304380562\\
19.09	0.00548805100375083\\
18.99	0.00554600256613536\\
18.89	0.00560487691550171\\
18.79	0.00566469374824574\\
18.69	0.00572547328909256\\
18.59	0.00578723630819467\\
18.49	0.00585000413887793\\
18.39	0.00591379869606656\\
18.29	0.00597864249541369\\
18.19	0.006044558673171\\
18.09	0.00611157100682935\\
17.99	0.00617970393656408\\
17.89	0.00624898258752249\\
17.79	0.00631943279299019\\
17.69	0.00639108111847753\\
17.59	0.0064639548867661\\
17.49	0.00653808220396135\\
17.39	0.00661349198659711\\
17.29	0.0066902139898396\\
17.19	0.00676827883684417\\
17.09	0.00684771804931794\\
16.99	0.0069285640793445\\
16.89	0.00701085034253257\\
16.79	0.00709461125255006\\
16.69	0.00717988225711086\\
16.59	0.00726669987548512\\
16.49	0.00735510173760601\\
16.39	0.00744512662485231\\
16.29	0.00753681451258859\\
16.19	0.00763020661454983\\
16.09	0.00772534542916414\\
15.99	0.00782227478790887\\
15.89	0.00792103990580306\\
15.79	0.008021687434147\\
15.69	0.00812426551562073\\
15.59	0.00822882384186507\\
15.49	0.00833541371367275\\
15.39	0.00844408810392649\\
15.29	0.00855490172342777\\
15.19	0.00866791108976883\\
15.09	0.00878317459941036\\
14.99	0.00890075260313639\\
14.89	0.00902070748506734\\
14.79	0.00914310374542671\\
14.69	0.00926800808726391\\
14.59	0.00939548950735225\\
14.49	0.00952561939149388\\
14.39	0.0096584716144763\\
14.29	0.00979412264494195\\
14.19	0.00993265165545021\\
14.09	0.0100741406380256\\
13.99	0.0102186745255086\\
13.89	0.0103663413190443\\
13.79	0.010517232222065\\
13.69	0.010671441781149\\
13.59	0.0108290680341614\\
13.49	0.0109902126661102\\
13.39	0.0111549811731805\\
13.29	0.0113234830354408\\
13.19	0.0114958318987494\\
13.09	0.0116721457664252\\
12.99	0.0118525472012877\\
12.89	0.0120371635387094\\
12.79	0.0122261271113757\\
12.69	0.0124195754864903\\
12.59	0.0126176517162214\\
12.49	0.0128205046022406\\
12.39	0.0130282889752664\\
12.29	0.0132411659905948\\
12.19	0.013459303440669\\
12.09	0.0136828760858217\\
11.99	0.013912066004406\\
11.89	0.0141470629636258\\
11.79	0.0143880648124768\\
11.69	0.0146352778983155\\
11.59	0.014888917508697\\
11.49	0.0151492083402451\\
11.39	0.0154163849964658\\
11.29	0.0156906925165595\\
11.19	0.0159723869374625\\
11.09	0.0162617358915212\\
10.99	0.0165590192424083\\
10.89	0.0168645297621005\\
10.79	0.0171785738519774\\
10.69	0.0175014723113582\\
10.59	0.0178335611570771\\
10.49	0.0181751924980076\\
10.39	0.0185267354687866\\
10.29	0.0188885772273646\\
10.19	0.0192611240214146\\
10.09	0.0196448023290877\\
9.99	0.0200400600801002\\
9.89	0.0204473679636773\\
9.79	0.0208672208304945\\
9.69	0.0213001391964096\\
9.59	0.0217466708565253\\
9.49	0.0222073926189289\\
9.39	0.0226829121683617\\
9.29	0.0231738700710626\\
9.19	0.0236809419331464\\
9.09	0.0242048407260969\\
8.99	0.024746319294334\\
8.89	0.0253061730613258\\
8.79	0.025885242952419\\
8.69	0.0264844185544539\\
8.59	0.0271046415343396\\
8.49	0.0277469093411358\\
8.39	0.0284122792188328\\
8.29	0.0291018725599899\\
8.19	0.0298168796337295\\
8.09	0.030558564725332\\
7.99	0.0313282717288977\\
7.89	0.0321274302393012\\
7.79	0.0329575621950396\\
7.69	0.0338202891296517\\
7.59	0.0347173400962713\\
7.49	0.0356505603376821\\
7.39	0.0366219207831232\\
7.29	0.0376335284631785\\
7.19	0.038687637945609\\
7.09	0.0397866639081247\\
6.99	0.0409331949791344\\
6.89	0.0421300089947569\\
6.79	0.0433800898401662\\
6.69	0.0446866460661228\\
6.59	0.0460531314978091\\
6.49	0.0474832680834092\\
6.39	0.0489810712650097\\
6.29	0.0505508781951315\\
6.19	0.0521973791695919\\
6.09	0.053925652702619\\
5.99	0.0557412047346579\\
5.89	0.0576500125388778\\
5.79	0.0596585739811062\\
5.69	0.0617739628923806\\
5.59	0.0640038914365995\\
5.49	0.0663567805017237\\
5.39	0.068841839316263\\
5.29	0.0714691556991293\\
5.19	0.0742497985974214\\
5.09	0.0771959348620703\\
4.99	0.0803209625664155\\
4.89	0.0836396636012731\\
4.79	0.087168378798907\\
4.69	0.0909252094689514\\
4.59	0.0949302499988135\\
4.49	0.0992058571138041\\
4.39	0.103776962552083\\
4.29	0.108671437342766\\
4.19	0.113920517654832\\
4.09	0.119559304403967\\
3.98999999999999	0.125627351586988\\
3.89	0.132169361820237\\
3.79	0.139236012002145\\
3.69	0.146884937684065\\
3.59	0.155181911996338\\
3.48999999999999	0.164202264349226\\
3.39	0.174032596305288\\
3.29	0.18477286795207\\
3.19	0.196538949106239\\
3.09	0.209465757585279\\
2.98999999999999	0.223711144170648\\
2.89	0.239460734426073\\
2.79	0.256934006500431\\
2.69	0.276391979104767\\
2.59	0.298147016293735\\
2.48999999999999	0.322575442331577\\
2.39	0.350133926226783\\
2.29	0.381380980530502\\
2.19	0.41700548362211\\
2.09	0.457864975618692\\
1.98999999999999	0.505037751571932\\
1.89	0.559894739788921\\
1.79	0.624200243438098\\
1.69	0.700255593291553\\
1.59	0.791107946679328\\
1.48999999999999	0.900860321607141\\
1.39	1.03514310853476\\
1.29	1.20185085030948\\
1.19	1.41232963773745\\
1.09	1.68335998653313\\
0.989999999999995	2.04060810121418\\
0.889999999999997	2.52493372048985\\
0.789999999999996	3.20461464508896\\
0.689999999999998	4.20079815164884\\
0.589999999999996	5.74547543809257\\
0.489999999999995	8.32986255726798\\
0.389999999999997	13.1492439184749\\
0.289999999999996	23.7812128418557\\
0.189999999999998	55.4016620498628\\
0.0899999999999963	246.913580246934\\
};
\addlegendentry{$D_{\rm KL}(\mathcal{N}(2,\sigma^2)\Vert\mathcal{N}(0,\sigma^2))$}

\addplot [color=red!70!black, line width=1.5pt, dashed]
  table[row sep=crcr]{%
19.99	2\\
19.89	2\\
19.79	2\\
19.69	2\\
19.59	2\\
19.49	2\\
19.39	2\\
19.29	2\\
19.19	2\\
19.09	2\\
18.99	2\\
18.89	2\\
18.79	2\\
18.69	2\\
18.59	2\\
18.49	2\\
18.39	2\\
18.29	2\\
18.19	2\\
18.09	2\\
17.99	2\\
17.89	2\\
17.79	2\\
17.69	2\\
17.59	2\\
17.49	2\\
17.39	2\\
17.29	2\\
17.19	2\\
17.09	2\\
16.99	2\\
16.89	2\\
16.79	2\\
16.69	2\\
16.59	2\\
16.49	2\\
16.39	2\\
16.29	2\\
16.19	2\\
16.09	2\\
15.99	2\\
15.89	2\\
15.79	2\\
15.69	2\\
15.59	2\\
15.49	2\\
15.39	2\\
15.29	2\\
15.19	2\\
15.09	2\\
14.99	2\\
14.89	2\\
14.79	2\\
14.69	2\\
14.59	2\\
14.49	2\\
14.39	2\\
14.29	2\\
14.19	2\\
14.09	2\\
13.99	2\\
13.89	2\\
13.79	2\\
13.69	2\\
13.59	2\\
13.49	2\\
13.39	2\\
13.29	2\\
13.19	2\\
13.09	2\\
12.99	2\\
12.89	2\\
12.79	2\\
12.69	2\\
12.59	2\\
12.49	2\\
12.39	2\\
12.29	2\\
12.19	2\\
12.09	2\\
11.99	2\\
11.89	2\\
11.79	2\\
11.69	2\\
11.59	2\\
11.49	2\\
11.39	2\\
11.29	2\\
11.19	2\\
11.09	2\\
10.99	2\\
10.89	2\\
10.79	2\\
10.69	2\\
10.59	2\\
10.49	2\\
10.39	2\\
10.29	2\\
10.19	2\\
10.09	2\\
9.99	2\\
9.89	2\\
9.79	2\\
9.69	2\\
9.59	2\\
9.49	2\\
9.39	2\\
9.29	2\\
9.19	2\\
9.09	2\\
8.99	2\\
8.89	2\\
8.79	2\\
8.69	2\\
8.59	2\\
8.49	2\\
8.39	2\\
8.29	2\\
8.19	2\\
8.09	2\\
7.99	2\\
7.89	2\\
7.79	2\\
7.69	2\\
7.59	2\\
7.49	2\\
7.39	2\\
7.29	2\\
7.19	2\\
7.09	2\\
6.99	2\\
6.89	2\\
6.79	2\\
6.69	2\\
6.59	2\\
6.49	2\\
6.39	2\\
6.29	2\\
6.19	2\\
6.09	2\\
5.99	2\\
5.89	2\\
5.79	2\\
5.69	2\\
5.59	2\\
5.49	2\\
5.39	2\\
5.29	2\\
5.19	2\\
5.09	2\\
4.99	2\\
4.89	2\\
4.79	2\\
4.69	2\\
4.59	2\\
4.49	2\\
4.39	2\\
4.29	2\\
4.19	2\\
4.09	2\\
3.98999999999999	2\\
3.89	2\\
3.79	2\\
3.69	2\\
3.59	2\\
3.48999999999999	2\\
3.39	2\\
3.29	2\\
3.19	2\\
3.09	2\\
2.98999999999999	2\\
2.89	2\\
2.79	2\\
2.69	2\\
2.59	2\\
2.48999999999999	2\\
2.39	2\\
2.29	2\\
2.19	2\\
2.09	2\\
1.98999999999999	2\\
1.89	2\\
1.79	2\\
1.69	2\\
1.59	2\\
1.48999999999999	2\\
1.39	2\\
1.29	2\\
1.19	2\\
1.09	2\\
0.989999999999995	2\\
0.889999999999997	2\\
0.789999999999996	2\\
0.689999999999998	2\\
0.589999999999996	2\\
0.489999999999995	2\\
0.389999999999997	2\\
0.289999999999996	2\\
0.189999999999998	2\\
0.0899999999999963	2\\
};
\addlegendentry{$D_{\rm WKL}\left(\mathcal{N}(0,\sigma^2)\Vert\mathcal{N}(2,\sigma^2)\right)$}

\end{axis}
\end{tikzpicture}%
	\end{minipage}
	\hfill
	\begin{minipage}{0.49\textwidth}
		\centering
%
%
\begin{tikzpicture}

\begin{axis}[%
width=0.46*4.521in,
height=0.45*3.566in,
at={(0.758in,0.481in)},
scale only axis,
axis y line=right,
unbounded coords=jump,
xmin=0,
xmax=4,
xlabel style={font=\color{white!15!black}},
xlabel={$\sigma$},
ymin=0,
ymax=0.85,
ylabel style={font=\color{white!15!black}},
ylabel={},
axis background/.style={fill=white},
xmajorgrids,
ymajorgrids,
legend style={legend cell align=left, align=left, draw=white!15!black}
]
\addplot [color=blue!70!black, line width=1.5pt, solid]
  table[row sep=crcr]{%
2	0.127687330330387\\
2.1	0.101674943938732\\
2.2	0.0790438171927283\\
2.3	0.0595920546218945\\
2.4	0.0431435513142098\\
2.5	0.0295437790161768\\
2.6	0.0186563991962289\\
2.7	0.0103605156578261\\
2.8	0.00454842704250691\\
2.9	0.00112377389790363\\
3	0\\
3.1	0.00109906606589805\\
3.2	0.00435036775131781\\
3.3	0.00968982019567521\\
3.4	0.0170590792682161\\
3.5	0.0264048757282972\\
3.6	0.0376784432060453\\
3.7	0.0508350245734865\\
3.8	0.0658334441579917\\
3.9	0.0826357355325089\\
4	0.101206816437108\\
};
\addlegendentry{$D_{\rm KL}\left(\mathcal{N}(0,\sigma^2)\Vert\mathcal{N}(0,\sigma_{\rm opt}^2)\right)$}

\addplot [color=red!70!black, line width=1.5pt, dashed]
  table[row sep=crcr]{%
2	0.439069783783671\\
2.1	0.361031748615095\\
2.2	0.289425073504708\\
2.3	0.2247151266362\\
2.4	0.167346572215076\\
2.5	0.117745135018892\\
2.6	0.076319148494524\\
2.7	0.043460920427223\\
2.8	0.0195479437711502\\
2.9	0.0049439752037598\\
3	0\\
3.1	0.00505509866447121\\
3.2	0.0204372282243645\\
3.3	0.0464639290345491\\
3.4	0.0834429662741547\\
3.5	0.131672913941957\\
3.6	0.191443688024826\\
3.7	0.263037034572263\\
3.8	0.346726977623743\\
3.9	0.44278023127527\\
4	0.551456579614247\\
};
\addlegendentry{$D_{\rm WKL}\left(\mathcal{N}(0,\sigma_{\rm opt}^2)\Vert\mathcal{N}(0,\sigma^2)\right)$}

\addplot [color=blue!70!black, line width=1.5pt, solid, forget plot]
  table[row sep=crcr]{%
0	inf\\
0.1	1.80758509299405\\
0.2	1.1294379124341\\
0.3	0.748972804325936\\
0.4	0.496290731874155\\
0.5	0.318147180559945\\
0.6	0.190825623765991\\
0.7	0.101674943938732\\
0.8	0.0431435513142098\\
0.9	0.0103605156578264\\
1	0\\
1.1	0.00968982019567521\\
1.2	0.0376784432060454\\
1.3	0.0826357355325089\\
1.4	0.143527763378787\\
1.5	0.219534891891836\\
1.6	0.309996370754265\\
1.7	0.414371748937829\\
1.8	0.532213335097881\\
1.9	0.663146113827605\\
2	0.806852819440055\\
};
\addplot [color=red!70!black, line width=1.5pt, dashed, forget plot]
  table[row sep=crcr]{%
0	0.25\\
0.1	0.23598707453503\\
0.2	0.207811241751318\\
0.3	0.173321223805333\\
0.4	0.136696741450068\\
0.5	0.100856602430007\\
0.6	0.0680513877221217\\
0.7	0.0401146387350106\\
0.8	0.0185940635794529\\
0.9	0.00482899115858035\\
1	0\\
1.1	0.00516265878161655\\
1.2	0.0212715208916473\\
1.3	0.04919780347503\\
1.4	0.0897427918887886\\
1.5	0.143648246621685\\
1.6	0.211604645434542\\
1.7	0.294257822784836\\
1.8	0.392214397141433\\
1.9	0.506046264541173\\
2	0.636294361119891\\
};
\end{axis}
\end{tikzpicture}%
	\end{minipage}
	\caption{Left: comparison of $D_{\rm WKL}\!\left(\mathcal{N}(0,\sigma^2)\Vert\mathcal{N}(2,\sigma^2)\right)$ (dashed) and $D_{\rm KL}\!\left(\mathcal{N}(2,\sigma^2)\Vert\mathcal{N}(0,\sigma^2)\right)$ (solid) as the common standard deviation decreases. Right: for $\mu=\mathcal{N}(0,\sigma_{\rm opt}^2)$ and $\nu_\sigma=\mathcal{N}(0,\sigma^2)$, local variation of $D_{\rm WKL}(\mu\Vert\nu_\sigma)$ (dashed) and $D_{\rm KL}(\nu_\sigma\Vert\mu)$ (solid) around $\sigma=\sigma_{\rm opt}\in\{1,3\}$. In the right panel, the plotted WKL divergence value at $\sigma=0$ is its continuous target-boundary limit; the infinite KL endpoint is omitted, and finite KL values above the displayed ordinate ranges are clipped.}
	\label{fig:zero_variance_limit}
\end{figure}
For the reversed argument order dictated by the $e_0$ canonical-divergence identity, the univariate KL divergence is
\begin{align*}
D_{\rm KL}\left(\nu\Vert\mu\right)
&=\log\left(\frac{\sigma_0}{\sigma_1}\right)
+\frac{\sigma_1^2}{2\sigma_0^2}
+\frac{(m_1-m_0)^2}{2\sigma_0^2}
-\frac12.
\end{align*}
Therefore, if we consider two nondegenerate univariate Gaussian measures $\mu=\mathcal{N}(m_0,\sigma^2)$ and $\nu=\mathcal{N}(m_1,\sigma^2)$ with equal variance $\sigma^2$ and possibly different means $m_0$ and $m_1$, the WKL-divergence and the classical KL-divergence between these measures can be computed as
\begin{align*}
      D_{\rm WKL}(\mu\lVert\nu)= \frac{(m_1-m_0)^2}{2} \quad\quad \textnormal{and}\quad\quad
      D_{\rm KL}(\nu\lVert\mu)= \frac{(m_1-m_0)^2}{2\sigma^2}. 
\end{align*}
The argument reversal is required because the $e_0$ canonical divergence equals $D_{\rm KL}(\nu\Vert\mu)$ \cite{ay2024information,ay2025correction}. If $m_1\neq m_0$, then along the equal-covariance path $\sigma_0=\sigma_1=\sigma\downarrow0$, the KL divergence tends to infinity, whereas
\[
D_{\rm WKL}(\mu\Vert\nu)=\frac12(m_1-m_0)^2
\]
remains constant.
Figure~\ref{fig:zero_variance_limit} (left) displays this contrast for $m_0=0$ and $m_1=2$.

Thus, along this particular equal-covariance approximation, the WKL-divergence has a finite limit proportional to the squared Euclidean distance between the support points of the limiting Dirac measures. This is a pathwise statement, not a path-independent continuous extension; the continuity analysis below shows that other covariance paths can yield different limiting values.

To compare the local standard-deviation curvature, fix
\[
\mu=\mathcal{N}(0,\sigma_{\rm opt}^2),
\qquad
\nu_\sigma=\mathcal{N}(0,\sigma^2).
\]
At the common minimizer $\sigma=\sigma_{\rm opt}$,
\begin{align*}
\left.\frac{d^2}{d\sigma^2}D_{\rm WKL}(\mu\Vert\nu_\sigma)\right|_{\sigma=\sigma_{\rm opt}}
&=1,
&
\left.\frac{d^2}{d\sigma^2}D_{\rm KL}(\nu_\sigma\Vert\mu)\right|_{\sigma=\sigma_{\rm opt}}
&=\frac{2}{\sigma_{\rm opt}^2}.
\end{align*}
Thus the local WKL divergence curvature in the standard-deviation coordinate is independent of $\sigma_{\rm opt}$, whereas the corresponding KL curvature diverges as $\sigma_{\rm opt}\downarrow0$. Figure~\ref{fig:zero_variance_limit} (right) illustrates these local profiles.

\begin{figure}[t!]
	\centering
     \includegraphics[width=0.78\linewidth]{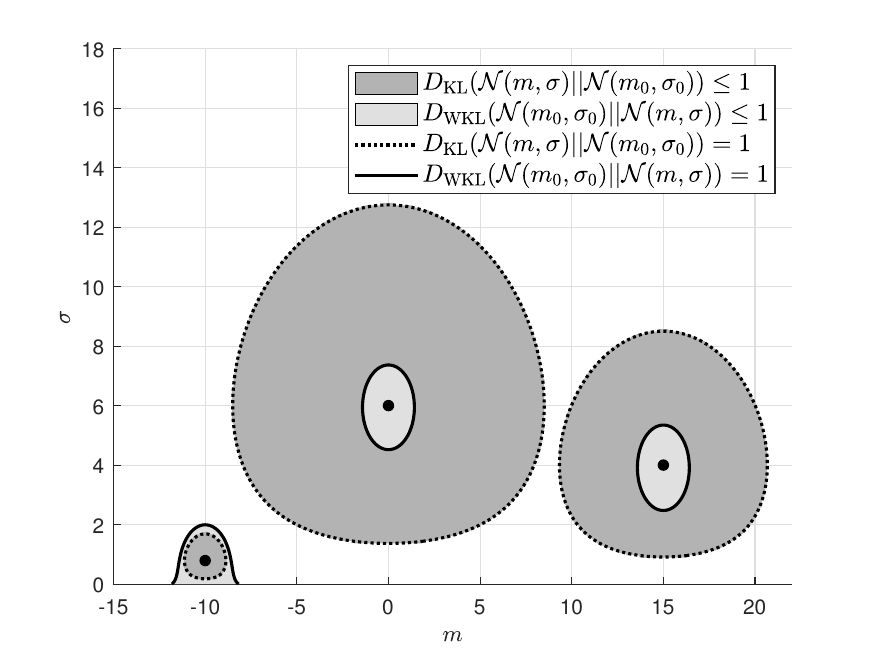}
	\caption{Level-one sublevel sets in the $(m,\sigma)$-plane for several reference measures $\mu_0=\mathcal{N}(m_0,\sigma_0^2)$, indicated by solid markers. The WKL divergence sets are defined by $D_{\rm WKL}(\mu_0\Vert\mathcal{N}(m,\sigma^2))\leq1$, while the KL sets use the reversed order $D_{\rm KL}(\mathcal{N}(m,\sigma^2)\Vert\mu_0)\leq1$. Displayed WKL divergence points with $\sigma=0$ use the continuous target-boundary extension proved in Theorem~\ref{theom:WKL_continuity}; reverse KL is infinite there and is omitted. Fill shades and solid/dotted boundaries are identified there.}
	\label{fig:divergence_balls}
\end{figure}

For the reference parameters and level displayed in Figure~\ref{fig:divergence_balls}, the KL sublevel sets vary substantially with $\sigma_0$, whereas the WKL sublevel sets are less sensitive to this change. The KL are seen to be convex, while the WKL sets become visibly non-convex for the smaller values of $\sigma_0$. 
These are observations about the plotted configurations, not a universal convexity threshold.
Furthermore, these observations are dependent on the chosen coordinate system. 
Related recent work studies non-accelerated and accelerated gradient flows of the classical KL functional on the Gaussian family under generalized bilinear-kernel Stein metrics \cite{stein2025accelerated}.
These metric-dependent optimization flows are distinct from both the coordinate-space experiment above and the $e_1$-connectors defining WKL. 
Detailed investigation of the convexity properties of WKL sublevel sets is left for future work.

\FloatBarrier

\begin{figure}[t!]
	\centering
     \includegraphics[width=\linewidth]{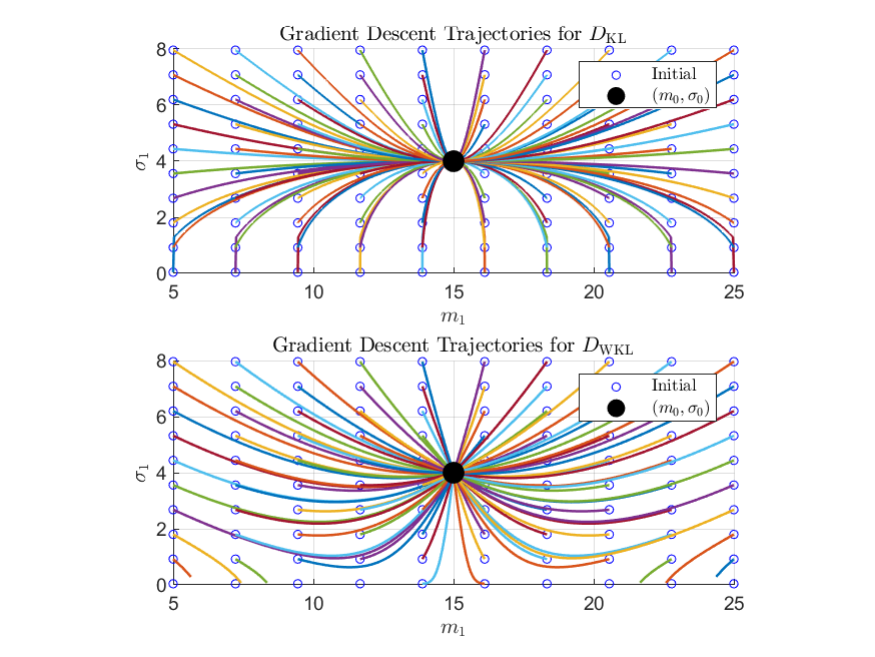}
	\caption{Archived coordinate-space gradient-descent trajectories for $D_{\rm KL}(\nu\Vert\mu)$ (top) and $D_{\rm WKL}(\mu\Vert\nu)$ (bottom), with reference $(m_0,\sigma_0)=(15,4)$. Open circles mark the displayed starting grid and the filled marker the reference. The axis includes $\sigma_1=0$, but reverse KL and its gradient are defined only for $\sigma_1>0$; boundary marks must therefore not be read as KL evaluations.}
    \label{fig:gradient_descent}
\end{figure}

Figure~\ref{fig:gradient_descent} depicts coordinate-space gradient-descent trajectories in the $(m_1,\sigma_1)$-plane, initialized from a grid around the reference parameters $(m_0,\sigma_0)$. In this finite experiment, all displayed KL trajectories move toward the reference point. Most WKL divergence trajectories do so as well, while some move toward smaller values of $\sigma_1$ and approach the covariance boundary. The plot therefore suggests that the chosen Euclidean-coordinate optimization scheme can exhibit boundary-attracting behavior for WKL.
Detailed investigations of the convergence properties of the WKL-divergence gradient flow are left for future work. 
\FloatBarrier
\subsection{Singular-Covariance Limits and Boundary Extensions}\label{sec:continuity}
In this subsection we collect the main continuity and boundary behaviour of the WKL-divergence \(D_{\mathrm{WKL}}\) on the Gaussian family.
In particular, we (i) derive a uniform lower bound that links \(D_{\mathrm{WKL}}\) to the Euclidean distance between the means, (ii) describe how the divergence behaves when one or both covariance matrices approach singularity, and (iii) contrast these phenomena with the classical KL-divergence.
We write
\[
\mathrm{cl}(S_n^+):=\{A\in S_n\mid A\succeq0\}
\]
for the positive-semidefinite cone.

We start by deriving a uniform lower bound on $D_{\rm WKL}$ in the next proposition.
\begin{proposition} \label{prop:continuity} The Wasserstein KL-divergence satisfies the inequality
\begin{align}\label{eq:uniform_lower_bound}  
    D_{\rm WKL}\left(\mathcal{N}(m_0,\Sigma_0) \| \mathcal{N}(m_1,\Sigma_1)\right) &\geq  \frac{1}{4}\left\lVert m_1-m_0 \right\rVert^2
\end{align}
for all $\Sigma_0\succ 0$ and $\Sigma_1\succ 0$.
Moreover, equality is attained in the sequential limit
\begin{align}
    \lim_{\Sigma_0 \rightarrow 0}\left(\lim_{\Sigma_1 \rightarrow 0}  D_{\rm WKL}\left(\mathcal{N}(m_0,\Sigma_0) \| \mathcal{N}(m_1,\Sigma_1)\right)\right) &=  \frac{1}{4}\left\lVert m_1-m_0 \right\rVert^2. \label{eq:limit1}
\end{align}
\end{proposition}
\begin{proof}
    Recall that we have already seen in the discussion after Theorem~\ref{theom:Wasserstein_KL_main_theorem} that the WKL-divergence is the sum of two non-negative terms.
    The second term thus already provides a lower bound on the WKL-divergence, that is,
    \begin{align*}
        D_{\rm WKL}(\mu \| \nu) \geq \frac{1}{4}\left\lVert \sqrt{Q+2\Perp{\log(R)}} (m_1-m_0) \right\rVert^2.
    \end{align*}
    Note that if $Q+2\Perp{\log(R)}-I\succeq 0$, then we have that
    \begin{align*}
       \frac{1}{4}\left\lVert \sqrt{Q+2\Perp{\log(R)}} (m_1-m_0) \right\rVert^2&= \frac{1}{4}(m_1-m_0)^\top\left(Q+2\Perp{\log(R)}\right)(m_1-m_0)\\
       &\geq \frac{1}{4}(m_1-m_0)^\top(m_1-m_0)
       =\frac{1}{4}\left\lVert m_1-m_0 \right\rVert^2
    \end{align*}
    which gives us the first statement.
    In order to see that $Q+2\Perp{\log(R)}-I\succeq 0$ is true, consider the orthogonal eigenvalue decomposition $R=P\Lambda P^\top$ and multiply the matrix $Q+2\Perp{\log(R)}$ from the left and right by $P^\top$ and $P$, respectively, to diagonalize it.
    This gives us
    \begin{align*}
        P^\top \left(Q+2\Perp{\log(R)}\right) P&=\Adag{(\Lambda-I)}\left(\log(\Lambda^2)\Lambda^2-\Lambda^2+I\right)\Adag{(\Lambda-I)}+2\Perp{\log(\Lambda)}
    \end{align*}
    which is a diagonal matrix with the $i^{\textnormal{th}}$ diagonal entry given by
    \begin{align*}
        g(\lambda_i):=\begin{cases}
            \frac{\log(\lambda_i^2)\lambda_i^2 - \lambda_i^2+1}{(\lambda_i-1)^2} \hspace{1cm} \textnormal{ if } \lambda_i \neq 1,\\
            \hspace{1cm} 2 \hspace{2.1cm} \textnormal{otherwise},
        \end{cases}
    \end{align*}
    where $\lambda_i$ is the $i^{\textnormal{th}}$ diagonal entry of $\Lambda$.
    We next show that $g(\lambda_i)\geq 1$ which implies that the eigenvalues of $\left(Q+2\Perp{\log(R)}\right)$ are greater than or equal to $1$ proving $Q+2\Perp{\log(R)}-I\succeq 0$.   
    Observe that convexity of the function $h:t\mapsto t\log t$ implies that
    \begin{align*}
        \lambda_i\log \lambda_i =h(\lambda_i)\geq h(1)+h'(1)(\lambda_i-1)=\lambda_i-1.
    \end{align*}
    Multiplying both sides by $2\lambda_i$ (note that $\lambda_i> 0$ since $R\succ 0$), then adding $1$ and subtracting $\lambda_i^2$ from both sides gives
        \begin{align*}
        \lambda_i^2 \log (\lambda_i^2)-\lambda_i^2+1=2\lambda_i^2\log \lambda_i+1-\lambda_i^2 \geq \lambda_i^2-2\lambda_i+1=(\lambda_i-1)^2.
    \end{align*}
    This completes the proof of the first statement.
    \\
    Let $\Sigma_0\succ 0$ be fixed and consider a sequence of symmetric positive definite matrices $ \Sigma_1^{(1)}, \Sigma_1^{(2)}, \cdots$ that converges to the zero matrix.
    Hence
    \begin{align*}
        R^{(i)}&=\Sigma_0^{-\frac{1}{2}}\left(\Sigma_0^{\frac{1}{2}}\Sigma_1^{(i)}\Sigma_0^{\frac{1}{2}}\right)^{\frac{1}{2}}\Sigma_0^{-\frac{1}{2}} \rightarrow 0,\\
        Q^{(i)}&=\Adag{\left(R^{(i)}-I\right)}\left(\log\left(\left(R^{(i)}\right)^2\right)\left(R^{(i)}\right)^2-\left(R^{(i)}\right)^2+I\right)\Adag{\left(R^{(i)}-I\right)} \rightarrow I.
    \end{align*}
    Furthermore, note that there exists a natural number $N$, such that for $i > N$, $\lVert R^{(i)}\rVert <1$ which implies that the largest eigenvalue of $R^{(i)}$ is less than $1$. 
    Furthermore, since $R^{(i)}\succ 0$, we get that $\Perp{\log\left(R^{(i)}\right)}=0$ for all $i>N$.
    This gives us that
    \begin{align*}
        \lim_{\Sigma_1 \rightarrow 0}  D_{\rm WKL}\left(\mathcal{N}(m_0,\Sigma_0) \| \mathcal{N}(m_1,\Sigma_1)\right) &= \frac{1}{4}\trace{\Sigma_0}+ \frac{1}{4}\left\lVert m_1-m_0 \right\rVert^2.
    \end{align*}
    Taking the limit as $\Sigma_0$ converges to the zero matrix gives us the desired result \eqref{eq:limit1}.
\end{proof}
When $\Sigma_0=\Sigma_1=\Sigma\succ0$, we have $R=I$, $Q=0$, and $\Perp{\log(R)}=I$. Hence
\begin{align}\label{eq:limit2}
D_{\rm WKL}\!\left(
\mathcal{N}(m_0,\Sigma)\,\middle\|\,\mathcal{N}(m_1,\Sigma)
\right)
=\frac12\lVert m_1-m_0\rVert^2,
\end{align}
independently of $\Sigma$. Together with \eqref{eq:limit1}, this gives
\begin{align*}
\lim_{\Sigma\to0}
D_{\rm WKL}\!\left(
\mathcal{N}(m_0,\Sigma)\,\middle\|\,\mathcal{N}(m_1,\Sigma)
\right)
&=\frac12\lVert m_1-m_0\rVert^2,\\
\lim_{\Sigma_0\to0}\lim_{\Sigma_1\to0}
D_{\rm WKL}\!\left(
\mathcal{N}(m_0,\Sigma_0)\,\middle\|\,\mathcal{N}(m_1,\Sigma_1)
\right)
&=\frac14\lVert m_1-m_0\rVert^2.
\end{align*}
Thus the two limits differ whenever $m_0\neq m_1$.

For comparison, fix $\Sigma_1\succ0$ and let $\Sigma_0=\sigma_0^2I$ with $\sigma_0\downarrow0$. In this case
\[
R=\frac{1}{\sigma_0}\Sigma_1^{1/2}.
\]
The covariance contribution to \eqref{eq:D_main_formula} is
\begin{align*}
\frac14\trace{\sigma_0^2I-\Sigma_1
    +\Sigma_1\log(\sigma_0^{-2}\Sigma_1)}
&=\frac14\bigl(n\sigma_0^2-\trace{\Sigma_1}
    +\trace{\Sigma_1\log\Sigma_1}\bigr)\\
&\quad-\frac12\log(\sigma_0)\trace{\Sigma_1},
\end{align*}
which tends to $+\infty$. The theorem below proves the corresponding statement for an arbitrary approach of $\Sigma_0$ to a singular positive-semidefinite matrix.

When $m_0\neq m_1$, these limits already show that no continuous value can be assigned at the corresponding pair of Dirac measures. The theorem below also treats the case $m_0=m_1$. The value $\frac14\lVert m_1-m_0\rVert^2$ will arise as the interior lower limit at such a pair.

For the boundary analysis, it is convenient to place the two scalar profiles governing the covariance and mean terms side by side. For $r\geq0$, define
\begin{align}\label{eq:boundary_profiles}
\chi(r)&:=
\begin{cases}
r^2\log(r^2),&r>0,\\
0,&r=0,
\end{cases}
&
\kappa(r)&:=
\begin{cases}
\dfrac{r^2\log(r^2)-r^2+1}{(r-1)^2},&r>0,\ r\neq1,\\[6pt]
2,&r=1,\\
1,&r=0.
\end{cases}
\end{align}
Thus $\chi(R)$ appears in the covariance contribution and $\kappa(R)$ in the quadratic mean contribution of the boundary formula below. The scalar limits at $0$ and $1$ show that both profiles are continuous on $[0,\infty)$.

We record explicitly the spectral-calculus fact used to pass from these scalar profiles to matrices; see also \cite[Chapter~1]{higham2008functions}.
\begin{lemma}[Continuity of the spectral functional calculus]\label{lemm:continuous_spectral_calculus}
Let $f:[0,\infty)\to\mathbb R$ be continuous. For $X\succeq0$, define $f(X)$ by applying $f$ to the eigenvalues of $X$. Then the map
\[
\mathrm{cl}(S_n^+)\longrightarrow S_n,
\qquad X\longmapsto f(X),
\]
is continuous. In particular, $X\mapsto\chi(X)$ and $X\mapsto\kappa(X)$ are continuous on $\mathrm{cl}(S_n^+)$.
\end{lemma}
\begin{proof}
Suppose $X_j\to X$ in operator norm. The spectra of $X$ and all sufficiently large $X_j$ lie in one compact interval $[0,M]$. Given $\varepsilon>0$, the Weierstrass approximation theorem supplies a polynomial $p$ satisfying
\[
\sup_{0\leq r\leq M}|f(r)-p(r)|<\varepsilon.
\]
The spectral theorem then gives $\lVert f(Y)-p(Y)\rVert_{\rm op}<\varepsilon$ for $Y=X$ and for all sufficiently large $Y=X_j$. Since polynomial evaluation $Y\mapsto p(Y)$ is continuous,
\[
\lVert f(X_j)-f(X)\rVert_{\rm op}
\leq 2\varepsilon+\lVert p(X_j)-p(X)\rVert_{\rm op},
\]
and the last term tends to zero. As $\varepsilon$ is arbitrary, $f(X_j)\to f(X)$. In finite dimension, this is equivalent to continuity for the Frobenius or entrywise Euclidean norm.
\end{proof}

\begin{theorem}\label{theom:WKL_continuity}
Equip $\mathbb R^n$ and $S_n$ with their Euclidean topologies, all product spaces below with the product topology, and $[0,\infty]$ with its usual extended-real topology. Define
\[
\mathcal{L}(m_0,\Sigma_0,m_1,\Sigma_1)
=D_{\rm WKL}\!\left(\mathcal{N}(m_0,\Sigma_0)\,\middle\|\,\mathcal{N}(m_1,\Sigma_1)\right)
\]
on $\mathbb{R}^n\times S_n^+\times\mathbb{R}^n\times S_n^+$. Then the following statements hold.

\begin{enumerate}
\item $\mathcal{L}$ is continuous on its domain.

\item $\mathcal{L}$ has a unique continuous extended-real-valued extension $\overline{\mathcal{L}}$ to
\begin{align*}
\mathcal{U}
={}&\left(\mathbb{R}^n\times S_n^+\times\mathbb{R}^n\times\mathrm{cl}(S_n^+)\right)\\
&\quad\cup\left(\mathbb{R}^n\times\mathrm{cl}(S_n^+)\times\mathbb{R}^n\times S_n^+\right).
\end{align*}
If $\Sigma_0\succ0$ and $\Sigma_1\succeq0$, set
\[
R=\Sigma_0^{-1/2}\left(\Sigma_0^{1/2}\Sigma_1\Sigma_0^{1/2}\right)^{1/2}\Sigma_0^{-1/2}.
\]
Then
\begin{align}\label{eq:continuous_boundary_extension}
\overline{\mathcal{L}}(m_0,\Sigma_0,m_1,\Sigma_1)
={}&\frac14\trace{\Sigma_0-\Sigma_1+\Sigma_0\chi(R)}\nonumber\\
&+\frac14(m_1-m_0)^\top\kappa(R)(m_1-m_0).
\end{align}
If $\Sigma_0$ is singular and $\Sigma_1\succ0$, then $\overline{\mathcal{L}}(m_0,\Sigma_0,m_1,\Sigma_1)=+\infty$.

\item There is no continuous extension of $\overline{\mathcal{L}}$ to
\[
\mathcal{X}=\mathbb{R}^n\times\mathrm{cl}(S_n^+)\times\mathbb{R}^n\times\mathrm{cl}(S_n^+).
\]
One lower-semicontinuous extension to $\mathcal{X}$ is
\begin{align}\label{eq:lsc_extension}
\overline{\mathcal{L}}_{\rm ext}(m_0,\Sigma_0,m_1,\Sigma_1)
=\begin{cases}
\overline{\mathcal{L}}(m_0,\Sigma_0,m_1,\Sigma_1),
& \Sigma_0\succ0\ \textnormal{or}\ \Sigma_1\succ0,\\[5pt]
\dfrac14\lVert m_1-m_0\rVert^2,
& \Sigma_0,\Sigma_1\ \textnormal{both singular}.
\end{cases}
\end{align}
At a pair of Dirac measures, this chosen value equals the interior lower limit:
\begin{align}\label{eq:dirac_liminf}
\liminf_{\substack{\Sigma_0,\Sigma_1\succ0\\(\Sigma_0,\Sigma_1)\to(0,0)}}
\mathcal{L}(m_0,\Sigma_0,m_1,\Sigma_1)
=\frac14\lVert m_1-m_0\rVert^2.
\end{align}
\end{enumerate}
\end{theorem}
\begin{proof}
\noindent\textit{Step 1: scalar reduction and interior continuity.}
For $R\succ0$, diagonalizing $R$ shows that
\begin{align}\label{eq:kappa_identity}
\Adag{(R-I)}\left(\log(R^2)R^2-R^2+I\right)\Adag{(R-I)}
+2\Perp{\log R}
=\kappa(R).
\end{align}
Indeed, the scalar multiplier is $\kappa(r)$ on every eigenspace of $R$. The relevant scalar limits are
\[
\lim_{r\downarrow0}r^2\log(r^2)=0,
\qquad
\lim_{r\downarrow0}\kappa(r)=1,
\qquad
\lim_{r\to1}\kappa(r)=2.
\]
Thus $\chi$ and $\kappa$ are continuous on $[0,\infty)$. Lemma~\ref{lemm:continuous_spectral_calculus} implies that $X\mapsto\chi(X)$ and $X\mapsto\kappa(X)$ are continuous on the cone of positive-semidefinite matrices.

\medskip
\noindent\textit{Step 2: target-covariance boundary.}
When $\Sigma_0\succ0$, the map
\[
(\Sigma_0,\Sigma_1)\longmapsto
\Sigma_0^{-1/2}\left(\Sigma_0^{1/2}\Sigma_1\Sigma_0^{1/2}\right)^{1/2}\Sigma_0^{-1/2}
\]
is continuous for $\Sigma_1\succeq0$. Hence \eqref{eq:continuous_boundary_extension} is finite and continuous when the first covariance is positive definite, including when $\Sigma_1$ is singular. On the interior, \eqref{eq:kappa_identity} shows that this formula agrees with \eqref{eq:D_main_formula}. This proves the first assertion and the target-covariance part of the second assertion.

\medskip
\noindent\textit{Step 3: source-covariance boundary.}
Consider a sequence
\[
\Sigma_{0,k}\to\overline{\Sigma}_0,
\qquad
\Sigma_{1,k}\to\overline{\Sigma}_1,
\]
where $\overline{\Sigma}_0$ is singular and $\overline{\Sigma}_1\succ0$. Let $R_k$ be the corresponding matrix. For all sufficiently large $k$, there is a constant $c>0$ such that $\Sigma_{1,k}\succeq cI$. Let $v_k$ be a unit eigenvector of $\Sigma_{0,k}$ associated with $\lambda_{\min}(\Sigma_{0,k})$. By operator monotonicity of the positive-semidefinite square root,
\[
\left(\Sigma_{0,k}^{1/2}\Sigma_{1,k}\Sigma_{0,k}^{1/2}\right)^{1/2}
\succeq\sqrt{c}\,\Sigma_{0,k}^{1/2}.
\]
It follows that
\[
v_k^\top R_kv_k
\geq\sqrt{\frac{c}{\lambda_{\min}(\Sigma_{0,k})}}
\longrightarrow+\infty.
\]
Consequently, if $r_k=\lambda_{\max}(R_k)$, then $r_k\to+\infty$.

Set $M_k=\log(R_k^2)R_k^2-R_k^2+I$. Since $R_k\Sigma_{0,k}R_k=\Sigma_{1,k}$, an orthonormal eigenbasis $\{u_{j,k}\}$ of $R_k$, with eigenvalues $r_{j,k}$, gives
\begin{align*}
\trace{M_k\Sigma_{0,k}}
&=\sum_j\left(2\log r_{j,k}-1+r_{j,k}^{-2}\right)u_{j,k}^\top\Sigma_{1,k}u_{j,k}.
\end{align*}
Every summand is nonnegative because $z-1-\log z\geq0$ for $z>0$, applied with $z=r_{j,k}^{-2}$. Keeping only an eigenvector associated with $r_k$ therefore yields
\[
\trace{M_k\Sigma_{0,k}}
\geq c\left(2\log r_k-1+r_k^{-2}\right)
\longrightarrow+\infty.
\]
The covariance contribution alone therefore diverges, proving continuity with value $+\infty$ on the source-singular part of $\mathcal{U}$. Uniqueness follows because the original domain is dense in $\mathcal{U}$ and $[0,\infty]$ is Hausdorff.

\medskip
\noindent\textit{Step 4: obstruction to a full continuous extension.}
Write $\Delta m=m_1-m_0$. Along $\Sigma_0=\Sigma_1=\varepsilon I$,
\[
\mathcal{L}=\frac12\lVert\Delta m\rVert^2.
\]
Along $\Sigma_0=\varepsilon I$ and $\Sigma_1=\varepsilon^2I$,
\[
\mathcal{L}
=\frac{n}{4}\left(\varepsilon-\varepsilon^2+\varepsilon^2\log\varepsilon\right)
+\frac14\kappa(\sqrt{\varepsilon})\lVert\Delta m\rVert^2
\longrightarrow\frac14\lVert\Delta m\rVert^2.
\]
These paths have different limits when $\Delta m\neq0$. If $\Delta m=0$, compare the diagonal path with
\[
\Sigma_0=\varepsilon e^{-1/\varepsilon}I,
\qquad
\Sigma_1=\varepsilon I.
\]
Along the latter path,
\[
\mathcal{L}
=\frac{n}{4}\left(\varepsilon e^{-1/\varepsilon}-\varepsilon+1\right)
\longrightarrow\frac n4,
\]
whereas the diagonal limit is zero. Thus continuity fails at every point $(m_0,0,m_1,0)$.

\medskip
\noindent\textit{Step 5: lower-semicontinuous extension and the Dirac lower limit.}
Since $\mathcal{U}$ is open in $\mathcal{X}$ and $\overline{\mathcal{L}}$ is continuous on $\mathcal{U}$, it remains only to check lower semicontinuity at points where both covariance matrices are singular. The uniform bound \eqref{eq:uniform_lower_bound} passes to the one-sided boundary extension:
\[
\overline{\mathcal{L}}(m_0,\Sigma_0,m_1,\Sigma_1)
\geq\frac14\lVert m_1-m_0\rVert^2
\]
whenever $(m_0,\Sigma_0,m_1,\Sigma_1)\in\mathcal{U}$; on $\mathcal X\setminus\mathcal U$, the same bound holds with equality by \eqref{eq:lsc_extension}. If a sequence converges to a point at which both covariance matrices are singular, this bound and the continuity of the means imply
\[
\liminf_k\overline{\mathcal{L}}_{\rm ext}(m_{0,k},\Sigma_{0,k},m_{1,k},\Sigma_{1,k})
\geq\frac14\lVert m_1-m_0\rVert^2.
\]
This is precisely the value assigned in \eqref{eq:lsc_extension}; hence the extension is lower semicontinuous. Equation \eqref{eq:dirac_liminf} follows from the uniform lower bound and the path $\Sigma_0=\varepsilon I$, $\Sigma_1=\varepsilon^2I$.
\end{proof}

\begin{remark}
The extension \eqref{eq:lsc_extension} is not asserted to be unique. Moreover, on the full singular--singular stratum it need not separate distinct degenerate Gaussian measures: if $m_0=m_1$, it assigns zero even when $\Sigma_0\neq\Sigma_1$. It is therefore used only as an extended lower-semicontinuous function. At pairs of Dirac measures, its value is distinguished by the exact lower-limit identity \eqref{eq:dirac_liminf}.
\end{remark}

We now specialize Theorem~\ref{theom:WKL_continuity} to the univariate case, using standard deviations as coordinates.
\begin{samepage}
\begin{corollary}\label{corr:univariate_boundary_extension}
Define
\begin{align*}
\mathcal{L}_1:\mathbb{R}\times(0,\infty)\times\mathbb{R}\times(0,\infty)&\longrightarrow[0,\infty),\\
\mathcal{L}_1(m_0,\sigma_0,m_1,\sigma_1)
&=D_{\rm WKL}\!\left(\mathcal{N}(m_0,\sigma_0^2)\,\middle\|\,\mathcal{N}(m_1,\sigma_1^2)\right).
\end{align*}
Set
\[
\mathcal{U}_1
:=\bigl(\mathbb{R}\times(0,\infty)\times\mathbb{R}\times[0,\infty)\bigr)
\cup\bigl(\mathbb{R}\times[0,\infty)\times\mathbb{R}\times(0,\infty)\bigr).
\]
Its unique continuous extended-real-valued extension to $\mathcal{U}_1$ is
\begin{align}\label{eq:continuous_boundary_extension_univariate}
\overline{\mathcal{L}}_1(m_0,\sigma_0,m_1,\sigma_1)
=\begin{cases}
\mathcal{L}_1(m_0,\sigma_0,m_1,\sigma_1),
& \sigma_0>0,\ \sigma_1>0,\\[5pt]
\dfrac14\left(\sigma_0^2+(m_1-m_0)^2\right),
& \sigma_0>0,\ \sigma_1=0,\\[7pt]
+\infty,
& \sigma_0=0,\ \sigma_1>0.
\end{cases}
\end{align}
One lower-semicontinuous extension to $\mathbb{R}\times[0,\infty)\times\mathbb{R}\times[0,\infty)$ is
\begin{align}\label{eq:lsc_extension_univariate}
\overline{\mathcal{L}}_{1,\rm ext}(m_0,\sigma_0,m_1,\sigma_1)
=\begin{cases}
\overline{\mathcal{L}}_1(m_0,\sigma_0,m_1,\sigma_1),
& \sigma_0>0\ \textnormal{or}\ \sigma_1>0,\\[5pt]
\dfrac14(m_1-m_0)^2,
& \sigma_0=\sigma_1=0.
\end{cases}
\end{align}
The chosen value at the origin equals the interior lower limit, but no continuous extension exists there.
\end{corollary}
\end{samepage}

\begin{figure}[p]
	\centering
	\begin{minipage}{0.47\textwidth}
		\centering
     \includegraphics[width=\linewidth]{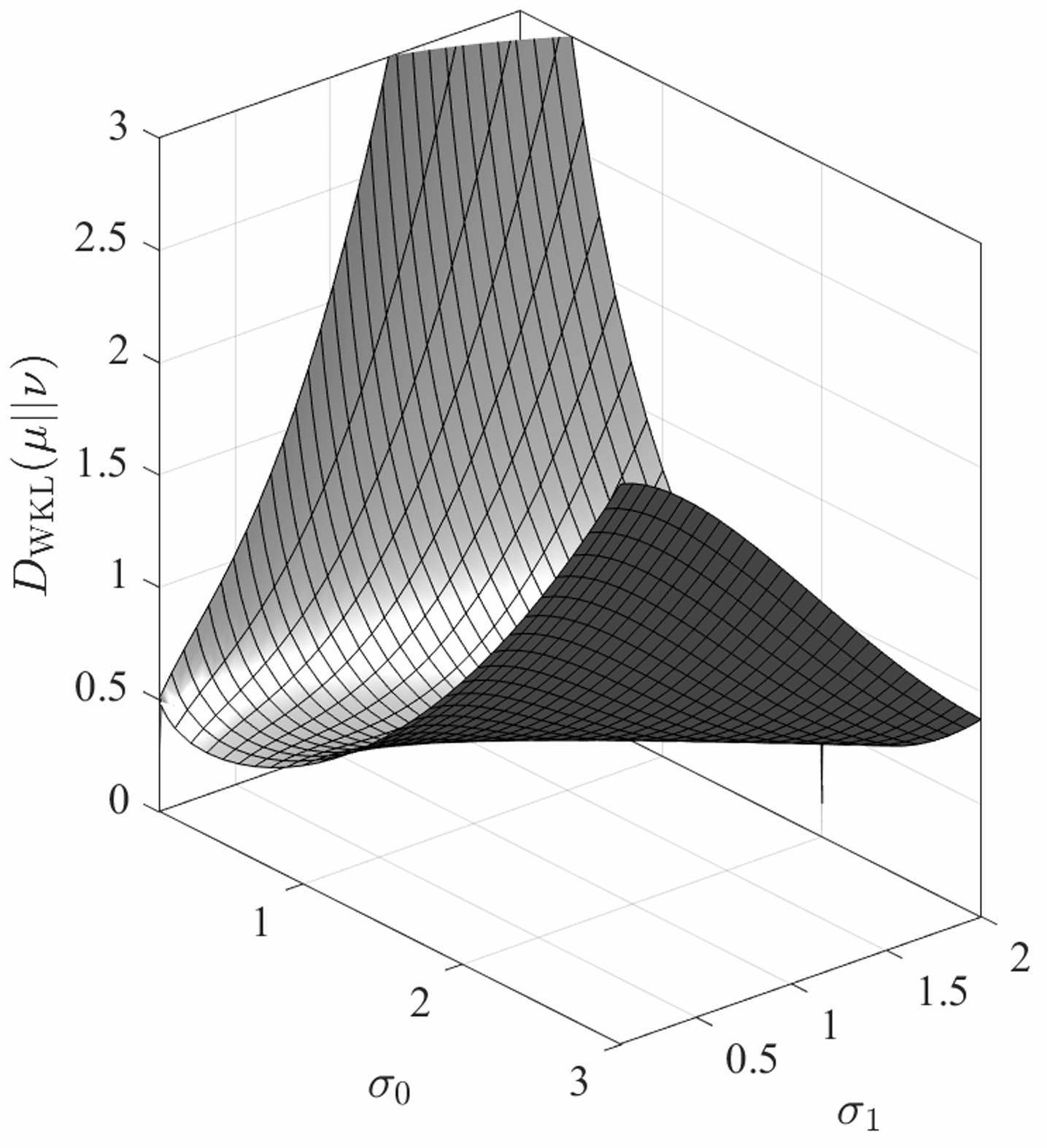}
	\end{minipage}
	\hfill
	\begin{minipage}{0.47\textwidth}
		\centering
     \includegraphics[width=\linewidth]{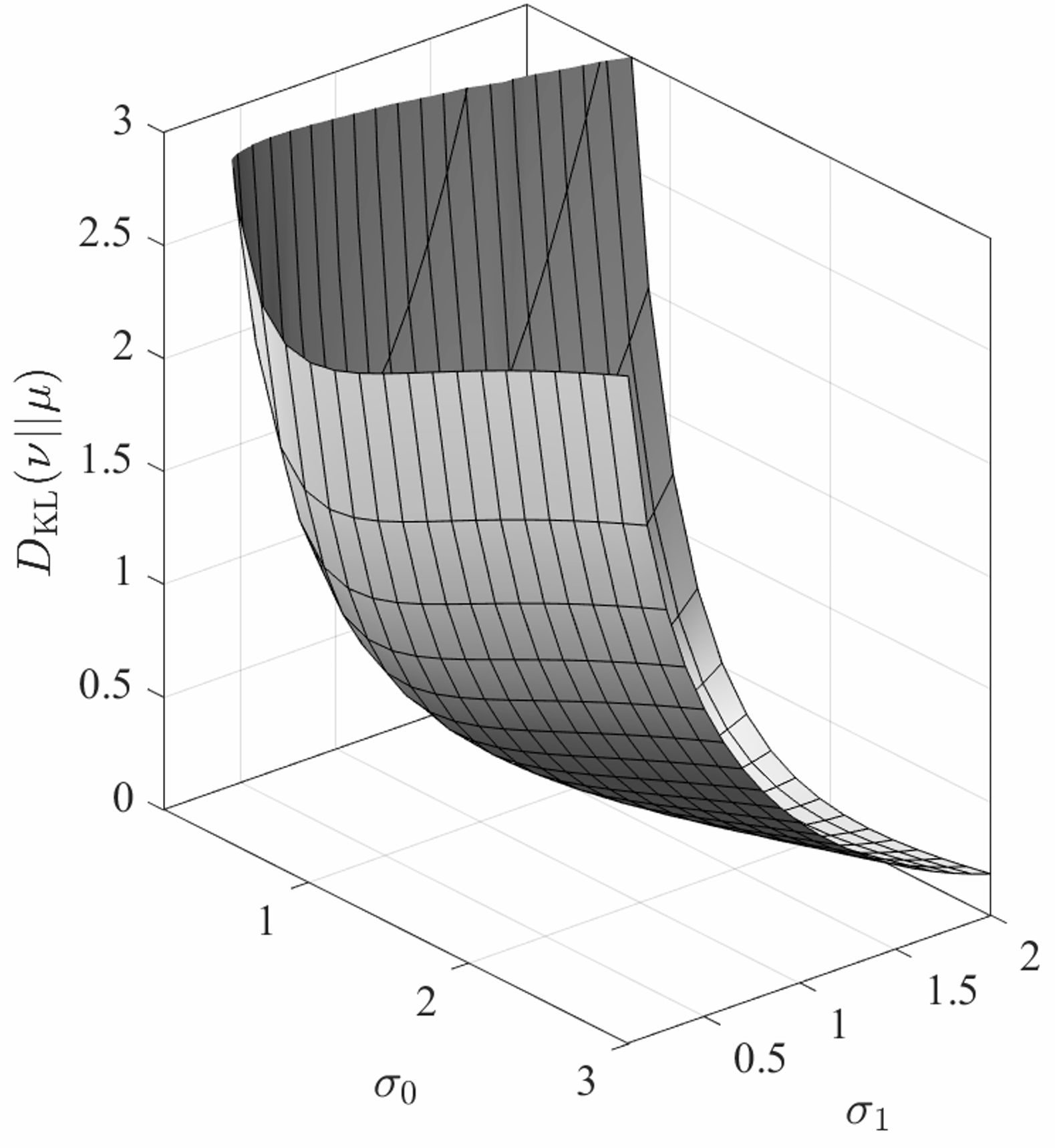}
	\end{minipage}
	\caption{Surface plots for $D_{\rm WKL}(\mu\lVert\nu)$ (left) and $D_{\rm KL}(\nu\lVert\mu)$ (right), where $\mu=\mathcal{N}(0,\sigma_0^2)$ and $\nu=\mathcal{N}(1,\sigma_1^2)$, over $0<\sigma_0\leq3$ and $0<\sigma_1\leq2$. Values above $3$ are clipped only for visualization.}
	\label{fig:surfaceplots}
\end{figure}

\begin{proof}
The finite target-boundary value follows from $\kappa(0)=1$ and $\chi(0)=0$. For fixed $\sigma_1>0$, the covariance contribution diverges as $\sigma_0\downarrow0$ because
\[
\sigma_0^2-\sigma_1^2+\sigma_1^2\log\!\left(\frac{\sigma_1^2}{\sigma_0^2}\right)\longrightarrow+\infty.
\]
Lower semicontinuity follows from Theorem~\ref{theom:WKL_continuity}.

For $m_0\neq m_1$, the paths $\sigma_0=\sigma_1=\varepsilon$ and $\sigma_0=\varepsilon$, $\sigma_1=\varepsilon^2$ have respective limits
\[
\frac12(m_1-m_0)^2
\qquad\textnormal{and}\qquad
\frac14(m_1-m_0)^2.
\]
If $m_0=m_1$, the diagonal limit is zero, whereas along
\[
\sigma_1=\varepsilon,
\qquad
\sigma_0=\varepsilon\exp\!\left(-\frac{1}{2\varepsilon^2}\right)
\]
the limit is $1/4$. Hence continuity fails at every point with $\sigma_0=\sigma_1=0$.
\end{proof}
Figure~\ref{fig:WKL_behavior} summarizes the resulting univariate boundary regimes.
\begin{figure}
    \centering
    \begin{tikzpicture}[thick, >=stealth]

    \draw[->] (-0.5,0) -- (5.5,0) node[right] {$\sigma_0$};
    \draw[->] (0,-0.5) -- (0,5.5) node[above] {$\sigma_1$};
    
    \draw[red!65!black, ultra thick, dashed] (0,0) -- (0,5);
    \node[red!65!black] at (-0.5,2.5) {$+\infty$};
    
    \draw[blue, ultra thick, dash dot] (0,0) -- (4,4);
    \node[blue] at (3.2,1.8) {$\frac{1}{2}(m_1 - m_0)^2$};
    
    \draw[green!50!black, ultra thick, dotted] (0,0) -- (5,0);
    \node[green!50!black] at (3,-0.4) {$\frac{1}{4}\sigma_0^2 + \frac{1}{4}(m_1 - m_0)^2$};
    
    \filldraw (0,0) circle (2pt);
    \draw[->, thick] (-0.5,-0.3) to[bend right=-10] (-0.1,-0.1)
    node[midway, above,xshift=-1.5cm, yshift=-0.8cm] {discontinuity};
    \end{tikzpicture}
\caption{Boundary values and path limits of $D_{\rm WKL}(\mu\lVert\nu)$ for $\mu=\mathcal{N}(m_0,\sigma_0^2)$ and $\nu=\mathcal{N}(m_1,\sigma_1^2)$. The target-boundary value at $\sigma_1=0$ is $\frac14\sigma_0^2+\frac14(m_1-m_0)^2$ (horizontal dotted line), whereas the source-boundary value at $\sigma_0=0$ is $+\infty$ (vertical dashed line). Along $\sigma_0=\sigma_1\downarrow0$, the limit is $\frac12(m_1-m_0)^2$ (dash-dotted diagonal). The black point marks the chosen lower-semicontinuous extension value $\frac14(m_1-m_0)^2$ at the origin, where continuity fails.}
    \label{fig:WKL_behavior}
\end{figure}
Figure~\ref{fig:surfaceplots} provides a complementary view of the asymmetric one-sided behavior in the interior.
\FloatBarrier
Finally, when $\Sigma_0=\Sigma_1=\Sigma\succ0$, the WKL-divergence equals $\frac12\lVert m_1-m_0\rVert^2$ and is independent of $\Sigma$. In contrast,
\begin{align*}
D_{\rm KL}\!\left(\mathcal{N}(m_1,\Sigma)\,\middle\|\,\mathcal{N}(m_0,\Sigma)\right)
=\frac12\left\lVert\Sigma^{-1/2}(m_1-m_0)\right\rVert^2.
\end{align*}
In the equal-covariance transition model of \cite{stein2026wklcontrol}, this covariance independence makes the WKL divergence regularizer remain finite as the process-noise covariance vanishes.
Within this equal-covariance regime, the reverse-KL diverges along paths for which $\lVert\Sigma^{-1/2}(m_1-m_0)\rVert\to\infty$; in particular, if $m_0\neq m_1$ and $\Sigma=\varepsilon^2I$ with $\varepsilon\downarrow0$, then it diverges to $\infty$.

\section{Conclusion}\label{sec:conclusion}
We have given a direct formulation of the Wasserstein Kullback--Leibler divergence on the Gaussian manifold over the noncompact sample space $\mathbb R^n$. Within the canonical-divergence construction, WKL divergence follows the same conjugate-connection route as reverse KL: the Fisher--Rao/$e_0$ pair is replaced by the Otto/$e_1$ pair, thereby incorporating Euclidean ground geometry. The $e_1$-connection is geodesically complete. Every ordered pair of nondegenerate Gaussian measures has a unique quadratic $e_1$ connector, and we classified the possible forward limits of these geodesics.

We give an explicit formula for the WKL divergence with separate covariance and mean contributions. WKL divergence is nonnegative and separates nondegenerate Gaussian measures. For equal covariances it reduces to $\frac12\lVert m_1-m_0\rVert^2$, independently of the common covariance.

The boundary analysis makes the asymmetry of WKL divergence precise. It extends finitely and continuously to positive-semidefinite target covariances from a positive-definite source, whereas it diverges when the source covariance becomes singular while the target remains positive definite. At every pair of Dirac measures, different joint covariance-collapse paths yield different limits, ruling out a continuous extension to the full product of positive-semidefinite covariance cones. The lower-semicontinuous extended-real extension constructed here attains the interior lower limit at Dirac pairs, but it is not separating on the full singular--singular stratum. Thus WKL divergence has a natural one-sided boundary extension, but admits no continuous extension to the full product of positive-semidefinite covariance cones.

\section*{Statements and Declarations}
\noindent\textbf{Funding.}
The authors did not receive support from any organization for the submitted work.

\medskip
\noindent\textbf{Competing interests.}
Nihat Ay serves as an Editor-in-Chief of \emph{Information Geometry}. Apart from this editorial role, the authors declare no relevant financial or non-financial competing interests. Nihat Ay will have no knowledge of the review process and will take no part in reviewer selection, editorial discussion, or the decision, and the authors request fully independent editorial handling.

\medskip
\noindent\textbf{Use of generative artificial intelligence.}
OpenAI's Codex was used during manuscript revision to assist with literature discovery, mathematical cross-checking, drafting, and editing. The authors independently reviewed and verified all AI-assisted content, checked the cited sources and derivations, and take full responsibility for the manuscript.

\medskip
\noindent\textbf{Data availability.}
No datasets were generated or analyzed during this study; consequently, data sharing is not applicable.

\appendix

\section{The Intrinsic Gaussian Mixture Connection and Its Canonical Divergences}\label{sec:app_mixture_duality}
On an ambient space of densities, the mixture connection is induced by the affine structure of signed measures, and its geodesics are ordinary density mixtures.  This description cannot be restricted literally to $\mathcal G$, because $(1-t)\mu_0+t\mu_1$ is generally not Gaussian.  The intrinsic mixture connection of the Gaussian exponential family is instead affine in expectation coordinates \cite{amari_2016}.  We now make that distinction explicit and verify directly, in the present notation, the corresponding conjugacy on $\mathcal G$ from \cite{ay2024information,ay2025correction}.  Strictly speaking, metric conjugacy is a relation between connections (equivalently, their parallel transports), not a pointwise duality between individual geodesic curves.

\subsection{Connection, Geodesics, and Parallel Transport}
For $\mu=\mathcal N(m,\Sigma)$, define the moment chart
\begin{equation}\label{eq:gaussian_moment_chart}
 \xi(\mu)=(m,M),\qquad M:=\mathbb E_\mu[XX^\top]=\Sigma+mm^\top .
\end{equation}
Its image is
\[
 \mathcal D=\{(m,M)\in\mathbb R^n\times S_n:M-mm^\top\succ0\}.
\]
This set is open and convex.  Indeed, for two points of $\mathcal D$, the covariance along their affine moment segment is
\begin{equation}\label{eq:intrinsic_mixture_geodesic}
 \begin{split}
 m_t&=(1-t)m_0+tm_1,\\
 \Sigma_t&=(1-t)\Sigma_0+t\Sigma_1
       +t(1-t)(m_1-m_0)(m_1-m_0)^\top,
 \end{split}
 \qquad 0\leq t\leq1,
\end{equation}
which is positive definite.  We define $\nabla^{(m)}$ to be the flat, torsion-free connection for which $\xi$ is affine,i.e., \eqref{eq:intrinsic_mixture_geodesic} is its unique geodesic segment between the endpoints.

There is a precise relation to the ambient density mixture
\(
 \bar\mu_t=(1-t)\mu_0+t\mu_1
\).
By linearity of integration,
\[
 \int x\,d\bar\mu_t=(1-t)m_0+tm_1,
 \qquad
 \int xx^\top\,d\bar\mu_t=(1-t)M_0+tM_1,
\]
where $M_j=\Sigma_j+m_jm_j^\top$.  These are the expectation coordinates of
the multivariate Gaussian family, and the inverse relation
$\Sigma=M-mm^\top$ shows that they determine a unique
Gaussian \cite{malago2015information}.
Thus, although $\bar\mu_t$ generally lies outside $\mathcal G$, the Gaussian
in \eqref{eq:intrinsic_mixture_geodesic} is the unique Gaussian having the
same first two moments as $\bar\mu_t$; it is not the density mixture itself.

If a tangent vector $a$ at $\mu=\mathcal N(m,\Sigma)$ has coordinate velocity $(u,U)\in\mathbb R^n\times S_n$, then
\[
 D\xi_\mu(u,U)=(u,U+um^\top+mu^\top).
\]
The parallel transport of $\nabla^{(m)}$ is therefore
\begin{equation}\label{eq:gaussian_mixture_transport}
 \Pi^{(m)}_{\mu,\nu}:=D\xi_\nu^{-1}D\xi_\mu,
 \qquad
 \Pi^{(m)}_{\mu,\nu}(u,U)
 =\bigl(u,U+u(m-m_\nu)^\top+(m-m_\nu)u^\top\bigr),
\end{equation}
where $\nu=\mathcal N(m_\nu,\Sigma_\nu)$.  Equivalently, this is the unique tangent transport preserving the action on all quadratic statistics.  Namely, for
\(
 q(x)=\frac12x^\top Qx+r^\top x+c,
 \qquad Q\in S_n,\quad r\in\mathbb R^n,\quad c\in\mathbb R,
\),
\begin{equation}\label{eq:quadratic_tangent_pairing}
 \int_{\mathbb R^n}q\,da
 =r^\top u+\frac12\operatorname{Tr}\!\left(
 Q(U+um^\top+mu^\top)\right),
\end{equation}
Writing $W:=U+um^\top+mu^\top$, the right-hand side depends only on the
affine tangent components $(u,W)=D\xi_\mu(u,U)$.  Since the parallel
transport in \eqref{eq:gaussian_mixture_transport} keeps $(u,W)$ unchanged,
the value of $\int q\,da$ remains unchanged under this transport; hence
\begin{equation}\label{eq:mixture_pairing_invariance}
 \int q\,d\bigl(\Pi^{(m)}_{\mu,\nu}a\bigr)=\int q\,da
 \qquad(q\in\mathcal Q).
\end{equation}

\subsection{Metric Conjugacy}
\begin{samepage}
\begin{proposition}[Metric conjugacy on the Gaussian manifold]\label{prop:gaussian_mixture_duality}
For every $\mu,\nu\in\mathcal G$, the intrinsic mixture connection $\nabla^{(m)}$ is conjugate to $\nabla^{(e_0)}$ with respect to the Fisher--Rao metric and to $\nabla^{(e_1)}$ with respect to the Otto metric.  Equivalently, for $i\in \{0,1\}$, $a,b\in T_\mu\mathcal G$, and the corresponding metric $g^{(i)}$,
\begin{equation}\label{eq:dual_parallel_transports}
 g^{(i)}_\nu\!\left(\Pi^{(m)}_{\mu,\nu}a,
              \Pi^{(e_i)}_{\mu,\nu}b\right)
 =g^{(i)}_\mu(a,b),
 \qquad
 g^{(0)}=g^{FR},\quad g^{(1)}=g^O.
\end{equation}
\end{proposition}
\end{samepage}

\begin{proof}
For $q+\mathbb R\in\mathcal Q/\mathbb R$, write the two metric identifications as
\[
 \phi^{(0)}_\mu(q+\mathbb R)
   =(q-\mathbb E_\mu q)\mu,
 \qquad
 \phi^{(1)}_\mu(q+\mathbb R)
   =-\Delta_\mu(q+\mathbb R)\mu.
\]
The definitions of the Fisher--Rao and Otto metrics, together with the calculation described in Section~\ref{sec:Otto_metric} in the second case, give the common pairing identity
\begin{equation}\label{eq:metric_cotangent_pairing}
 g^{(i)}_\mu\bigl(a,\phi^{(i)}_\mu(q+\mathbb R)\bigr)
 =\int q\,da,
 \qquad i=0,1.
\end{equation}
The $e_i$ transport keeps the quadratic cotangent representative fixed:
\(
 \Pi^{(e_i)}_{\mu,\nu}=\phi^{(i)}_\nu(\phi^{(i)}_\mu)^{-1}
\).
Take $q+\mathbb R=(\phi^{(i)}_\mu)^{-1}b$.  Applying \eqref{eq:metric_cotangent_pairing} at $\nu$, then \eqref{eq:mixture_pairing_invariance}, and finally \eqref{eq:metric_cotangent_pairing} at $\mu$ yields
\[
 g^{(i)}_\nu\!\left(\Pi^{(m)}_{\mu,\nu}a,
                    \Pi^{(e_i)}_{\mu,\nu}b\right)
 =\int q\,d(\Pi^{(m)}_{\mu,\nu}a)
 =\int q\,da
 =g^{(i)}_\mu(a,b).
\]
Constancy of this pairing along every curve is precisely conjugacy of the two connections.
\end{proof}

For the Fisher--Rao metric, this recovers the classical duality of Gaussian natural coordinates $(\Sigma^{-1}m,-\frac12\Sigma^{-1})$ and expectation coordinates $(m,M)$ \cite{amari_2016}.  For the Otto metric, the same intrinsic mixture connection is conjugate to $\nabla^{(e_1)}$, but the conclusion is not dual flatness because $\nabla^{(e_1)}$ has nonzero torsion \cite{ay2026torsion}.

\subsection{Canonical Divergences along Mixture Geodesics}\label{sec:app_mixture_canonical}
We now evaluate the canonical-divergence functional \eqref{eq:KL_energy_formula}
along the intrinsic Gaussian mixture geodesic \eqref{eq:intrinsic_mixture_geodesic}.
Because the connection is fixed while the metric is varied, we make both choices
explicit and write
\begin{equation}\label{eq:mixture_canonical_definition}
 D_g^{(m)}(\mu\Vert\nu)
 :=\int_0^1 t\,
 g_{\gamma_m(t)}\!\left(\dot\gamma_m(t),\dot\gamma_m(t)\right)dt,
 \qquad g\in\{g^{FR},g^O\}.
\end{equation}
Here $\gamma_m$ is oriented from
$\mu=\mathcal N(m_0,\Sigma_0)$ to
$\nu=\mathcal N(m_1,\Sigma_1)$, and both covariance matrices are positive
definite. Throughout this subsection, ``mixture geodesic'' means the intrinsic
Gaussian geodesic, which is affine in the moment coordinates $(m,M)$. Set
\begin{equation}\label{eq:mixture_path_derivatives}
 \begin{aligned}
 d&:=m_1-m_0,\\
 \Sigma_t&=(1-t)\Sigma_0+t\Sigma_1+t(1-t)dd^\top,\\
 V_t&:=\dot\Sigma_t=\Sigma_1-\Sigma_0+(1-2t)dd^\top.
 \end{aligned}
\end{equation}

\subsubsection{Fisher--Rao Metric: Forward KL}
The following identity is the classical exponential--mixture duality on the
Gaussian exponential family \cite{amari_2016}; we include the calculation to fix
the endpoint orientation under the convention in \eqref{eq:KL_energy_formula}.
\begin{proposition}[Mixture canonical divergence for the Fisher--Rao metric]
\label{prop:mixture_fr_forward_kl}
Along the intrinsic Gaussian mixture geodesic,
\begin{equation}\label{eq:mixture_fr_forward_kl}
 D_{\mathrm{FR}}^{(m)}(\mu\Vert\nu)
 =D_{\rm KL}(\mu\Vert\nu).
\end{equation}
Thus the orientation is the forward order, in contrast to the reverse order
of the Fisher--Rao/$e_0$ construction in
\eqref{eq:KL_energy_formula_specialize_to_KL}.
\end{proposition}

\begin{proof}
In the moment chart $\xi=(m,M)$ of \eqref{eq:gaussian_moment_chart}, consider
the negative-entropy potential, up to an irrelevant additive constant,
\[
 \Phi(m,M):=-\frac12\log\det(M-mm^\top)=-\frac12\log\det\Sigma.
\]
The curve $\xi_t=(1-t)\xi_0+t\xi_1$ is affine, while
$\ddot\Sigma_t=-2dd^\top$.  The Fisher--Rao norm of a Gaussian coordinate
velocity $(\dot m,\dot\Sigma)$ is
\[
 \lVert(\dot m,\dot\Sigma)\rVert_{\mathrm{FR}}^2
 =\dot m^\top\Sigma^{-1}\dot m
  +\frac12\operatorname{Tr}\!\left(
      \Sigma^{-1}\dot\Sigma\Sigma^{-1}\dot\Sigma\right).
\]
Consequently, along $\gamma_m$,
\begin{align*}
 \frac{d^2}{dt^2}\Phi(m_t,M_t)
 &=\frac12\operatorname{Tr}\!\left(
      \Sigma_t^{-1}V_t\Sigma_t^{-1}V_t\right)
   -\frac12\operatorname{Tr}(\Sigma_t^{-1}\ddot\Sigma_t)\\
 &=\frac12\operatorname{Tr}\!\left(
      \Sigma_t^{-1}V_t\Sigma_t^{-1}V_t\right)
   +d^\top\Sigma_t^{-1}d
 =\lVert\dot\gamma_m(t)\rVert_{\gamma_m(t),\mathrm{FR}}^2.
\end{align*}
Integration by parts therefore gives
\begin{align*}
 D_{\mathrm{FR}}^{(m)}(\mu\Vert\nu)
 &=\int_0^1t\,\Phi''(t)\,dt
 =\Phi'(1)-\Phi(1)+\Phi(0).
\end{align*}
Since $V_1=\Sigma_1-\Sigma_0-dd^\top$,
\[
 \Phi'(1)
 =-\frac12\operatorname{Tr}(\Sigma_1^{-1}V_1)
 =\frac12\left[
   -n+\operatorname{Tr}(\Sigma_1^{-1}\Sigma_0)
   +d^\top\Sigma_1^{-1}d\right].
\]
Combining this with
$\Phi(0)-\Phi(1)=\frac12\log(\det\Sigma_1/\det\Sigma_0)$ yields
\begin{align*}
 D_{\mathrm{FR}}^{(m)}(\mu\Vert\nu)
 =\frac12\left[
  \operatorname{Tr}(\Sigma_1^{-1}\Sigma_0)
  +d^\top\Sigma_1^{-1}d-n
  +\log\frac{\det\Sigma_1}{\det\Sigma_0}\right],
\end{align*}
which is precisely $D_{\rm KL}(\mu\Vert\nu)$.
\end{proof}

\subsubsection{Otto Metric: General Formula}
\begin{samepage}
Recall the Lyapunov isomorphism defined in \eqref{eq:lyapunov_operator},
and define $A_t:=\mathcal L_{\Sigma_t}^{-1}(V_t)$.  The affine gradient
field representing $\dot\gamma_m(t)$ is
\[
 x\longmapsto A_t(x-m_t)+d.
\]
\end{samepage}
The Gaussian Otto metric \eqref{eq:otto_gaussian_metric} therefore gives
\begin{align}
 \lVert\dot\gamma_m(t)\rVert_{\gamma_m(t),O}^2
 &=\lVert d\rVert^2+\operatorname{Tr}(A_t^2\Sigma_t)\notag\\
 &=\lVert d\rVert^2
   +\frac12\operatorname{Tr}\!\left(
      V_t\mathcal L_{\Sigma_t}^{-1}(V_t)\right),
 \label{eq:mixture_otto_speed}
\end{align}
where the second equality follows from
$\operatorname{Tr}(A_tV_t)=2\operatorname{Tr}(A_t^2\Sigma_t)$.
Hence the Otto/mixture canonical divergence has the exact representation
\begin{equation}\label{eq:mixture_otto_general}
 \boxed{
 D_O^{(m)}(\mu\Vert\nu)
 =\frac12\lVert d\rVert^2
  +\frac12\int_0^1t\,
   \operatorname{Tr}\!\left(
    V_t\mathcal L_{\Sigma_t}^{-1}(V_t)\right)dt .}
\end{equation}
Because $\Sigma_t\succ0$ uniformly on $[0,1]$, the integrand is smooth and
finite.  Its trace term is nonnegative and vanishes at a fixed $t$ exactly
when $V_t=0$.  Moreover,
$D_O^{(m)}(\mu\Vert\nu)\geq\frac12\lVert d\rVert^2$.  Equality of the full
divergence to zero therefore forces $d=0$, in which case
$V_t\equiv\Sigma_1-\Sigma_0$; hence it also forces $\Sigma_1=\Sigma_0$.
Thus \eqref{eq:mixture_otto_general} is a separating divergence on the
nondegenerate Gaussian manifold.


Suppose $d=0$ and $\Sigma_0$ and $\Sigma_1$ commute.  If
$\lambda_{0i}$ and $\lambda_{1i}$ are their eigenvalues in a common
orthonormal basis, set
\[
 \delta_i:=\lambda_{1i}-\lambda_{0i},
 \qquad
 \lambda_i(t):=\lambda_{0i}+t\delta_i.
\]
In that basis, $V_t=\operatorname{diag}(\delta_i)$ and
$\Sigma_t=\operatorname{diag}(\lambda_i(t))$, so the Lyapunov equation gives
\[
 \mathcal L_{\Sigma_t}^{-1}(V_t)
 =\operatorname{diag}\!\left(\frac{\delta_i}{2\lambda_i(t)}\right),
 \qquad
 \operatorname{Tr}\!\left(V_t\mathcal L_{\Sigma_t}^{-1}(V_t)\right)
 =\sum_{i=1}^n\frac{\delta_i^2}{2\lambda_i(t)}.
\]
Consequently,
\[
 D_O^{(m)}(\mu\Vert\nu)
 =\frac14\sum_{i=1}^n\int_0^1
   \frac{t\delta_i^2}{\lambda_{0i}+t\delta_i}\,dt,
 \qquad
 \int_0^1\frac{t\delta_i^2}{\lambda_{0i}+t\delta_i}\,dt
 =\delta_i-\lambda_{0i}\log\frac{\lambda_{1i}}{\lambda_{0i}},
\]
with the last identity interpreted by continuity when $\delta_i=0$.
Substitution gives
\begin{equation}\label{eq:mixture_otto_commuting_equal_mean}
 D_O^{(m)}(\mu\Vert\nu)
 =\frac14\sum_{i=1}^n\left[
   \lambda_{1i}-\lambda_{0i}
   -\lambda_{0i}\log\frac{\lambda_{1i}}{\lambda_{0i}}
  \right].
\end{equation}
This happens to equal $D_{\rm WKL}(\nu\Vert\mu)$ in this restricted
equal-mean, commuting case.  The equality does not persist in general.

\subsubsection{Otto Metric: Univariate Formula}
Let $v_i=\sigma_i^2>0$, let $d=m_1-m_0$, and set $c=d^2$.  Along the
intrinsic mixture geodesic,
\[
 v_t=v_0+t(v_1-v_0)+t(1-t)c,
 \qquad
 \dot v_t=v_1-v_0+(1-2t)c.
\]
Since the one-dimensional Otto metric is
$dm^2+\frac1{4v}dv^2$, \eqref{eq:mixture_otto_general} becomes
\begin{equation}\label{eq:mixture_otto_univariate_integral}
 D_O^{(m)}(\mu\Vert\nu)
 =\int_0^1t\left(c+\frac{\dot v_t^2}{4v_t}\right)dt.
\end{equation}
If $c=0$, this reduces to
\begin{equation}\label{eq:mixture_otto_univariate_equal_mean}
 D_O^{(m)}(\mu\Vert\nu)
 =\frac14\left[
   v_1-v_0-v_0\log\frac{v_1}{v_0}\right].
\end{equation}
For $c>0$, define
\[
 a:=v_1-v_0+c,
 \qquad
 h:=\sqrt{a^2+4cv_0}.
\]
Then $v_t=v_0+at-ct^2$ and
$\dot v_t^2=h^2-4cv_t$, so the two terms in
\eqref{eq:mixture_otto_univariate_integral} simplify to
\[
 D_O^{(m)}(\mu\Vert\nu)=\frac{h^2}{4}\int_0^1\frac{t}{v_t}\,dt.
\]
To evaluate the remaining integral, write
\[
 I:=\int_0^1\frac{t}{v_t}\,dt,
 \qquad
 J:=\int_0^1\frac{dt}{v_t}.
\]
Since $\dot v_t=a-2ct$, it follows that
\[
 I=\frac{a}{2c}J-\frac{1}{2c}\int_0^1\frac{\dot v_t}{v_t}\,dt
  =\frac{a}{2c}J-\frac{1}{2c}\log\frac{v_1}{v_0}.
\]
Moreover,
\[
 v_t=\frac{h^2-(2ct-a)^2}{4c},
\]
and the substitution $z=(2ct-a)/h$ to change the variable of integration gives
\[
 J=\frac{2}{h}\left[
   \operatorname{artanh}\!\left(\frac{a}{h}\right)
  +\operatorname{artanh}\!\left(\frac{2c-a}{h}\right)
 \right].
\]
Substitution into $D_O^{(m)}=h^2I/4$ yields the closed form
\begin{align}
 D_O^{(m)}(\mu\Vert\nu)
 &=\frac{ah}{4c}\left[
   \operatorname{artanh}\!\left(\frac{a}{h}\right)
  +\operatorname{artanh}\!\left(\frac{2c-a}{h}\right)
  \right]
  -\frac{h^2}{8c}\log\frac{v_1}{v_0}.
 \label{eq:mixture_otto_univariate_closed}
\end{align}
Both inverse-hyperbolic-tangent arguments lie in $(-1,1)$ because
$h^2-a^2=4cv_0>0$ and $h^2-(2c-a)^2=4cv_1>0$.


\section{Canonical Divergences along the \texorpdfstring{$e_0$- and $e_1$-Geodesics}{e0- and e1-Geodesics}}\label{sec:app_exponential_canonical}
We conclude with the two corresponding reductions along the $e_0$- and $e_1$-geodesics. For the Fisher--Rao metric, the $e_0$-geodesic yields reverse KL; for the Otto metric, the $e_1$-geodesic yields WKL.

\subsection{Fisher--Rao/\texorpdfstring{$e_0$}{e0}: Reverse KL}\label{sec:app_reverse_KL}
We give a direct derivation of the identity in \eqref{eq:KL_energy_formula_specialize_to_KL}, including its argument order. Let $\mu$ and $\nu$ be nondegenerate Gaussian measures and set
\[
 \ell:=\log\frac{d\nu}{d\mu},
 \qquad
 Z(t):=\int_{\Rn}e^{t\ell}\,d\mu,
 \qquad
 \psi(t):=\log Z(t).
\]
The $e_0$-geodesic from $\mu$ to $\nu$ given in \eqref{eq:e0-geodesic} can be written using the above notation as
\begin{equation}\label{eq:app_e0_exponential_arc}
 \gamma_0(t)=e^{t\ell-\psi(t)}\mu,
 \qquad 0\leq t\leq1.
\end{equation}
If $p_0$ and $p_1$ are the Lebesgue densities of $\mu=\mathcal N(m_0,\Sigma_0)$ and $\nu=\mathcal N(m_1,\Sigma_1)$, respectively, then $p_0^{1-t}p_1^t$ is an unnormalized Gaussian kernel with precision
\[
 P_t=(1-t)\Sigma_0^{-1}+t\Sigma_1^{-1}\succ0
 \qquad (0\leq t\leq1).
\]
The positive-definite cone is open and the segment is compact, so $P_t$ remains positive definite on a neighborhood of $[0,1]$. Hence $Z$ and $\psi$ are smooth there, which justifies the differentiations below.

Differentiating under the integral sign gives
\begin{align}
 \psi'(t)&=\int_{\Rn}\ell\,d\gamma_0(t),\label{eq:app_psi_prime}\\
 \psi''(t)&=\int_{\Rn}\bigl(\ell-\psi'(t)\bigr)^2\,d\gamma_0(t)
 =\operatorname{Var}_{\gamma_0(t)}(\ell).\label{eq:app_psi_second}
\end{align}
On the other hand, differentiating \eqref{eq:app_e0_exponential_arc} yields
\[
 \frac{d\dot\gamma_0(t)}{d\gamma_0(t)}=\ell-\psi'(t).
\]
By the definition of the Fisher--Rao metric and \eqref{eq:app_psi_second},
\[
 \left\lVert\dot\gamma_0(t)\right\rVert_{\gamma_0(t),\mathrm{FR}}^2
 =\int_{\Rn}\left(\frac{d\dot\gamma_0(t)}{d\gamma_0(t)}\right)^2d\gamma_0(t)
 =\psi''(t).
\]
Since $Z(0)=Z(1)=1$, we have $\psi(0)=\psi(1)=0$. Integration by parts and \eqref{eq:app_psi_prime} therefore give
\begin{align*}
 D^{(e_0)}(\mu\Vert\nu)
 &=\int_0^1t\left\lVert\dot\gamma_0(t)\right\rVert_{\gamma_0(t),\mathrm{FR}}^2dt\\
 &=\int_0^1t\psi''(t)\,dt
 =\bigl[t\psi'(t)\bigr]_0^1-\bigl[\psi(t)\bigr]_0^1\\
 &=\psi'(1)
 =\int_{\Rn}\log\frac{d\nu}{d\mu}\,d\nu
 =D_{\rm KL}(\nu\Vert\mu).
\end{align*}
Thus the reverse order is a consequence of following the canonical prescription along the ordered geodesic from $\mu$ to $\nu$. 

\subsection{Otto/\texorpdfstring{$e_1$}{e1}: Reduction to WKL}\label{sec:WKL_divergence_formula}
In this subsection we prove that, when the Otto metric and the $e_1$-geodesic are used, the canonical-divergence functional reduces exactly to the Wasserstein KL divergence formula \eqref{eq:contrast_func}, i.e.,
\begin{align*}
     D^{(\rm e_1)}(\mu \| \nu):=\int_0^1 t \langle\dot{\gamma}_1(t),\dot{\gamma}_1(t) \rangle_{\gamma_1(t)}^{\rm O} dt=\int_{\Rn} \int_0^1 \left(f \circ \varphi_1 - f \circ \varphi_t\right) \, dt \,  d\mu.
\end{align*}
Since $\gamma_1(t)=(\varphi_t)_*\mu$ and the autonomous flow satisfies $\varphi_{t+s}=\varphi_s\circ\varphi_t$, the definition of the weighted divergence gives
\[
\dot{\gamma}_1(t)
=-\divg{\gamma_1(t)}(\grad f)\,\gamma_1(t).
\]
Equivalently, under the Otto isomorphism
\[
\phi_{\gamma_1(t)}(f+\mathbb{R})
=-\Delta_{\gamma_1(t)}(f+\mathbb{R})\,\gamma_1(t),
\]
the tangent vector $\dot{\gamma}_1(t)$ is represented by $f+\mathbb{R}$. Therefore,
\[
\left\langle\dot{\gamma}_1(t),\dot{\gamma}_1(t)\right\rangle_{\gamma_1(t)}^{\rm O}
=\int_{\Rn}\lVert\grad f(y)\rVert^2\,d\gamma_1(t)(y).
\]
Using the pushforward identity,
\begin{align*}
\left\langle\dot{\gamma}_1(t),\dot{\gamma}_1(t)\right\rangle_{\gamma_1(t)}^{\rm O}
&=\int_{\Rn}\left\lVert\grad f(\varphi_t(x))\right\rVert^2\,d\mu(x)\\
&=\int_{\Rn}\left\langle\grad f(\varphi_t(x)),\frac{d}{dt}\varphi_t(x)\right\rangle\,d\mu(x)\\
&=\int_{\Rn}\frac{d}{dt}f(\varphi_t(x))\,d\mu(x).
\end{align*}
Consequently, Fubini's theorem and integration by parts in $t$ yield
\begin{align}\label{eq:KL_energy_formula_specialize_to_WKL}
D^{(\rm e_1)}(\mu\Vert\nu)
&=\int_{\Rn}\int_0^1 t\,\frac{d}{dt}f(\varphi_t(x))\,dt\,d\mu(x)\nonumber\\
&=\int_{\Rn}\left(f(\varphi_1(x))-\int_0^1f(\varphi_t(x))\,dt\right)d\mu(x)\nonumber\\
&=\int_{\Rn}\int_0^1\left(f\circ\varphi_1-f\circ\varphi_t\right)(x)\,dt\,d\mu(x).
\end{align}
Because $\varphi_t$ is affine, $f\circ\varphi_t$ and its time derivative are quadratic in $x$, uniformly for $t\in[0,1]$; hence Gaussian moment bounds justify all interchanges above. This completes the reduction to the WKL-divergence formula \eqref{eq:contrast_func}.
\setlength{\bibsep}{0pt}
\bibliography{mybibliography}

\end{document}